\documentclass[reqno,11pt]{amsart}
\usepackage{amscd,amssymb}
\usepackage{rotating}
\usepackage{lscape}

\theoremstyle{plain}
\newtheorem{Thm}[subsection]{Theorem}
\newtheorem{Cor}[subsection]{Corollary}
\newtheorem{Lem}[subsection]{Lemma}
\newtheorem{Prop}[subsection]{Proposition}
\newtheorem{Conj}[subsection]{Conjecture}

\theoremstyle{definition}
\newtheorem{Def}[subsection]{Definition}

\theoremstyle{remark}

\newtheorem{Rem}[subsection]{Remark}

\newtheorem{conj}{Conjecture}

\numberwithin{equation}{section}
\renewcommand{\rm}{\normalshape}

\newif\ifShowLabels
\ShowLabelstrue
\newdimen\theight
\def\TeXref#1{%
    \leavevmode\vadjust{\setbox0=\hbox{{\tt
        \quad\quad  {\small \rm #1}}}%
    \theight=\ht0
    \advance\theight by \lineskip
    \kern -\theight \vbox to
    \theight{\rightline{\rlap{\box0}}%
    \vss}%
    }}%

\ShowLabelsfalse

\newenvironment{thm}[1]%
    { \begin{Thm} \label{T:#1}  \ifShowLabels \TeXref{T:#1} \fi }%
    { \end{Thm} }

\renewcommand{\th}[1]{\begin{thm}{#1} \sl }
\renewcommand{\eth}{\end{thm} }

\newenvironment{lemma}[1]%
    { \begin{Lem} \label{L:#1}  \ifShowLabels \TeXref{L:#1} \fi }%
    { \end{Lem} }
\newcommand{\lem}[1]{\begin{lemma}{#1} \sl}
\newcommand{\elem}{\end{lemma}}

\newenvironment{propos}[1]%
    { \begin{Prop} \label{P:#1}  \ifShowLabels \TeXref{P:#1} \fi }%
    { \end{Prop} }
\newcommand{\prop}[1]{\begin{propos}{#1}\sl }
\newcommand{\eprop}{\end{propos}}

\newenvironment{corol}[1]%
    { \begin{Cor} \label{C:#1}  \ifShowLabels \TeXref{C:#1} \fi }%
    { \end{Cor} }
\newcommand{\cor}[1]{\begin{corol}{#1} \sl }
\newcommand{\ecor}{\end{corol}}

\newenvironment{defeni}[1]%
    { \begin{Def} \label{D:#1}  \ifShowLabels \TeXref{D:#1} \fi }%
    { \end{Def} }
\newcommand{\defe}[1]{\begin{defeni}{#1} \sl }
\newcommand{\edefe}{\end{defeni}}

\newenvironment{remark}[1]%
    { \begin{Rem} \label{R:#1}  \ifShowLabels \TeXref{R:#1} \fi }%
    { \end{Rem} }
\newcommand{\rem}[1]{\begin{remark}{#1}}
\newcommand{\erem}{\end{remark}}

\newenvironment{conjec}[1]%
    { \begin{Conj} \label{Co:#1}  \ifShowLabels \TeXref{Co:#1} \fi }%
    { \end{Conj} }
\renewcommand{\conj}[1]{\begin{conjec}{#1} \sl }
\newcommand{\econj}{\end{conjec}}

\newcommand{\eq}[1]%
    { \ifShowLabels \TeXref{E:#1} \fi
       \begin{equation} \label{E:#1} }
\newcommand{\eeq}{ \end{equation} }

\newcommand{\prf}{ \begin{proof} }
\newcommand{\epr}{ \end{proof} }

\newcommand\lam{\lambda}        \newcommand\Lam{\Lambda}

\newcommand\calW{{\mathcal{W}}}

 \newcommand\grg{{\mathfrak{g}}}

\newcommand\sdp{\times \hskip -0.3em {\raise 0.3ex
\hbox{$\scriptscriptstyle |$}}} 

\newcommand\Aut{\operatorname{Aut}}

\newcommand\const{\operatorname{const}}

\newcommand\End{\operatorname{End\,}}

\newcommand\Gr{\operatorname{Gr}}

\newcommand\Hom{\operatorname {Hom}}

\newcommand\id{\operatorname{id}}
\newcommand\Id{\operatorname{Id}}

\newcommand\Int{\operatorname{Int}}

\newcommand\Sym{\operatorname{Sym}}

\newcommand\Tr{\operatorname{Tr}}

\renewcommand{\Id}{\text{Id}}

\newcommand\nc{\newcommand}
\nc{\opn}{\operatorname}

\nc\aff{\operatorname{aff}}
\nc\oGr{\overline{\Gr}}
\nc\Bun{\operatorname{Bun}}
\nc\hgrg{\widehat{\grg}}
\renewcommand\Int{\operatorname{Int}}
\nc\bInt{\overline{\Int}}
\nc\hatLam{\widehat{\Lam}}
\nc\bmu{\overline{\mu}}
\nc\bnu{\overline{\nu}}
\nc\blambda{\overline{\lam}}

\nc\ocalW{\overline{\calW}}
\nc\pos{\operatorname{pos}}
\nc\IH{\operatorname{IH}}
\nc\fsl{\mathfrak{sl}}
\nc\fgl{\mathfrak{gl}}
\nc\Rep{\operatorname{Rep}}
\nc\Gal{\operatorname{Gal}}
\nc{\tilGr}{\widetilde{\Gr}}

\nc\Pic{\operatorname{Pic}}
\nc\hgl{\widehat{\fgl}}
\nc\hsl{\widehat{\fsl}}
\begin{document}
\title{On the size and structure of finite linear groups}
\author{Boris Weisfeiler}
\maketitle



{\em This is a nearly complete, previously unpublished 
manuscript\footnote{Parts of  the manuscript were circulated among the experts. Its results were discussed by
Walter Feit [Weisfeiler's work on finite linear groups, Abstracts of Papers Presented to the AMS, 1994, 897-20-299]  and
 Michael Collins [On Jordan's theorem for complex linear groups, J. Group Theory 10 (2007), 411-423; Bounds for finite primitive complex linear groups, J. Algebra 319 (2008), 759-776; Modular analogues of Jordan's theorem for finite linear groups, arXiv:0709.3245, J. Reine Angew. Math.]. ).} by Boris Weisfeiler. 
The results were announced by him in August 1984. Soon after, in early January 1985, he disappeared during a hiking trip in Chile.

 The investigation into Boris Weisfeiler disappearance is still ongoing in Chile, see}

 https://www.weisfeiler.com/boris

{\em 
I would like to thank Tamara Kurdyaeva for typesetting the document in 
 LaTex which made the arXiv posting possible,  and Roman Bezrukavnikov and Michael Finkelberg for organizing the process. 

\medskip  

\ \ \ \ Olga Weisfeiler}

\bigskip

\section{\textbf{Introduction.}}

Our object here is to study linear, and, except in Section 14, finite groups. Our results concern mostly the size of such groups although some other, structural, results are obtained as well. The departure point for the present work was a result of M. Nori \cite{21} which can be considered as a conceptual refinement of Jordan's theorem of linear groups. It turned out that the methods used in \cite{20} and \cite{28} (and based on classification of finite simple groups) can be used to generalize, extend, and strengthen both Nori's [Theorem B] and Jordan's theorem. Most of the present work is dedicated to obtaining the best bounds for our version of Jordan's theorem. This turned out to be quite difficult, especially because of necessity to specially handle groups in small dimensions. A qualitative result is much easier to obtain, see B. Weisfeiler [NAS].

Before going on to the statements of our results let us introduce some terminology. For a field $k$ we denote by $p(k)$ the characteristic exponent of $k$, that is $p(k)=$ char $k$ if char $ k > 0$ and $p(k)=1$ if char $ k = 0$. A $l$-group and a group of Lie $l$-type are both trivial. If $p\neq 1$ is a prime then a group of Lie $p$-type is a group of Lie type of characteristic $p$ (see Section 4 for more detail). A group is said here to be centrally simple if its quotient by the center is simple. Two groups are centrally isomorphic if their quotients by the centers are isomorphic. $O_p(G)$ is, as usual, the largest normal $p$-subgroup of $G$.

Let $k$ denote an algebraically closed field of characteristic exponent $p$. Our version of Jordan's theorem is \\

\th{}(see (13.1)). Let $G$ be a finite subgroup of $GL_n(k)$. Then $G$ contains\\

i) a normal subgroup $T \supseteq O_p(G)$,\\

ii) a normal subgroup $L \supseteq T$\\

such that\\

(a) $T/O_p(G)$ is a commutative $p'$-group isomorphic to a product of $\leq n$ cyclic groups,\\

(b) $L/T$ is isomorphic to a direct product of finite simple groups of Lie $p$-type,\\

(c) $\left|G/LT\right| \leq $
$\left[ 
\begin{aligned}
& n^4\,(n+2)! \ \text {if $n\leq 63$} \\ 
& (n+2)! \ \text {if $n > 63$}
\end{aligned}
\right.$	\\
																 			
\eth

If $p=1$ then $O_p(G)=1$ and $L/T=1$ and we obtain the usual statement of Jordan's theorem:\\

\textbf{Theorem. } \textit{If $p=1$ then $G$ contains a normal commutative subgroup $B$ that $\left|G/B\right|\leq f(n)$.}\\

The best $f(n)$ known until our paper (B. Weisfeiler \cite{?}) was obtained by G. Frobenius (see A. Speiser [\cite{23}, Satz 201]); it was $f(n)=n!\,n\,12^{n(\pi\,(n+1)+1)}$ where $\pi(n)$ is the number of primes $\leq n$. Recall that $\pi \sim n/\ln{n}$. Thus our estimate is of the type $n^{\const\cdot n}$ and G. Frobenius' is of the type $n^{\const\cdot n^2/(\ln{n})^2}$. 

When $p\neq 1$ our result implies that of R. Brauer and W. Feit:\\

\textbf{Theorem} \textit{(see W. Feit [book, Theorem XI.1.2]). If $p^m$ is the order of the Sylow $p$-subgroup of $G$ then $G$ contains a normal commutative $p'$-subgroup $B$ such that $\left|G/B\right| \leq f(p,m,n)$ for an appropriate function of three variables.}\\

Our Theorem 1.1 gives the above theorem with $f(p,m,n)=p^{3\,m}\,n^4\,(n+2)!$ (see (13.2)). This, of course, improves the estimate of R. Brauer and W. Feit. But it seems also noteworthy that our function $f(p,m,n)$ shows once again that the deviation of characteristic $> 0$ case from characteristic $0$ case is concentrated in the $p$-subgroup. Thus our $f(p,m,n)$ is a product of the cube of the order of the Sylow $p$-subgroup with a function independent of $p$ and $m$. 

Following H. Bass [J. Alg.] and using results of J. Tits [Free \cite{27}] and very recent classification of periodic simple linear groups (see, e.g., S. Thomas [vol. 41]) we can obtain the following structure result\\

\th{} (see (14.1)). Let $G$ be a subgroup of $GL_n(k)$. Then $G$ contains\\

1) a triangulizable normal subgroup $T$,\\

2) a normal subgroup $P\supseteq T$,\\

3) a normal subgroup $F\supseteq T$,\\

4) a normal subgroup $L\supseteq T$\\

such that\\

(a) the Zariski closures of $P/T$ and $F/T$ are connected and semi-simple,\\

(b) $P/T$ is simple periodic of Lie $p$-type,\\

(c) $[P,F]\subseteq T$ and $F$ has a certain minimality property,\\

(d) $L/T$ is simple finite of Lie $p$-type,\\

(e) $\left|G/PFL\right| \leq n^4\,(n+2)!$\\

\eth

An interesting feature of this result (except for the estimate) is that it exhibits a decomposition of linear groups into $P$ and $F$ parts.\\

A version of (1.1) for primitive groups is more precise:\\

\th{} (see (11.1)). Let $G$ be a primitive subgroup of $GL_n(k)$ with center $C$. Then $G$ contains\\

i) a normal subgroup $A$ isomorphic to a direct product of alternating groups $Alt_{m_i}$, $m_i \geq 10$,\\

ii) a normal perfect subgroup $L$ centrally isomorphic to a direct product of finite simple groups of Lie $p$-type\\

such that 
$$
\left|G/ACL\right| \leq n^{2\,\log_2{n}+5}.
$$
\ \\
\eth

In other words this theorem says that unlimited growth of $\left|G/C\right|$ comes from groups of Lie $p$-type and the growth of type $n^{c\cdot n}$ comes from the alternating groups. The order of the remaining part is only of the type $n^{c\cdot\ln{n}}$, i.e., incomparably smaller. This result should be, perhaps, compared with a result of P. Cameron [\cite{7} , Theorem 6.1], where one sees an estimate $n^{c\cdot\ln{\ln{n}}}$ on the order of a primitive group, other than some specified groups.\\

Our proof of (1.3) begins with a study of centrally simple linear groups.\\

\prop{} (see (?)). Let $G$ be a finite centrally simple non-commutative subgroup of $GL_n(k)$. Suppose that $G$ is not of Lie $p$-type and not isomorphic to an alternating group. Then

$$
\left| \Aut{G}\right| \leq n^4\,(2\,n+1)^{2\,\log_3{(2\,n+1)}+1}.
$$
\ \\
\eprop

This result has relevance to the study of maximal subgroups of finite groups of Lie $p$-type, $p > 1$ or, the same, of the primitive permutation representations of the latter. In particular, one can combine (1.4) with the results of M. Aschbacher \cite{1} and M. Liebeck \cite{19}. (I am grateful to M. Liebeck for making his paper available to me before publication):\\

\th{} Let $H_0 $ be a classical simple group of Lie $p$-type $^{c}X_a(q^c)$ and $H$ a subgroup of $\mathrm{Aut}\ H_0$ with $H\supseteq H_0$. Let $G$ be a maximal subgroup of $H$. Then one of the following holds:\\

(a) $G$ is "known" (a list, called $C_H$, oh these $G$ is given in M. Aschbacher \cite{1}),\\

(b) the socle of $G$ is simple of Lie $p$-type and $\left|G\right|\leq q^{3\,c\,n}$, where $n$ is the dimension of the natural representation of $^{c}X_a$,\\

(c) the socle of $G$ is an alternating group,\\

(d) the socle of $G$ is simple and $\left|G\right| \leq n^4\,(2\,n+1)^{2\,\log_3{(2\,n+1)}+1}$.\\

\eth

Note that the estimate in (b) has the type $q^{\const\cdot n}=n^{\const\cdot n/{\ln{n}}}$ (with the latter constant increasing with $q$). Thus again there is a wide gap between the estimates in (b) and (d). Of course, the type of the estimate in (b) can not be substantially improved. Thus it seems desirable to separate cases (b), (c), and (d).\\

Another implication for maximal subgroups gives the following, rougher\\

\prop{} Let $H_0 $ be a simple group of Lie $p$-type $^{c}X_a(q^c)$ and $H$ a subgroup of $\mathrm{Aut}\ H_0$ with $H\supseteq H_0$. Let $G_0$ be a perfect centrally simple group and $G$ a subgroup of $\mathrm{Aut}\ G_0$ with $G\supseteq G_0$. If $G$ is a maximal subgroup of $H$ then one of the following holds\\

\ \ \ \ (a) $G_0$ is of Lie $p$-type,\\

or (b) $G_0$ is an alternating group,\\

or (c) $\left|G\right|\leq r^4\,(2\,r+1)^{2\,\log_3{(2\,r+1)}+1}$\\

\begin{tabular}{cccccccccc}
where & $r =$ & $n+1$ & $2\,n$ & $2\,n+1$ & 7 & 26 & 27 & 56 & 248 \\
if & $X_a=$ & $A_n$ & $C_n,D_n$ & $B_n$ & $G_2$ & $F_4$ & $E_G$ & $E_7$ & $E_8$ \\
\end{tabular}
\\
\eprop

This follows directly from (1.4) applied to the composite of $G_0\rightarrow$ $^{c}X_a(\mathbb{\bar{F}}_p) \rightarrow GL_r(\mathbb{\bar{F}}_p)$. Since, when $q$ varies, the tower of groups $^{c}X_a(q^c)$ is infinite, it follows that groups in (b) and (c) can be maximal only for finitely many $q$:\\

\th{} (see ?). In the assumptions of (1.6) there exists $r$, depending on $^{c}X_a$ (and not on $p$), such that if $q > p^r$ and $G$ (as in (1.6)) is maximal in $H$ then $G_0$ is of Lie $p$-type.\\
\eth

We give in (?) explicit values for $r$ when $H$ is of classical type.

\newpage
\section{\textbf{Notation and preliminaries.}}

\subsection{}Some of terminology and notation was introduced in Section 1.\\

\subsection{}If $M$ is a group and $X$ a subset of $M$ then $Z_M(X)$ (resp. $N_M(X), C(M), \mathrm{Aut}\ M,\mathrm{Aut}_c\ M, \mathrm{Out}\ M$) denotes the centralizer of $X$ in $M$ (resp. the normalizer of $X$ in $M$, the center $M$, the group of automorphisms of $M$, the group of automorphisms of $M$ trivial on $C(M)$, the group of outer automorphisms of $M$). Occasionally we also write $NZ_M(X)$ for $N_M(X)/Z_M(X)$ and, in the case when $X$ is a group, $\overline{NZ}_M(x)$ for $N_M(X)/X\cdot Z_M(X)$ . \\

\subsection{}The symmetric and alternating groups on $n$ letters (resp. on a set $X$) are denoted $\Sym_n$ and Alt$_n$ (resp. $\Sym{X}$ and Alt $X$).\\

\subsection{}$\mathbb{N}$ is the set of non-negative integers. $\mathbb{Z}/a$ is the cyclic group of order $a$.\\

\subsection{}$\log{x}$ (resp. $\ln{x}$) denotes $\log_2{x}$ (resp. the natural logarithm of $x$). $\Gamma(x)$ denotes the $\Gamma$-function of $x$, so that $\Gamma(n+1)=n!$ when $n\in\mathbb{N}$. We often use notation $f(x)$ for 
$$
f(x)=(2\,x+1)^{2\,\log_3(2\,x+1)+1},
$$ 
it is one of our main functions.\\

\subsection{}Our notation for the parameter of twisted groups of Lie type agrees with that of R. Steinberg \cite{?} and, therefore, differs from that of other authors, see (4.1.1) below.\\

\subsection{}In our study we will need repeatedly the precise knowledge of centrally simple groups having faithful linear representations of the given degree. These are listed below in Table T2.7. In this table $a\cdot G$ denotes a perfect central extension of $G$ by $\mathbb{Z}/a$; however, if $a\cdot G$ appears in characteristic $p$ with $p|a$ then it should be read as $(a/p)\cdot G; a\cdot G (p=l)$ means that this group appears only in characteristic $l$. This Table is complied from W. Feit [Nice \cite{11}, \S 8.4] and A. Zalessky [\cite{33} 1981, \S 13]. See Table T6.3 for different isomorphisms.

\newpage

\begin{center}
Table T2.7.\\
\textit{Centrally simple linear groups of small degree.}\\
\end{center}
\ \\
\begin{center}
\begin{tabular}{|c|c|c|c|}
\hline
$n$ & \multicolumn{3}{c|}{$G \leq GL_n(k), p(k)=p$} \\
& (almost) any $p$ & sporadic $p$ & Lie $p$-type \\
\hline

& \ & \ & \ \\ 

$2$ & $2\cdot$Alt$_5$ & \ & $A_1(p^a)$  \\
& \ & \ & \ \\ 

\hline
& \ & \ & \ \\ 

$3$ & Alt$_5$, $6\cdot$Alt$_6$, $\bar{A}_1(7)$ & $3\cdot$Alt$_7(p=5)$ & $A_2(p^a), \bar{A}_1(p^a),$ $^{2}A_2(p^{2\,a})$ \\
& \ & \ & \ \\ 

\hline
& \ & \ & \ \\ 

$4$ & Alt$_5(p\neq 5)$, $6\cdot$Alt$_6$,  $2\cdot$Alt$_7$, & Alt$_6(p=2)$, Alt$_7(p=2)$, &
$A_3(p^a),$ $^{2}A_3(p^{2\,a}), B_2(p^a),$ \\

& $2\cdot$Alt$_5$, $2\cdot\bar{A}_1(7)$, $2\cdot\bar{B}_2(3)$ & $4\cdot\bar{A}_2(4)(p=3)$ &
$^{2}B_2(2^{2\,a+1})(p=2), $ \\

& \ & \  & $A_1(p^a)(p > 3), \bar{A}_1(p^a)(a\geq 2)$ \\
& \ & \ & \ \\ 

\hline
& \ & \ & \ \\ 

$5$ & Alt$_5(p\neq 2)$, Alt$_6(p\neq 2, 3)$, & Alt$_7(p=7)$, $M_{11}(p=3)$ &
$A_4(p^a),$ $ ^{2}A_4(p^{2\,a}), \bar{B}_2(p^a),$ \\

& $\bar{A}_1(11), \bar{B}_2(3)$ & \  & $\bar{A}_1(p^a)(p \geq 5)$ \\

& \ & \ & \ \\
\hline
\end{tabular}
\end{center}

\newpage
\section{\textbf{Estimates for the alternating groups.}}
Let $k$ be a field and $p=p(k)$ its characteristic exponent. We quote here some results of I. Schur, L. E. Dickson, and A. Wagner.\\

\prop{} Let $H\simeq$ Alt$_m$. Let $\varphi: H\rightarrow GL_n(k)$ be faithful irreducible representation.\\

i) If $p=1$ and $m\geq 4, m\neq 5$, then $n\geq m-1$; for $m=5, n\geq 3$,\\

ii) if $m\geq 9$ or $p\neq 2$ and $m\geq 7$ then $n\geq m-2$; moreover, if $p\nmid m$ then $n\geq m-1$,\\

iii) if $m=5$ (resp. 6, 7, 8) then $n\geq 2$ (resp. 3, 4, 4).\\

\eprop

\prf{} (3.1) (i) is a result of I. Schur [\cite{22}, \S 44]; (3.1) (ii) and (iii) are results of A. Wagner [\cite{29}, ?] (although essentially known from L. E. Dickson \cite{9}). 

See G. D. James [\cite{17}, Theorem 6 (ii)] for (ii) when $m\geq 10$. When $m=9$ see A. Wagner [\cite{29}, ?] and when $p\neq2, m\geq 7$, see A. Wagner \cite{29}. 

To see (iii) we note that Alt$_5 \simeq SL_2(\mathbb{F}_4)$ and has therefore a $2$-dimensional representation in characteristic $2$; Alt$_6 \simeq PSL_2(\mathbb{F}_9)$ and has therefore an irreducible $3$-dimensional (=adjoint) representation in characteristic $3$, but it does not have, by Table T2.7 (the list of linear groups of small degree), representations of dimension $2$; Alt$_8\simeq \mathcal{D}Sp_4(\mathbb{F}_2)$ and has therefore a representation of dimension $4$ in characteristic $2$, but, again by Table T2.7 it does not have smaller representations; Alt$_7$ has by above a $4$-dimensional representation in characteristic $2$, but by Table T2.7 it has no smaller representations.\\
\epr

\cor{} Let $H, \varphi,$ and $n$ be as in (3.1). If $n\geq 8$ then $\left|H\right|\leq (n+2)!/2.$\\

\ecor

\prop{} Let $H\simeq$ Alt$_m$. Let $\varphi: H \rightarrow PGL_n(k)$ be a faithful projective irreducible representation. Suppose that $\varphi$ does not lift to a linear representation of $H$.\\

i) If $p=1, m\geq 4, m\neq 6$, then $m\leq 2+2\,\log{n}$; if $p=1, m=6$, then $n \geq 3$,\\

ii) for all $p$ and $m>7$ we have $m \leq (81+32\,\log{n})/15 \leq 5.4+2.134\,\log{n}$,\\

iii) for all $p$ and $4\leq m \leq 16$ we have \\

\begin{tabular}{cccccccccccccc}
m& 4 & 5 & 6 & 7 & 8 & 9 & 10 & 11 & 12 & 13 & 14 & 15 & 16 \\
$n\geq$ & 2 & 2 & 3 & 3 & 8 & 8 & 8 & 16 & 16 & 16 & 32 & 32 & 128 \\
\end{tabular}
\ \\
\eprop

\prf{} (3.3) (i) is a result of I. Schur [\cite{22}, \S 44]. 

To prove (ii) and (iii) write $m = 2^{w_1}+2^{w_2}+\ldots + 2^{w_s}, w_1 > \ldots > w_s$, for the $2$-adic decomposition of $m$. Then by A. Wagner [\cite{29} Theorem 1.3 (ii)] we see that $2^{\frac{m-s-1}{2}}\mid n$ if $m > 7$. 

For $7<m\leq16, m\neq 11$, this gives us the estimate in (iii). 

For $m=5,6,7$ the estimates in (iii) follow from the Table T2.7 (of linear groups of small degree). 

For $m=4$ clearly $n\geq 2$ since Alt$_4$ is not commutative. 

For $m=11$ we get $s=3$ whence $8\mid n$. Suppose $n=8$. Consider in $H\simeq$ Alt$_{11}$ the subgroup $H_1\times H_2\simeq$ Alt$_8\times$ Alt$_3$. We know that $\varphi$ lifts to a linear representation $\tilde{\varphi}: \tilde{\mathrm{Alt}}_{11}\rightarrow GL_8(k)$ where $\tilde{\mathrm{Alt}}_{11}$ is the (non-split) double cover of Alt$_{11}$ (I. Schur [\cite{22}, \S 5, Theorem II]). Let $\pi(n): \tilde{\mathrm{Alt}}_{11}\rightarrow$ Alt$_{11}$ be the covering map. Then the relations (I. Schur [\cite{22}, \S 5, relations (IV)]) show that $\pi: \pi^{-1}(H_1)\rightarrow H_1$ is the non-split double cover of Alt$_8$. Since $n=8$ the representation will be irreducible for $H_1$ (by (3.1)(3) with $m=8$). Since $H_2\simeq \mathbb{Z}/3$ and $\ker{\pi}\simeq \mathbb{Z}/2$ it follows that $\pi^{-1}(H_2)\simeq \mathbb{Z}/3\times \mathbb{Z}/2$. Therefore $\tilde{\varphi}(\pi^{-1}(H_2))$ commutes with $\tilde{\varphi}(\pi^{-1}(H_1))$. Since the latter is irreducible, $\tilde{\varphi}(\pi^{-1}(H_2))\subseteq kId_8$. But then $\tilde{\varphi}(\pi^{-1}(H_2))$ is in the center of $\tilde{\varphi}(\pi^{-1}(H))$, an impossibility. Returning to the general $m$ we see that $s\leq \log{(m+1)}$. Since $\log{(x+1)}$ is a convex function it is bounded from above by a tangent line at any point. Thus $s\leq \log{(m+1)}\leq (m+49)/16$ (where we took the tangent at $x=15$). Thus $n\geq 2^{\frac{m-\log{(m+1)-1}}{2}}$, or $\frac{m-\log{(m+1)}-1}{2}\leq \log{n}$ or $2+2\,\log{n}\geq m-\log{(m+1)}\geq m-(m+49)/16=(15\,m-49)/16$. This gives $m\leq (81+32\,\log{n})/15$ whence (ii).\\
\epr

\cor{} Let $H,\varphi, n, m$ be as in (3.3). If $n\geq 2, m\geq 4$, then $$\left|\Aut H\right|\leq (2\,n+1)^{2\,\log_3{(2\,n+1)+1}}$$.
\ecor

\prf{} Recall that $\Aut$ Alt$_m\simeq \Sym_m$ if $m\geq 4, m\neq 6$, (see B. Huppert \cite{16}) and $\left|\Aut\mathrm{Alt}_6\right|=2\cdot 6!$. Now our claim is verified directly for cases of (3.3)(iii) using Table TA (values of functions for small arguments). If $m>16$ we have by (3.3)(ii) that 
$$
\left|\Aut{H}\right|=m!\leq m\cdot (m/e)^m=e\cdot (m/e)^{(m+1)}\leq 2.72\cdot(2+0.8\,\log{n})^{6.4+2.134\,\log{n}}.\\
$$

So it is sufficient to check that
$$
2.72\cdot(2+0.8\,\log{n})^{6.4+2.134\,\log{n}}< (2\,n+1)^{2\,\log_3{(2\,n+1)+1}}
$$
for $n\geq 128$.\\

Now $\log_3{x}=\log{x}/\log{3}$. Therefore
$$
(2\,n+1)^{2\,\log_3{(2\,n+1)+1}}> (2\,n+1)^{1.26\,\log{(2\,n+1)+1}}> (2\,n)^{1.26\,\log{(2\,n)+1}}=(2\,n)^{1.26\,\log{n}+2.26}=
$$
$$
2^{2.26}\cdot n^{1.26}\cdot n^{1.26\,\log{n}+2.26}> 4.79\,n^{1.26\,\log{n}+3.52}> 4.79\,n^{1.26\,(\log{n}+2.79)}.\\
$$
On the other hand 
$$
2.72\,(2+0.8\,\log{n})^{6.4+2.134\,\log{n}}< 2.72\,(2+0.8\,\log{n})^{2.3\,(\log{n}+2.79)}.\\
$$
Thus it suffices to establish that
$$
2.72\,(2+0.8\,\log{n})^{2.3\,(\log{n}+2.79)}< 4.79\,n^{1.26\,(\log{n}+2.79)}
$$
or 
$$
(2+0.8\,\log{n})^{2.3}< n^{1.26}
$$
or $2+0.8\,\log{n}< n^{0.547}$. This holds for $n=128$. On the other hand if $f(x):= 2+0.8\,\log{x}-x^{0.547}$ then $f'(x)=0.8/x-0.547\,x^{-0.453}< 0$ for $x\geq 128$. Thus $f(x)< 0$ for $x\geq 128$ whence our claim.
\epr

\newpage
\section{\textbf{Recollections and preliminaries about groups of Lie type.}}
We use R. Steinberg \cite{?} as basic reference for groups of Lie type. In particular, we denote by $^{c}X_a(m^c)$ the universal group of Lie type $^{c}X_a(m^c)$. As usual when $c=1$ we just we just write $X_a(m)$. The groups $\mathcal{D}^{c}X_a(m^c)$ are also considered of Lie type (see (4.3.1)(b) below).

\subsection{}Here $m$ is the parameter associated to our group. This $m$ is an integral power $m=q^s$ of a prime $q$ if $^{c}X_a(m^c)\simeq G(\mathbb{F}_m)$ for some simply connected algebraic $\mathbb{F}_m$ - group. For groups of type $^{2}B_2=$ $^{2}C_2$, $^{2}F_4,$ and $^{2}G_2 \ \  m^c$ is an \textit{odd} power of a prime, $m=q^s, 2\,s\in \mathbb{N}, s\notin \mathbb{N}.$ 

\subsubsection{}N.B. Some authors to whom we refer (Gorenstein, Landazuri, and Seitz among them) use $m$ differently. For them $"^{c}X_a(m)"$ is our $^{c}X_a(m^c)$ except when $^{c}X_a=$ $^{2}B_2$, $^{2}F_4$, $^{2}G_2$. For these latter groups their notation $^{2}X_a(m^2)$ coincides with ours but then they write all related expressions (e.g. the order) as functions of $m^2$ (and not of $m$ as we do).

\subsection{}When $m^c$ is a power of a prime $q$ we say that $q$ is a characteristic of $^{c}X_a(m^c)$ or of a perfect group centrally isomorphic to $^{c}X_a(m^c)$ and $\mathcal{D}^{c}X_a(m^c)$ or that these are of Lie $q$-type. Note that $q$ generally depends not on the central isomorphism class of (an abstract group) $^{c}X_s(m^c)$ but on its representation as $^{c}X_a(m^c)$, see (4.3.2) below. We write $q=q(^{c}X_a(m^c)),  q=q(\mathcal{D}^{c}X_a(m^c))$ etc.

\subsection{\,}We denote by $^{c}\bar{X}_a(m^c)$ the central quotient of $^{c}X_a(m^c)$.

\subsubsection{}$^{c}\bar{X}_a(m^c)$ is simple non-commutative except in the following cases\\

(a) $\bar{A}_1(2), \bar{A}_1(3),$ $^{2}\bar{A}_2(4),$ $^{2}\bar{B}_2(2)$ are solvable;\\

(b) the derived group of $B_2(2), G_2(2),$ $^{2}F_4(2),$ $^{2}G_2(3)$ is simple non-commutative of prime index $q$ (where $q$ is as in (4.2)).\\

See R. Steinberg [\cite{?}, Theorems 5 and 34 and comments on them].

\subsubsection{}The central quotients of the following groups are (sporadically) isomorphic to groups of Lie type in different characteristic or to alternating groups:
$$
A_1(4), A_1(5), A_1(7), A_1(8), A_1(9),
$$
$$
A_2(2), A_3(2), B_2(3),\ ^{2}A_2(9),\ ^{2}A_3(4).
$$

The isomorphisms are given in Table T6.3.

See R. Steinberg [\cite{?}, Theorem 37].

When $^{c}\bar{X}_a(m^c)$ is simple non-commutative we denote by $^{c}\tilde{X}_a(m^c)$ the universal cover of $^{c}\bar{X}_a(m^c)$. The kernel of the canonical map $^{c}\tilde{X}_a(m^c)\rightarrow$ $^{c}X_a(m^c)$ is called the Schur multiplier. 

\subsubsection{}(a) $^{c}\tilde{X}_a(m^c)$ is isomorphic to $^{c}X_a(m^c)$ except in the following cases 
$$
A_1(4), A_1(9), A_2(2), A_2(4), A_3(2),
$$
$$
B_2(2), B_3(2), B_3(3), D_4(2), F_4(2),
$$
$$
G_2(2), G_2(4), \ ^{2}A_3(4), \ ^{2}A_3(9),\ 
^{2}A_5(4), \  ^{2}E_6(4).
$$

\ \ \ \ \ \ \ (b) In any case (exceptional or not) the kernel of $^{c}\tilde{X}_a(m^c)\rightarrow$ $^{c}X_a(m^c)$ is a $q$-group where $m^c=q^s$ and $q$ is a prime. It holds for any $q$ for which $^{c}\bar{X}_a(m^c)$ happens to be of Lie $q$-type.

See R. Steinberg [\cite{?}, 1] or D. Gorenstein [\cite{G}, Table 4.1, p.302], where kernel of $^{c}\tilde{X}_a(m^c)\rightarrow$ $^{c}\bar{X}_a(m^c)$ is also explicitly given.

\newpage
\pagestyle{empty}

\begin{center}
Table T4.4.\\
\end{center}
\ \\
\footnotesize
\begin{tabular}{|c|c|c|c|c|c|}
\hline
& \ & \ & \ & \ & upperbound\\
$^{c}X_a$ & $d(X_a)$ & $b$ & $A_g$ & $A_d$ & on\\
& see 4.4.1 & see 4.4.2 & see 4.5 & see 4.5 &  $ c\,\left|A_d\right|\,\left|A_g\right|$ \\
\hline

& \ & \ & \ & \ & \\ 

$A_1$ & $3$ & $1$ & $1$ & $\mathbb{Z}/(2, m-1)$ & $2$ \\

\hline
& \ & \ & \ & \ & \\ 

$A_n, n\geq 2$ & $n^2+2\,n$ & $n$ & $\mathbb{Z}/2$ & $\mathbb{Z}/(n+1, m-1)$ & $2\,(n+1)$ \\

\hline
& \ & \ & \ & \ & \\ 

$B_2, q=2,$ & $10$ & $2.5$ & $\left|A_g\right|=2$ & $1$ & $2$ \\
$m\neq 2$ & \ & \ & \ & \ & \\  

\hline
& \ & \ & \ & \ & \\ 

$B_n, q\neq 2,$ & $2\,n^2+n$ & $2\,n-2$ & $1$ & $\mathbb{Z}/2$ & $2$ \\
$n\geq 2$ & \ & \ & \ & \ & \\ 

\hline
& \ & \ & \ & \ & \\ 

$C_n, q=2,$ & $2\,n^2+n$ & $2\,n-2.3$ & $1$ & $1$ & $1$ \\
$n\geq 3$ & \ & \ & \ & \ & \\ 

\hline
& \ & \ & \ & \ & \\ 

$C_n, q\neq 2,$ & $2\,n^2+n$ & $n$ & $1$ & $\mathbb{Z}/2$ & $2$ \\
$n\geq 3$ & \ & \ & \ & \ & \\ 

\hline
& \ & \ & \ & \ & \\ 

$D_4$ & $28$ & $5$ & $\Sym_3$ & $\left(\mathbb{Z}/(2, m-1)\right)^2$ & $24$ \\

\hline 

$D_n, n\geq 5$ & $2\,n^2-n$ & $2\,n-3$ & $\mathbb{Z}/2$ & 
$\begin{cases}
\left(\mathbb{Z}/(2, m-1)\right)^2,& n\ \text{even} \\
\mathbb{Z}/(4, m^n-1), & n\ \text{odd} \\
\end{cases}$ & 
$8$\\

\hline
& \ & \ & \ & \ & \\ 

$E_6$ & $78$ & $11$ & $\mathbb{Z}/2$ & $\mathbb{Z}/(3, m-1)$ & $6$ \\

\hline
& \ & \ & \ & \ & \\ 

$E_7$ & $133$ & $17$ & $1$ & $\mathbb{Z}/(2, m-1)$ & $2$ \\

\hline
& \ & \ & \ & \ & \\ 

$E_8$ & $248$ & $29$ & $1$ & $1$ & $1$ \\

\hline

$F_4$ & $52$ & $10$ & 
$\begin{cases}
\left|A_g\right|=2 & \text{if $q=2$} \\
1 & \text{if $q\neq 2$} \\
\end{cases}$ & 
$1$ & $2$ \\

\hline 

$G_2$ & $14$ & $3$ & 
$\begin{cases}
\left|A_g\right|=2 & \text{if $q=3$} \\
1 & \text{if $q\neq 3$} \\
\end{cases}$ & 
$1$ & $2$ \\

\hline
& \ & \ & \ & \ & \\ 

$^{2}A_n, n\geq 2$ & $n^2+2\,n$ & $n$ & $1$ & $\mathbb{Z}/(n+1, m+1)$ & $2\,(n+1)$ \\

\hline
& \ & \ & \ & \ & \\ 

$^{3}D_4$ & $28$ & $5$ & $1$ & $1$ & $3$ \\

\hline 

$^{2}D_n$ & $2\,n^2-n$ & $2\,n-3$ & $1$ & 
$\begin{cases}
\left(\mathbb{Z}/(2, m+1)\right)^2,& n\ \text{even} \\
\mathbb{Z}/(4, m^n+1), & n\ \text{odd} \\
\end{cases}$ & 
$8$ \\

\hline
& \ & \ & \ & \ & \\ 

$^{2}E_6$ & $78$ & $15$ & $1$ & $\mathbb{Z}/(3, m+1)$ & $6$ \\

\hline
& \ & \ & \ & \ & \\ 

$^{2}B_2$ & $10$ & $3$ & $1$ & $1$ & $2$ \\

\hline
& \ & \ & \ & \ & \\ 

$^{2}F_4$ & $52$ & $10$ & $1$ & $1$ & $2$ \\

\hline
& \ & \ & \ & \ & \\ 

$^{2}G_2$ & $14$ & $4$ & $1$ & $1$ & $2$ \\

\hline
\end{tabular}

\newpage
\normalsize
\pagestyle{plain} 
\subsection{}For a group $^{c}X_a(m^c)$ of Lie type we denote by $d=d(X_a)$ the dimension of the corresponding algebraic group of type $X_a$. The numbers $d=d(X_a)$ are listed in Table T4.4.

\subsubsection{}
$$
\left|^{c}X_a(m^c)\right|\leq m^d.\\
$$

This is known (and can be easily checked by looking at a table of the orders e.g. in D. Gorenstein [\cite{13}, Table 2.4, p. 135]). 

For groups \textit{not} listed in (4.3.2) we denote by $l=l(^{c}X_a(m^c))$ the smallest degree of centrally faithful irreducible representations of $^{c}\tilde{X}_a(m^c)$ (or, the same, faithful irreducible projective representations of $^{c}\bar{X}_a(m^c)$) over all fields of charecteristic different from $q=q(^{c}X_a(m^c))$.

\subsubsection{}Except for cases listed below $l(^{c}X_a(m^c))\geq (m^b-1)/2$ where $b=b(^{c}X_a(m^c))$ is the number of given in $3^{\underline{d}}$ column of Table T4.4. 

Exceptions: $A_1(4), A_1(9), A_2(4), B_2(2), B_3(3), D_4(2), F_4(2),$ $^{2}A_3(9),$ $^{2}B_2(8),$ $^{2}E_6(4)$.

This can be readily deduced from V. Landazuri and G. Seitz [\cite{18}, p. 419]. The estimates we give are generally worse than the ones given there. The advantage (for us) of our form for degrees is that they are given by a uniform expression.

Using Table T4.4 one can verify now that (with $d$ and $b$ as in (4.4.1) and (4.4.2))

\subsubsection{}(a) $d\leq 2\,b^2+b$, \\

\ \ \ \ \ \ \ (b) $d\leq b^2+2\,b$ if either $^{c}X_a$ is different from $C_n(C_2=B_2, C_1=A_1)$ or $m^c$ is even.

\subsection{}Suppose that $^{c}X_a(m^c)$ is not listed in (4.3.1)(a). Let $A=A(^{c}X_a(m^c))$ denote the group Out$(^{c}X_a(m^c))$ of outer automorphisms of $^{c}\bar{X}_a(m^c)$ i.e. $A:=($Aut$(^{c}X_a(m^c)))/{^{c}\bar{X}_a(m^c)}$. By R. Steinberg \cite{?} (see D. Gorenstein and R. Lyons [\cite{14}, 7] for explicit information) we know that $A$ contains two subgroups (possibly trivial): $A_d$ and $A_f$, and a subset $A_g$. 

$A_g$ is the set of graph automorphisms of $^{c}X_a$ (see R. Steinberg [\cite{?}, Corollary to Theorem 29, Theorem 36, and subsequent remarks to both]); $A_q$ is given in column $4$ of Table T4.4; it is a group unless $^{c}X_a$ is $B_2, F_4$ or $G_2$ and characteristic is $2.2,$ or $3$ respectively in which case $A_g\,A_f$ is a cyclic group generated by the non-trivial element of $A_g$; 

$A_d$ the group of diagonal automorphisms (see R. Steinberg [\cite{?}, Lemma 58 and proof of Theorem 36]); $A_d$ is given in column $5$ of Table T4.4;

$A_f$ the group of field automorphisms, (see R. Steinberg [\cite{?}, just above Theorem 30]).

\subsubsection{}(a) $A_f$ is isomorphic to the Galois group of $\mathbb{F}_{m^c}$ over its prime field $\mathbb{F}_q$,\\ 

\ \ \ \ \ \ \ (b) $A_f\simeq \mathbb{Z}/s$ where $s=c\cdot \log_q{m}$.

\prf
(a) is the definition. Writing $m^c=q^s$ we get (b).
\epr

\subsubsection{}(a) $A=A_d\,A_f\,A_g$,\\

\ \ \ \ \ \ (b) $A_d$ is normal in $A$,\\

\ \ \ \ \ \ (c) $A_g$ can be taken to commute with $A_f$,\\

\ \ \ \ \ \ (d) $\mathcal{D}^3 A=\{1\},$ and $\mathcal{D}^2 A=\{1\}$ unless $^{c}X_a=D_4.$

\prf
(a) is known from R. Steinberg [\cite{?}, Theorems 30 and 36]. (b) is evident from definitions (since $A_d$ can always be chosen to come from a maximal twist-invariant torus, see J. Tits [\cite{27}]). Since $A_g\neq 1$ implies that $c=1$ we can assume in the proof of (c) that $X_a(m)=G(\mathbb{F}_m)$ where $G$ is defined over the prime field $\mathbb{F}_q$. Then $A_g\subseteq (\End {G})(\mathbb{F}_q)/($Inn $G)(\mathbb{F}_q)$ whence (c) evidently follows. Now (d) follows from (c) if one inspects columns $4$ and $5$ of Table T4.4 
\epr

\subsubsection{} 
$$
\left|^{c}X_a(m^c)\right|=\left|^{c}\bar{X}_a(m^c)\right|\cdot\left|A_d\right|.\\
$$

See R. Steinberg [\cite{?}, Exercise (b) in the end of \S 10 and Corollary to Theorem 35] or D. Gorenstein and R. Lyons [\cite{14}, (7-1)(g)].

\subsubsection{}In the notation of (4.4.2) we have except for groups from (4.3.1)(a):\\

(a) $\left|A\right|\leq 
\begin{cases}
4.8\,\log{((m^b-1)/2)}& \text{if $^{c}X_a=D_4$ and $q=2$}  \\
3.03\,\log{((m^b-1)/2)} & \text{if $^{c}X_a\neq D_4$ or $q\neq 2$ } \\
2 & \text{if $m^b\leq 8$}
\end{cases}$ 
\\
\ \\

(b) $\left|A_g\,A_f\right|\leq 
\begin{cases}
1.2\,\log{((m^b-1)/2)} \\
2 & \text{if $m^b\leq 8$}
\end{cases}$

\prf{}
Set $x:=(m^b-1)/2.$ Comparing columns $3$ and $6$ of Table T4.4 we see that 
$$
c\,\left|A_d\right|\cdot\left|A_g\right|\leq 
\begin{cases}
4.8\,b & \text{if $^{c}X_a=D_4$} \\
3\,b & \text{otherwise}
\end{cases}
$$

For $^{c}X_a=D_4$ we have $c\,\left|A_d\right|\,\left|A_g\right|=4.8\,b=4.8\,\log_m{(2\,x+1)}.$ Therefore using (4.5.1)(b) we have 
$$
\left|A\right|=\left|A_f\right|\,\left|A_d\right|\,\left|A_g\right|=c\,\log_q{m}\cdot \left|A_d\right|\,\left|A_g\right|=4.8\,(\log_q{m})\,\log_m({2\,x+1)}=4.8\,\log_q{(2\,x+1)}\leq
$$
$$
\leq
\begin{cases}
4.8\,\log{(2\,x+1)} & \text{if $q=2$} \\
4.8\,\log_3{(2\,x+1)} & \text{if $q\geq 3$}
\end{cases}
$$

Since $4.8\,\log_3{(2\,x+1)}=(4.8/\log 3)\,\log{(2\,x+1)}< 3.03\,\log{(2\,x+1)}$ we have (for $^{c}X_a=D_4$)
$$
\left|A\right|\leq 
\begin{cases}
4.8\,\log{(2\,x+1)} & \text{if $q=2$} \\
3.03\,\log{(2\,x+1)} & \text{if $q\geq 3$}
\end{cases}
$$

If $^{c}X_a\neq D_4$ then $c\,\left|A_d\right|\,\left|A_g\right|\leq 3\,b$ whence $c\,\left|A_d\right|\,\left|A_g\right|\leq 3\,\log_m{(2\,x+1)}$ and 
$$
\left|A\right|=\left|A_f\right|\,\left|A_d\right|\,\left|A_g\right|=c\,\log_q{m}\cdot \left|A_d\right|\,\left|A_g\right|\leq 3\,(\log_q{m})\,\log_m{(2\,x+1)}=
$$
$$
=3\,\log_q{(2\,x+1)}\leq 3\,\log{(2\,x+1)}.
$$

This establishes the first two lines of (a).

The cases $m^b\leq 8$ for exception of those listed in (4.3.1)(a) are treated in Table T4.5.4 below which directly follows from Columns 4, 5, and 3 of Table T4.4. The remaining part of (a) follows from Table T4.5.4.

\begin{center}
Table T4.5.4\\
\ \\
\end{center}

\begin{center}
\begin{tabular}{|c|c|c|c|c|c|c|c|c|}
\hline
& \ & \ & \ & \ & \ & \ & \ & \ \\ 
$G$ & $A_1(4)$ & $A_1(5)$ & $A_1(7)$ & $A_2(2)$ & $A_2(2)$ & $A_3(2)$ & $G_2(2)$ & $^{2}A_3(4)$ \\
\hline

& \ & \ & \ & \ & \ & \ & \ & \ \\ 
$A$ & $A_f$ & $A_d$ & $A_d$ & $A_g$ & $B_g$ & $A_g$ & $\{1\}$ & $A_f$ \\  
\hline 

& \ & \ & \ & \ & \ & \ & \ & \ \\ 
$\left|A\right|$ & $2$ & $2$ & $2$ & $2$ & $2$ & $2$ & $1$ & $2$ \\  
\hline 

& \ & \ & \ & \ & \ & \ & \ & \ \\ 
$x$ & $1.5$ & $2$ & $3$ & $1.5$ & $1.5$ & $3.5$ & $3.5$ & $3.5$ \\ 

\hline
\end{tabular}
\end{center}

\ \\

\ \\

To prove (b) we proceed similarly.  We have $c\,\left|A_g\right|\leq 1.2\,b= 1.2\,\log_m{(2\,x+1)}$ by comparing columns 3 and 4 of Table T4.4. Then
$$
\left|A_g\,A_f\right|=\left|A_f\right|\,\left|A_g\right|= (c\,\log_q{m})\,\left|A_g\right|\leq
$$
$$
\leq 1.2(\log_q{m})\,(\log_m{(2\,x+1)})=1.2\,\log_q{(2\,x+1)}\leq 1.2\,\log{(2\,x+1)}.\\
$$

This together with a glance at Table T4.5.4 proves (b).
\epr

\textbf{Corollary 4.5.5.} In the notations of (4.4.2) we have, except for groups from (4.3.1)(a):\\

(a) $\left|A\right|<((m^b-1)/2)^2$,\\

(b) $\left|A_g\,A_f\right|\leq (m^b-1)/2$ except for $A_1(4), A_2(2), B_2(2)$.\\

\prf{} Set $x=(m^b-1)/2$. The inequality $4.8\,\log{(2\,x+1)}< x^2$ for $x\geq 4$ (i.e. $m^b\geq 9$) together with a glance at Table T4.5.4 implies (a). Similarly, the inequality $1.2\,\log{(2\,x+1)}< x$ for $x\geq 4$ and another glance at Table T4.5.4 yield (b).
\epr

\textbf{4.5.6.} Remark. (4.5.5) is much rougher than (4.5.4). However when we try to extend our estimates to products of groups (in Section 9) the use of logarithmic estimates for factors still leads (at least by our methods) to power estimates for the product.

\newpage

\section{\textbf{Estimates for groups of Lie type in their characteristic.}}

Let $k$ be an algebraically closed field of characteristic $p\geq 2$ and $L\simeq$ $^{c}X_a(m^c), m^c=p^s,$ a universal finite group of Lie $p$-type. 

Consider a non-trivial irreducible representation $\varphi: L\rightarrow GL_n(k)$ and set $N: N_{GL_n(k)}L$. We have the following chain of natural homomorphisms:
$$
N \rightarrow \Aut{L} \rightarrow \textnormal{Out}\,L \rightarrow A_f
$$
(where we use the notation of (4.5)). Let $N_f$ be the image of $N$ in $A_f$. 

\prop{} $n\geq d^t$ where $d$ is given in Table T5.1 below and $t:=\max{\left\{1, \left|N_f\right|/c\right\}}.$
\eprop

\ \\

\begin{center}
Table T5.1\\
\ \\
\end{center}

\begin{center}
\begin{tabular}{|c|c|c|c|c|c|c|c|c|c|}
\hline
& \ & \ & \ & \ & \ & \ & \ & \ & \\ 
$X_a$ & $A_a$ & $B_a$ & $C_a$ & $D_a$ & $E_6$ & $E_7$ & $E_8$ & $F_4$ & $G_2$ \\
\hline

& \ & \ & \ & \ & \ & \ & \ & \ & \\ 
$d$ & $a+1$ & $2\,a$ & $2\,a$ & $2\,a$ & $27$ & $56$ & $240$ & $24$ & $6$ \\  
\hline 
\end{tabular}
\end{center}

\ \\

\ \\
\cor{} 
The image of $N$ in Out $L$ has order $\leq 
\begin{cases}
6& \text{if $X_a=A_2, n=3$} \\
n\,\log{n} & \text{otherwise} \\
\end{cases}$
\ecor

\prf{} 
By (5.1) we have $\left|N_f\right|/c\leq \log_d{n}$. Therefore by (4.5.2)
$$
\left|N\right|\leq \left|A_d\right|\cdot\left|A_g\right|\cdot\left|N_f\right|\leq c\,\left|A_d\right|\,\left|A_g\right|\cdot\log_d{n}.
$$ 
Comparing tables T4.4 (where we have to take $n=a$) and T5.1 we see
$$
c\,\left|A_d\right|\,\left|A_g\right|\leq 
\begin{cases}
3\,d & \text{if $X_a=D_4$} \\
2\,d & \text{otherwise} \\
\end{cases}
$$
Thus for $X_a=D_4$ we have
$$
\left|N\right|\leq 3\,d\,\log_d{n}=24\cdot \log_8{n}= (24/\log{8})\,\log{n}=8\,\log{n}\leq n\,\log{n}.
$$
in the remaining cases 
$$
\left|N\right|\leq 2\,d\,\log_d{n}=2\,(d/\log{d})\,\log{n}\leq (2/\log{d})\,n\,\log{n}.
$$
When $d\geq 4$ this gives $\left|N\right|\leq n\,\log{n}$. If $d=3$ then the type is $A_2$ and $\left|N\right|\leq 2\,(d/\log{d})\,\log{n}=(2/\log{3})\,3\,\log{n}$
$$
\begin{cases}
=6 & \text{if $n=d=3$} \\
\leq (2/\log{3})\cdot (3\,n/4)\,\log{n} < n\,\log{n} & \text{if $n\geq 4$} \\
\end{cases}
$$
Finally, if $d=2$ then the type is $A_1$ and $c\,\left|A_d\right|\,\left|A_g\right|=2$ whence $\left|N\right|\leq 2\,\log{n}\leq n\,\log{n}$. This concludes the proof of (5.2). 
\epr

\lem{}
Let $G$ be an algebraic $k$-group of type $X_a$ and $\varphi: G\rightarrow GL_n$ its non-trivial irreducible rational representation over $k$. Then $n\geq d$ where $d$ is as in Table T5.1.
\elem

\prf{}
Let $b$ be the highest weight of $\varphi$ and let $R:=R(\varphi)$ be the set of all weights of $\varphi$. Then the Weyl group of $G$ acts on $R$ and, therefore, $\left|Wb\right|\leq \left|R\right|\leq n$. Let $P$ be a parabolic subgroup of $G$ corresponding to $b$ (the stabilizer of the weight space of weight $b$) and $W_p$ the subgroup of $W$ corresponding to $P$. Then $W_p$ is the stabilizer of $b$. Hence $n\geq \left|W/W_p\right|$. an easy case analysis gives that $\min{\left\{\left|W/W_p\right|, P \ \textnormal{parabolic} \right\}}=d$ whence our claim.
\epr

\subsection{}\textit{Proof of 5.1.} For a representation $\varphi$ of $L:=$ $^{c}X_a(m^c)$ and a homomorphism $\alpha: \mathbb{F}_{m^c}\rightarrow k$ (such $\alpha$ is, automatically, a power of the Frobenius $\mathbb{F}r$) one can define a new (in general) representation $\varphi\circ\alpha$ of $L$ by $\varphi\circ\alpha(l)=\varphi(\alpha(l))$ (see R. Steinberg [\cite{?}, 5] or [\cite{?}, 12.13]). 

By R. Steinberg [\cite{?}, Theorems 7.4, 9.3, 12.2] there exists a set $M$ of irreducible representations of $L$ over $k$ such that every other representation $\varphi$ can be \textit{uniquely} obtained as a tensor product 
$\varphi \simeq$ $\otimes_{i=0}^r\varphi_i\circ \mathbb{F}r^{i}$, where $\varphi_i\in M, i=0, 1, \ldots, r$ and $r=(s/c)-1$ unless $^{c}X_a=$ $^{2}B_2,$ $^{2}F_4,$ $^{2}G_2,$ and in this latter case $r=s-1$. 

Now let $\bar{x}$ be a generator of $N_f$ and $x$ its preimage in $N$. The $x$ acts on $L$ as $y\cdot\mathbb{F}r^z$ where $y\in\ker{(N\rightarrow A_f)}$ and $0\leq z\leq s$. It is clear that replacing in the above decomposition $\varphi=\otimes\varphi_i\circ \mathbb{F}r^i$ the maps $\mathbb{F}r^i$ by $y_i\circ\mathbb{F}r^i$ where $y_i$ are \textit{fixed} (for every $i=0, 1, \ldots, r$) automorphisms of $L$ from $\ker{(\Aut{L}\rightarrow A_f)}$ does not affect the claim. Thus there is still uniqueness and existence of decompositions $\varphi=\otimes_{i=0}^r\varphi_i\circ y_i\circ \mathbb{F}r^i$.

Let us take $\varphi$ to be our representation from the beginning of this Section. Every $i=0, 1, \ldots, r$ we write as $i=i_1+z\,i_2$ with $0\leq i_1 < z$. Then we set $y_i\circ \mathbb{F}r^i:= (y\circ\mathbb{F}r^z)^{i_2}\circ\mathbb{F}r^{i_1}.$ Since the action of $x$ normalizes $\varphi$ and in view of uniqueness of the tensor product decomposition we must have $\varphi_{j_1}\simeq \varphi_{j_2}$ of $j_1\equiv j_2 \left(\mod{z}\right)$ and $0\leq j_1, j_2\leq r$. This shows that the tensor product, if non-trivial, contains at least as many non-trivial as there are integral multiples of $z$ between $0$ and $r$. In view of the expression for $r$ given above we see that this number is $\geq z/c$. 

Thus $\dim{\varphi} \geq \max{\left\{d_1^{z/c}, d_1\right\}}$ where $d_1$ is the minimal dimension of a non-trivial irreducible $k$-representation of $L$. By R. Steinberg [\cite{?}, Theorem 43] each irreducible $k$-representation of $L$ is a restriction of one of $G$ (where $G$ is as in (5.3)) whence by (5.3) $d_1\geq d$ and (5.1) is proved.$\Box$

\subsection{}Our proof gives an apparently stronger statement. Define the action of Out $L$ on the set of equivalence classes of representations of $L$ by $(\varphi\circ\alpha)(l):= \varphi(\tilde{\alpha}(l))$ for $\varphi$ an representation, $\alpha\in$ Out $L, l\in L;$ here $\tilde{\alpha}$ is a lift of $\alpha$ to $\Aut{L}$. Let $($Out $L)_{\varphi}$ be the stabilizer of the equivalence class of $\varphi$ in Out $L$.

\prop{} 
In the notation of (5.1)
$$
\left|\left(\textnormal{Out} L\right)_\varphi\right|\leq 
\begin{cases}
6& \text{if $X_a=A_2, n=3$} \\
n\,\log{n} & \text{otherwise.} \\
\end{cases}
$$
\eprop

\newpage
\section{\textbf{Estimates for groups of Lie type in non-equal characteristic.}}

Let $k$ be an algebraically closed field of characteristic exponent $p=p(k)$. Let $G$ be a finite simple group of Lie $q$-type, $q\neq p, q$ a prime. To avoid trouble with different characteristics (see (4.3.2)) we fix an isomorphism $G\simeq$ $^{c}\bar{X}_a(m^c)$ and write $m^c=q^s.$ Set $H:=\mathcal{D}G$. Set $f(t):= (2\,t+1)^{2\,\log_3{(2\,t+1)}+1}.$

\prop{}
Let $\varphi:= H\rightarrow PGL_n(k)$ be a faithful irreducible projective representation. Then \\

(a) $\left|H\right|\leq f(n)$ except for the cases
$$
H\simeq \text{$^{2}\bar{A}_3(9)$ and $n=6$ when $\left|H\right|\leq 1.58\,f(6)$},
$$
$$
\text{and $H\simeq D_4(2)$ and $n=8$ when $\left|H\right|\leq 4.62\,f(8)$}\\
$$

(b) $\left|\Aut{H}\right|\leq n\,f(n)$ except for the following cases\\

\begin{tabular}{cccccc}
$n =$ & $2$ & $2$ & $4$ & $6$ & $8$ \\
$H\simeq$ & \textnormal{Alt}$_6$ & $A_1(8)$ & $\bar{A}_2(4)$ & $^{2}\bar{A}_3(9)$ & $D_4(2)$ \\
$\left|\Aut{H}\right|\leq$ & $2.6\,f(2)$ & $2.71\,f(2)$ & $4.1\,f(4)$ & $12.61\,f(6)$ & $27.69\,f(8)$ \\
\end{tabular}
\ \\

(c) $\left| \textnormal{Out}\,H \right|\leq n^2$,\\

(d) $\left|N_{PGL_n(k)}(\varphi(H))\right|\leq n\,f(n)$ except when $H\simeq D_4(2)$ and $n=8$.
\eprop

\lem{}
If $H$ is not centrally isomorphic to one of the groups listed in (4.3.1), (4.3.2), and (4.4.2) then
$$
\left|^{c}X_a(m^c)\right|\leq f(n).
$$
\elem

\prf{}
By (4.4.2) $n\geq (m^b-1)/2$, i.e. $m^b\leq 2\,n+1$. Since $^{2}B_2(2),$ $^{2}F_4(2),$ and $^{2}G_2(3)$ are excluded by (4.3.1) we have $m\geq 2$ and, therefore, $b\leq \log{(2\,n+1)}$. If we exclude type $C_r, r\geq 2,$ in odd characteristic we have by (4.4.1) and (4.4.3)(b)
$$
\left|^{c}X_a(m^c)\right|\leq m^d\leq m^{b^2+2\,b}=(m^b)^{b+2}\leq (2\,n+1)^{\log{(2\,n+1)}+2}.
$$

For the type $C_r, r\geq 2$ in odd characteristic we have $m\geq 3$ whence $b\leq \log_3{(2\,n+1)}$ whence by (4.4.1) and (4.4.3)(a)
$$
\left|^{c}X_a(m^c)\right|\leq m^d\leq m^{2\,b^2+b}=(m^b)^{2\,b+1}\leq (2\,n+1)^{2\,\log_3{(2\,n+1)}+1}.
$$

Note that the first occurrences of this latter case are for $(r, m)=(2, 3)$ (resp. $(2, 5), (3, 3)$), $(m^b-1)/2=4$ (resp. $12, 13$).

In general we have, therefore, that 
$$
\left|^{c}X_a(m^c)\right|\leq \max{\left\{(2\,n+1)^{2\,\log_3{(2\,n+1)}+1},
(2\,n+1)^{\log{(2\,n+1)}+2}\right\}}.
$$
One easily sees that $(2\,n+1)^{2\,\log_3{(2\,n+1)}+1} \geq (2\,n+1)^{\log{(2\,n+1)}+2}$ if $n\geq 7$. For $n\leq 6$ the difference can come only from groups of type $C_r, r\geq 2$, in odd characteristic. As we remarked the first such groups occur at $n=4, 12, 13$. One has $\left|B_2(3)\right|=2\cdot25,920\leq f(4)$. Therefore both functions give valid estimates for $n\leq 12$. Therefore $f(n)$ can be taken for an estimate for all $n$.
\epr
 
\newpage
\pagestyle{empty} 

\begin{center}
Table T6.3.\\
\end{center}
\ \\

\tiny
\begin{tabular}{|c|c|c|c|c|c|c|c|c|}
\hline
& $M$ is & Order & \  & \  & \  & $\min$ & \ & Schur \\
Group & centrally & modulo & min & $LS$ & Adjusted & non-Lie- & Out $M$ & multiplier \\
$M$ & isomorphic to & center & $n$ & estimate & estimate & $p$-type $n$ & \ & of $M/$center \\
\hline
& \ & \ & \ & \ & \ & \ & \ & \ \\
$A_1(4)$ & $A_1(5),$ Alt$_5$ & 60 & 2 & 2 & 2 & 2 & $\mathbb{Z}/2(=A_f)$ & $\mathbb{Z}/2$ \\
\hline 

& \ & \ & \ & \ & \ & \ & \ & \ \\
$A_1(5)$ & $A_1(4),$ Alt$_5$ & 60 & 2 & 2 & 2 & 2 & $\mathbb{Z}/2(=A_d)$ & $\mathbb{Z}/2$ \\
\hline 

& \ & \ & \ & \ & \ & \ & \ & \ \\
$A_1(7)$ & $A_2(2)$ & 168 & 2 & $3(p\neq 7)$ & $2^{b)}$ & $3^{d)}$ & $\mathbb{Z}/2(=A_d)$ & $\mathbb{Z}/2$ \\
\hline 

& \ & \ & \ & \ & \ & \ & \ & \ \\
$A_1(8)$ & $\mathcal{D}$$(^{2}G_2(3))$ & 504 & 2 & $7(p\neq 2)$ & $2^{b)}$ & $7$ & $\mathbb{Z}/3(=A_f)$ & $\mathbb{Z}/2$ \\
\hline

& \ & \ & \ & \ & \ & \ & \ & \ \\
$A_1(9)$ & $\mathcal{D}B_2(2),$ Alt$_6$ & 360 & 2 & $3(p\neq 3)$ & $2^{b)}$ & $3$ & $\mathbb{Z}/2\times\mathbb{Z}/2$ & $\mathbb{Z}/2\times\mathbb{Z}/3$ \\
& \ & \ & \ & \ & \ & \ & $(=A_d\times A_f)$ & \ \\
\hline

& \ & \ & \ & \ & \ & \ & \ & \ \\
$A_2(2)$ & $A_1(7)$ & 168 & 2 & $2(p\neq 2)$ & $2^{b)}$ & $3^{d)}$ & $\mathbb{Z}/2(=A_g)$ & $\mathbb{Z}/2$ \\
\hline

& \ & \ & \ & \ & \ & \ & \ & \ \\
$A_3(2)$ & Alt$_8$ & 20,160 & 4 & $7(p\neq 2)$ & $4^{c)}$ & $7$ & $\mathbb{Z}/2(=A_g)$ & $\mathbb{Z}/2$ \\
\hline

& \ & \ & \ & \ & \ & \ & \ & \ \\
$\mathcal{D}B_2(2)$ & $A_1(9),$ Alt$_6$ & 360 & 2 & $2(p\neq 2)$ & $2^{b)}$ & $3^{d)}$ & $\mathbb{Z}/2\times\mathbb{Z}/2$ & $\mathbb{Z}/2$ \\
& \ & \ & \ & \ & \ & \ & $(=A_g\times(H/\mathcal{D}H))$ & \ \\
\hline

& \ & \ & \ & \ & \ & \ & \ & \ \\
$B_2(3)$ & $^{2}A_3(4)$ & 25,920 & 4 & $4(p\neq 3)$ & $4^{b)}$ & $4^{d)}$ & $\mathbb{Z}/2(=A_d)$ & $\mathbb{Z}/2$ \\
\hline

& \ & \ & \ & \ & \ & \ & \ & \ \\
$\mathcal{D}G_2(2)$ & $^{2}A_2(9)$ & 6,048 & 3 & $3(p\neq 2)^{a)}$ & $3^{b)}$ & $6^{d)}$ & $\mathbb{Z}/3(=H/\mathcal{D}H)$ & $1$ \\
\hline

& \ & \ & \ & \ & \ & \ & \ & \ \\
$^{2}A_2(9)$ & $\mathcal{D}G_2(2)$ & 6,048 & 3 & $3(p\neq 3)$ & $3^{b)}$ & $6^{d)}$ & $\mathbb{Z}/3(=A_f)$ & $1$ \\
\hline

& \ & \ & \ & \ & \ & \ & \ & \ \\
$^{2}A_3(4)$ & $B_2(3)$ & 25,920 & 4 & $4(p\neq 2)$ & $4^{b)}$ & $4^{d)}$ & $\mathbb{Z}/2(=A_f)$ & $\mathbb{Z}/2$ \\
\hline

& \ & \ & \ & \ & \ & \ & \ & \ \\
$\mathcal{D}$$(^{2}G_2(3))$ & $A_1(8)$ & 504 & 2 & $2(p\neq 3)^{a)}$ & $2^{b)}$ & $7$ & $\mathbb{Z}/3(=H/\mathcal{D}H)$ & $\mathbb{Z}/2$ \\
\hline

& \ & \ & \ & \ & \ & \ & \ & \ \\
$A_2(4)$ & \  & 20,160 & 4 & $4(p\neq 2)$ & same & same  & $\Sym_3\times\mathbb{Z}/2$ & $\mathbb{Z}/3\times\mathbb{Z}/4$ \\
& \ & \ & \ & \ & as $LS$ & as $LS$ & $(=A_d\,A_g\,A_f)$ & $\times\mathbb{Z}/4$ \\
\hline

& \ & \ & \ & \ & \ & \ & \ & \ \\
$B_3(3)$ & \ & $9.17\cdot 10^9$ & 13 & $27(p\neq 3)$ & " & " & $\mathbb{Z}/2(=A_d)$ & $\mathbb{Z}/2\times\mathbb{Z}/3$ \\
\hline

& \ & \ & \ & \ & \ & \ & \ & \ \\
$D_4(2)$ & \ & $1.74\cdot 10^8$ & 10 & $8(p\neq 2)$ & same & $8^{f)}$ & $\Sym_3(=A_g)$ & $\mathbb{Z}/2\times\mathbb{Z}/2$ \\
& \ & \ & \ & \ & as $LS$ & \ & \ & \ \\
\hline

& \ & \ & \ & \ & \ & \ & \ & \ \\
$F_4(2)$ & \ & $33\cdot 10^{15}$ & 33 & $44(p\neq 2)$ & " & same & $\mathbb{Z}/2(=A_g)$ & $\mathbb{Z}/2$ \\
& \ & \ & \ & \ & \ & as $LS$  & \ & \ \\
\hline

& \ & \ & \ & \ & \ & \ & \ & \ \\
$G_2(4)$ & \ & $2.5\cdot 10^8$ & 10 & $12(p\neq 2)^{e)}$ & " & " & $\mathbb{Z}/2(=A_f)$ & $\mathbb{Z}/2$ \\
\hline

& \ & \ & \ & \ & \ & \ & \ & \ \\
$^{2}A_3(9)$ & \  & $3.26\cdot 10^6$ & 7 & $6(p\neq 3)$ & " & $6^{d)}$  & $\mathbb{Z}/4\times\mathbb{Z}/2$ & $\mathbb{Z}/4\times\mathbb{Z}/3$ \\
& \ & \ & \ & \ & \ & \  & $(=A_d\cdot A_f)$ & $\times\mathbb{Z}/3$ \\
\hline

& \ & \ & \ & \ & \ & \ & \ & \ \\
$^{2}B_2(8)$ & \ & $29,120$ & 4 & $8(p\neq 2)$ & " & same & $\mathbb{Z}/3(=A_f)$ & $\mathbb{Z}/2\times\mathbb{Z}/2$ \\
& \ & \ & \ & \ & \ & as $LS$  & \ & \ \\
\hline

& \ & \ & \ & \ & \ & \ & \ & \ \\
$^{2}E_6(4)$ & \  & $2.3\cdot 10^{23}$ & 88 & $1500(p\neq 2)$ & " & " & $\mathbb{Z}/2(=A_f)$ & $\mathbb{Z}/3\times\mathbb{Z}/2$ \\
& \ & \ & \ & \ & \ & \ & \ & $\times\mathbb{Z}/2$ \\
\hline

& \ & \ & \ & \ & \ & \ & \ & \ \\
$\mathcal{D}$$(^{2}F_4(2))$ & \ & $1.8\cdot 10^7$ & 8 & $8(p\neq 2)^{a)}$ & " & " & $\mathbb{Z}/2(=H_d/\mathcal{D}H)$ & $1$ \\
\hline
\end{tabular}

\newpage
\normalsize
\pagestyle{plain}
\subsection{}\textit{Proof of (6.1)(a)} for the cases omitted in (6.2) is contained in Table T6.3 and for the case $D_4(2)$ which requires special attention in (6.3.2). The upper portion of Table T6.3 handles groups of Lie type which can appear in two different characteristics (see (4.3.2)) or are isomorphic to alternating groups. The lower portion treats groups for which $LS$-estimates ($LS$ stands for Landazuri-Seitz, see (4.4.2)) do not have the form $(m^b-1)/2$ with $b$ from Table T4.4 (see exceptions in (4.4.2)). Column 4 gives minimal $n$ for which $\left|M/\textnormal{center}\right|\leq f(n)$ (see (6.3.1) below for explicit formula). Column 6 gives an adjusted estimate on dimensions of projective representations of $M$; explanations are given in notes a)-c) below. If column 6 is larger than $\min {n}$ then (6.1)(a) holds for $M/\textnormal{center}$. Column 7 gives an estimate (still from below) on dimensions of projective representations of $M$ in such characteristics $p$ for which $M/\textnormal{center}$ is not isomorphic to a group of Lie $p$-type; explanations are given in notes below.\\

Explanations:\\

a) We took an estimate for $G$ in V. Landazuri and G. Seitz [\cite{18}, p.419] and divided it by $\left|G/\mathcal{D}G\right|$ to obtain an estimate for $M=G$.\\

b) The estimate is the minimum over isomorphic groups on the same line as $M$.\\

c) See (5.1).\\

d) See Table T2.7 (groups of small degree).\\

e) The estimate given in V. Landazuri and G. Seitz [\cite{18}, p.419] is incorrect for $G_2(4)$. The correct estimate is $12$ (personal communication of G. Seitz who also explained how a slip in the proof of Lemma 5.6(b) of the above paper should be corrected).\\

f) $D_4(2)$ has an $8$-dimensional projective representation as the derived group of the Weyl group of $E_8$, see R. Steinberg [\cite{?}, $\S 11$, after Theorem 37].\\

\textbf{Lemma 6.3.1.} $\min {n}$ is the smallest integer $\geq (3^x-1)/2$ where
$$
x:=\frac{1}{4}(-1+\sqrt{0.5+8\,\log_3{\left|M/\textnormal{center}\right|}}).
$$   

\prf{}
Setting $x=2\,n+1$ we have to solve $2\,x^2+x \leq \log_3{\left|M/\textnormal{center}\right|}$ whence our claim.
\epr

\textbf{Lemma 6.3.2.} $D_4(2)$ does not have faithful irreducible projective representations of dimension 9 over fields of characteristic $p\neq 2$. 

\prf{}
Suppose $\varphi: \tilde{D}_4(2)\rightarrow GL_9(k), p\neq 2,$ is an irreducible representation. Since the center of $\varphi(\tilde{D}_4(2))$ will be contained then in the center of $SL_9(k)$ which is isomorphic to a subgroup of 9-th roots of 1 and since the center of $\tilde{D}_4(2)$ is isomorphic to $\mathbb{Z}/2\times\mathbb{Z}/2$, it follows that the center of $\varphi(\tilde{D}_4(2))$ is trivial. Thus $\varphi$ is, in fact, a representation of $D_4(2)$. But then the proof of V. Landazuri and G. Seitz [\cite{18}, Lemma 3.3(2)] gives that $D_4(2)$ has no representations over $k$ of degree $\leq 27$, whence our claim.
\epr 

\subsection{}To prove (6.1)(b) we use (4.5.3). Except when $H$ is in the Table T6.3 this gives us that $\left|\Aut{H}\right|=|H|\cdot |A_d|\cdot |A_f|\cdot |A_g|=$$\left|^{c}X_a(m^c)\right|\cdot |A_f|\cdot |A_g|$ whence by (4.5.5)(b), (4.4.2) and (6.2) $\left|\Aut{H}\right|\leq n\,\left|^{c}X_a(m^c)\right|$ whence (6.1)(b). When $H$ is in the Table T6.3 one has to verify (6.1)(b) directly (see Table TA for values of $f(n)$ for small $n$).

\subsection{}The claim of (6.1)(c) is contained in (4.5.4)(a) when (4.4.2) holds for $H$. in the remaining cases one uses Table T6.3 to verify the claim directly.

\subsection{}Of course, (6.1)(d) follows from (6.1)(b) for all but four cases.\\ 
Set $N:= \left|NZ_{PGL_n(k)}(\varphi(H))\right|.$\\

\begin{tabular}{ccccccc}
\textbf{Lemma }& If & $n=$ & 2 & 2 & 4 & 6 \\
& and & $H\simeq$ & Alt$_6$ & $A_1(8)$ & $\bar{A}_2(4)$ & $^{2}\bar{A}_3(9)$ \\
& then & $N\leq$ & 720 & 504 & 40,320 & 6,531,840 \\
\end{tabular}
\\ \

\subsubsection{} If $H\simeq \textnormal{Alt}_6$ and $n=2$ then comparing rows $A_1(9)$ and $\mathcal{D}B_2(2)$ of Table T6.3 we see that $p=3$ so that ours is the natural representation of $SL_2(\mathbb{F}_9)$. The normalizer of $SL_2(\mathbb{F}_9)$ is $GL_2(\mathbb{F}_9)\cdot k^{*}$, whence $NZ_{PGL_2(k)}\varphi(H)\simeq PGL_2(\mathbb{F}_9)$ and $N=2\,|H|$.

\subsubsection{} If $H\simeq A_1(8)$ and $n=2$ then comparing rows $A_1(8)$ and $\mathcal{D}(^{2}G_2(3))$ of Table T6.3 we see that $p=2$ and ours is the natural representation of $SL_2(\mathbb{F}_8)$. The normalizer of $SL_2(\mathbb{F}_8)$ is $GL_2(\mathbb{F}_8)\cdot k^{*}$, whence $N=|H|$. 

\subsubsection{} Let now $n=4$ and $\tilde{H}\subseteq GL_4(k)$ be a perfect group centrally isomorphic to $H\simeq \bar{A}_2(4)\simeq PSL_3(\mathbb{F}_4)$. By Table T2.7 (groups of small degree) we have $p=\textnormal{char}{k}=3$. Then $H$ contains a subgroup $T$ isomorphic to $\ker{\left\{N_{\mathbb{F}_{64}/\mathbb{F}_4}: \mathbb{F}_{64}^{*}\rightarrow \mathbb{F}_4^{*}\right\}}$. Clearly $T\simeq \mathbb{Z}/7$ and $N_H(T)/T$ acts on $T$ as Gal$(\mathbb{F}_{64}/\mathbb{F}_4)\simeq \mathbb{Z}/3$. Let $\tilde{T}$ be the 7-component of the preimage of $T$ in $\tilde{H}$. Since $\tilde{T}$ is a Sylow 7-subgroups of $\tilde{H}$ we have (by Frattini argument that $N_{GL_4(k)}(\tilde{H})=\tilde{H}\cdot N_{GL_4(k)}(\tilde{T})$). We have that $NZ_{GL_n(k)}(\tilde{T})/\tilde{T}$ contains $\mathbb{Z}/3$ and is contained in $\mathbb{Z}/6$. By representation theory of Frobenius groups we have that $|NZ(\tilde{T})/\tilde{T}|\leq 4$ whence $|NZ(\tilde{T})/\tilde{T}|=3$ and $NZ(\tilde{T})\subseteq H$. Since $A_g$ can be assumed to be the transpose-inverse on $T$ we have that then the image of $NZ(\tilde{T})$ in Out $\tilde{H}=A_d\cdot A_g\cdot A_f$ does not contain $A_g$. Thus this image is subgroup of $A_d\cdot A_f\simeq \Sym_3$. If it is the whole group $\simeq\Sym_3$ then $Z_{GL_4(k)}(\tilde{T})$ contains a subgroup centrally isomorphic to $\Sym_3$. Since all eigenspaces of $\tilde{T}$ are $1$-dimensional this is impossible so that the image $S$ of $NZ(\tilde{T})$ in Out $\tilde{H}$ is of order 2 or 3. 

If $|S|=3$ then take the 3-component $\tilde{S}$ of the preimage of $S$ in $H$. It commutes with $\tilde{T}$ whence by multiplicity 1 of eigenspaces of $\tilde{T}$ we have that $\tilde{S}$ is diagonalizable. But since $|\tilde{S}|=3$ and char $k=3$ it follows that $\tilde{S}$ is unipotent. Thus $\tilde{S}=\{1\}$.

\subsubsection{} Let $n=6$ and let $\tilde{H}\subseteq GL_6(k)$ be a perfect group centrally isomorphic to $H\simeq$ $^{2}\bar{A}_3(9)\simeq PSU_4(\mathbb{F}_9)$. We have $|H|=2^7\cdot 3^6\cdot 5\cdot 7=3,265,920$. We assume that $p\neq 3$, i.e., char $k\neq 3$. 

Let $\tilde{H}\simeq$ $^{2}\tilde{A}_3(9)$. Then $Z:=\ker{(^{2}\tilde{A}_3(9)\rightarrow {^{2}A_3(9)})}\simeq (\mathbb{Z}/3)^2$. By D. Gorenstein and R. Lyons [\cite{18}, (7-8)(3)] Out $\tilde{H}\simeq (\mathbb{Z}/4)\cdot(\mathbb{Z}/2)$ (the dihedral group) acts faithfully on $Z$. It is easy to check then, using elementary representation theory of Out $\tilde{H}$ on $Z(:=$ the dual of $Z)$, that the stabilizer of a point $z\in Z$ in Out $\tilde{H}$ is isomorphic to $\mathbb{Z}/2$.

Let $\varphi: \tilde{H}\rightarrow GL_6(k), \textnormal{char}\ k\neq 3$, be an irreducible representation. If $\varphi(Z)=\Id$ then the 3-Sylow subgroup of $\varphi(\tilde{H})$ is the same as that of $H\simeq PSU_4(\mathbb{F}_9)$ and it contains an extraspecial 3-subgroup $P$ of order 3 (consisting of matrices 
$
\begin{bmatrix} 
1 & \  & \ & 0 \\ 
a & 1 & \  & \  \\
b & 0 & 1 & \ \\
c & \bar{b} & \bar{c} & 1 \\
\end{bmatrix}
$
, $a, b \in \mathbb{F}_9, c\in\mathbb{F}_3$, in an appropriate basis). By representation theory of $P$ any of its irreducible faithful representations is of degree $9 > 4$. Thus the case $\varphi(Z)=\Id$ is impossible. 

Since $\varphi(\tilde{H})$ is irreducible  it follows then that $\varphi(Z)\simeq \mathbb{Z}/3$ and consists of scalar matrices. In particular, $\varphi(Z)$ is in the center of $N:=N_{GL_6(k)}(\varphi(\tilde{H}))$ whence the image $S$ of $N$ in Out $\tilde{H}$ acts trivially on $\varphi(Z)$ whence $\varphi|Z(\in Z)$ is stable under $S$. By one of the remarks above we have then $|S|\leq 2$, as desired. 
\subsection{}One may need (and we shall in Section 16) a variation of (6.1)(a):\\

\textbf{Proposition.} Let $G$ be centrally simple perfect finite subgroup of $GL_n(k)$ of Lie $q$-type, $q\neq p$. Then $|G|\leq f(n)$ except in the following cases\\

\begin{tabular}{ccccc}
$n=$ & 2 & 4 & 6 & 8 \\
$G\simeq$ & $2\cdot\mathcal{D}B_2(2)$ & $4\cdot\bar{A}_2(4)$ & $6\cdot$$^{2}\bar{A}_3(9)$ & $2\cdot D_4(2)$ \\
$G=$ & 720 & 80720 & $1.96\cdot 10^7$ & $3.48\cdot 10^8$ \\
$p=$ & 3 & $\neq 2$ & $\neq 3$ & $\neq 2$\\
$q=$ & 2 & 2 & 3 & 2\\ 
\end{tabular}

\prf{}
If $G$ is not centrally isomorphic to a group from (4.3.3)(a) or from Table T6.3 then our claim follows from (6.2). Let $\tilde{G}$ be the universal cover of $G$. If $G$ is in the Table T6.3 one readily verifies that $|\tilde{G}|\leq f(n)$ except for the pairs $(\mathcal{D}(^{2}\tilde{G}_2(3)), 2)$(for $p=2$), $(\mathcal{D}\tilde{B}_2(2), 2)$(for $p=3$), $(\tilde{A}_2(4), 4), (\tilde{A}_2(4), 5), (^{2}\tilde{A}_3(9), 6), (^{2}\tilde{A}_3(9), 7), (^{2}\tilde{A}_3(9), 8),$ $ (\tilde{D}_4(2), 8), (\tilde{D}_4(2), 9), (\tilde{D}_4(2), 10).$ Now note that since $G$ is perfect, its center $C$ is contained in the group of $n$-th roots of 1. Thus $C\simeq \mathbb{Z}/m$, where $m|(n/p^a)$ where $p^a$ is the highest power of $p$ dividing $n$. Applying this information to the pairs $(\tilde{G},n)$ given above and taking into account (6.3.2) one is left only with pairs listed in (6.7). 

Now it remains to check the groups from (4.3.3)(a) which are not contained in Table T6.3. These groups, their orders, their centers, and the estimates for the from V. Landazuri and G. Seitz [\cite{18}, p.419] are \\

\begin{tabular}{cccc}
$\tilde{G}$ & $\tilde{B}_3(2)$ & $\tilde{G}_2(3)$ & $^{2}\tilde{A}_5(4)$ \\
$|\tilde{G}|$ & $2.9\cdot 10^6$ & $1.27\cdot 10^7$ & $1.1\cdot 10^{11}$ \\
Center $(\tilde{G})$ & $\mathbb{Z}/2$ & $\mathbb{Z}/3$ & $\mathbb{Z}/3\times\mathbb{Z}/2\times\mathbb{Z}/2$ \\
$LS$-estimate & 7 & 14 & 21 \\
$f(LS-\text{estimate})$ & $9.42\cdot 10^6$ & $2.67\cdot 10^7$ & $6.58\cdot 10^{12}$ 
\end{tabular}
\epr

\newpage
\section{\textbf{Estimates for sporadic groups.}}

Our purpose here is to obtain estimates from below on the minimal degrees of faithful irreducible projective representations of the 26 sporadic groups. This Section uses references which, often, are not used in other parts of the paper. For this reason we referred directly to many papers without placing them in the list of reference at the end of paper. 

My knowledge of sporadic groups is very scanty and, in many cases, the literature on them contains enormous gaps covered by references to unpublished work, lectures, and personal communications. I had good fortune of being helped in many cases where I was lost by Robert Griess. I am extremely grateful to him; the arguments he supplied are marked so, some other data refers directly to his \textit{personal communications}. However some of his suggestions were later superseded. Therefore we refer directly only to portion of the many arguments he offered. 

We want to keep the estimate $f(x):=(2\,x+1)^{2\,\log_3{(2\,x+1)}+1}$. This is, however, impossible (as it was for $^{2}A_3(9)$ and $D_4(2),$ see (6.1)), for Suzuki's group Suz and Conway's groups $\cdot 1$ and $\cdot 2$. To give an estimate for these groups we introduce cumbersome functions $F(H,n)$. To define them we use notation (for $a_1, a_2, b\in \mathbb{R}$)
$$
y_{a_1, a_2, b}(t)=
\begin{cases}
1 & t < a_1 \\
b & a_2 > t \geq a_1 \\
f(t) & t \geq a_2
\end{cases}
$$  

We also use $\delta_{a, b}$ for the Christoffel symbol. 

\th{}
Let $k$ be an algebraically closed field of characteristic exponent $p=p(k)$ and $\varphi : H\rightarrow PGL_n(k)$ a faithful irreducible projective representation of a sporadic simple group $H$.\\
 
(a) Unless $H$ is centrally isomorphic to Suz, $\cdot 1, \cdot 2,$ we have 
$$
\left|H\right|\leq f(n),\\
$$

(b) $\left|\textnormal{Suz}\right|\leq y_{12, 18, \textnormal{Suz}}(n)=: F(\textnormal{Suz}, n)$
$$
\left|\cdot 1\right|\leq y_{24, 49, \left|\cdot 1\right|}(n)=: F(\cdot 1, n)
$$
$$
\left|\cdot 2\right|\leq y_{20, 24, \left|\cdot 2\right|}(n)=: F(\cdot 2, n)
$$
\eth  

\subsection{}Our proof of (7.1) is mostly contained in Table T7.2. In this table $\min{n}$ (resp. $\min{\tilde{n}}$) is the smallest integer $m$ such that $f(m)\geq |H|$ (resp. $m^{2\,\log{m}+4.32}\geq |H|$). We will need $\min{n}$ in $\S 10$. To obtain estimates we try to find two subgroups $H_1$ and $H_2$ and to use their lower estimate to estimate $\mathcal{C}d(H)$. Column "Subgroup" says what subgroups we use and the next column gives a reference to a source where existence of these subgroups is pointed out. The next three columns describe the estimate obtained from $H_1$ and $H_2$ and give a reason why such an estimate holds. The possible reasons are listed in (7.3) below. The "adjusted estimate" is generally equal to the minimum of the estimates for $H_1$ and $H_2$. If there is another reason for taking the indicated estimate it is given next to the "adjusted estimate". The next three columns describe precise results if such are known to me. The last column gives sometimes an indirect reference; e.g. c) [51] is paper [51] in the list of references of c).

\newpage
\pagestyle{empty} 

\begin{landscape}

\tiny
\begin{center}
Table T7.2.\\
\end{center}
\ \\

\begin{center}
\begin{tabular}{|c|c|c|c|c|c|c|c|c|c|c|c|c|c|c|}
\hline

& \ & \ & \ & \ & \ & \ & \ & \ & \ & \ & \ & precise & \ & \ \\
& \ & min & min & Schur & \ & \ & \ & \ & \ & adjusted & \ & result & \ & \ \\
$H$ & $|H|$ & $n$ & $\tilde{n}$ & multi & subgroup & source & esti & condi & reason & estimate & additional & when & condi & source \\
& $(:=b)$ & $(:=a_2)$ & \ & plier & \ & \ & mate & tions & \ & $(:=a_1)$ & reasons & known & tions & \ \\
& \ & \ & \ & \ & \ & \ & \ & \ & \ & \ & \ & to me & \ & \ \\
\hline

& \ & \ & \ & \ & \ & \ & \ & \ & \ & \ & \ & \ & \ & \ \\
$M_{11}$ & 7920 & 4 & 4 & 1 & \ & \ & \ & \ & \ & 5 & see$\rightarrow$ & 5 & $p=3$ & a) \\
\hline

& \ & \ & \ & \ & \ & \ & \ & \ & \ & \ & \ & \ & \ & \ \\
$M_{12}$ & $9.5\cdot 10^4$ & 5 & 4 & 2 & \ & \ & \ & \ & \ & 6 & b) & 6 & $p=3$, proj & ? \\
\hline

& \ & \ & \ & \ & \ & \ & \ & \ & \ & \ & \ & \ & \ & Conway \\
$M_{22}$ & $4.43\cdot 10^5$ & 6 & 5 & 12 & \ & \ & \ & \ & \ & 6 & b) & 6 & $p=2$, proj & p.242 \\
\hline

& \ & \ & \ & \ & \ & \ & \ & \ & \ & \ & \ & \ & \ & \ \\
$M_{23}$ & $1.02\cdot 10^7$ & 8 & 6 & 1 & \ & \ & \ & \ & \ & 11 & see$\rightarrow$ & 11 & \ & a) \\
\hline

& \ & \ & \ & \ & \ & \ & \ & \ & \ & \ & \ & \ & \ & \ \\
$M_{24}$ & $2.45\cdot 10^8$ & 10 & 7 & 1 & \ & \ & \ & \ & \ & 11 & see$\rightarrow$ & 11 & \ & a) \\
\hline

& \ & \ & \ & \ & \ & \ & \ & \ & \ & \ & \ & \ & \ & \ \\
$J_1$ & $1.76\cdot 10^5$ & 5 & 5 & 1 & \ & \ & \ & \ & \ & 6 & b) & 56 & $p=1$ Higman & $^{\textnormal{c})}[71]$ \\
& \ & \ & \ & \ & \ & \ & \ & \ & \ & \ & \ & 7 & $p=11$ & \ \\
\hline

& \ & \ & \ & \ & \ & \ & \ & \ & \ & \ & \ & \ & \ & \ \\
$J_2$ & $6.05\cdot 10^5$ & 6 & 5 & 2 & \ & \ & \ & \ & \ & 6 & b) & 6 & $p=1$, proj & Feit? \\
\hline

& \ & \ & \ & \ & \ & \ & \ & \ & \ & \ & \ & \ & \ & \ \\
$J_3$ & $5.02\cdot 10^7$ & 9 & 7 & 1 & $19\cdot 9$ & $^{\textnormal{c})}$p.214 & 9 & $p\neq 19$ & f) & 9 & \ & 85 & $p=1$ & $^{\textnormal{c})}[73]$ \\
& \ & \ & \ & \ & $A_1(16)$ & \ & 15 & $p\neq 2$ & d) & \ & \ & \ & \ & \ \\
\hline

& \ & \ & \ & \ & \ & \ & \ & \ & \ & \ & \ & \ & \ & \ \\
$J_4$ & $8.68\cdot 10^{19}$ & 58 & 28 & 1 & $2^{1+12}$ & c)pp.215, & 64 & $p\neq 2$ & e) & 64 & \ & \ & \ & \ \\
& \ & \ & \ & \ & $N(11^{1+2})$ & 238 & 110 & $p\neq 11$ & e) & \ & \ & \ & \ & \ \\
\hline

& \ & \ & \ & \ & \ & \ & \ & \ & \ & \ & \ & \ & \ & \ \\
HS & $4.44\cdot 10^7$ & 9 & 7 & 2 & $^{2}A_2(25)$ & $^{\textnormal{c})}$p.220 & 20 & $p\neq 5$ & d) & 10 & \ & 22 & $p=1,$ linear & $^{\textnormal{c})}[42]$ \\
& \ & \ & \ & \ & $M_{11}$ & \ & 10 & $p\neq 3,11$ & a) & \ & \ & \ & \ & \ \\
\hline

& \ & \ & \ & \ & \ & \ & \ & \ & \ & \ & \ & \ & \ & \ \\
Mc & $1.28\cdot 10^8$ & 9 & 7 & 3 & $^{2}A_2(25)$ & $^{\textnormal{c})}$p.222 & 20 & $p\neq 5$ & d) & 10 & \ & \ & \ & \ \\
& \ & \ & \ & \ & $M_{11}$ & \ & 10 & $p\neq 3,11$ & a) & \ & \ & \ & \ & \ \\
\hline

& \ & \ & \ & \ & \ & \ & \ & \ & \ & \ & \ & \ & \ & \ \\
Suz & $4.48\cdot 10^{11}$ & 18 & 11 & G & $G_2(4)$ & c)pp.221, & 12 & $p\neq 2$ & r1) & 12 & \ & 12 & $p=1,$ proj & $^{\textnormal{c})}[90]$ \\
& \ & \ & \ & \ & $5^2\cdot(4\times\Sym_3)$ & 222 & 12 & $p\neq 5$ & r2) & \ & \ & 143 & $p=1,$ linear & $^{\textnormal{c})}[154]$ \\
\hline

& \ & \ & \ & \ & \ & \ & \ & \ & \ & \ & \ & \ & \ & \ \\
Ru & $1.46\cdot 10^{11}$ & 16 & 11 & 2 & $^{2}A_2(25)$ & $^{\textnormal{c})}$p.224 & 20 & $p\neq 5$ & d) & 16 & \ & \ & \ & \ \\
& \ & \ & \ & \ & $^{2}F_4(2)$ & \ & 16 & $p\neq 2$ & d) & \ & \ & \ & \ & \ \\
\hline

& \ & \ & \ & \ & \ & \ & \ & \ & \ & \ & \ & \ & \ & \ \\
He & $4.03\cdot 10^9$ & 13 & 9 & 1 & $7^2\cdot A_1(7)$ & $^{\textnormal{c})}$p.221 & 48 & $p\neq 7$ & f) & 18 & \ & \ & \ & \ \\
& \ & \ & \ & \ & $C_2(4)$ & \ & 18 & $p\neq 2$ & d) & \ & \ & \ & \ & \ \\
\hline

& \ & \ & \ & \ & \ & \ & \ & \ & \ & \ & \ & \ & \ & \ \\
Ly & $5.18\cdot 10^{16}$ & 38 & 20 & 1 & $G_2(5)$ & $^{\textnormal{c})}$p.223 & 120 & $p\neq 5$ & d) & 110 & m) & 2480 & $p=1$ & $^{\textnormal{c})}[93]$ \\
& \ & \ & \ & \ & $N(67)$ & \ & 22 & $p\neq 67$ & f) & \ & \ & \ & \ & \ \\
\hline

\end{tabular}
\end{center}

\end{landscape}

\newpage

\begin{landscape}

\tiny
\begin{center}
Table T7.2 continued.\\
\end{center}
\ \\

\begin{center}
\begin{tabular}{|c|c|c|c|c|c|c|c|c|c|c|c|c|c|c|}
\hline

& \ & \ & \ & \ & \ & \ & \ & \ & \ & \ & \ & precise & \ & \ \\
& \ & min & min & Schur & \ & \ & \ & \ & \ & adjusted & \ & result & \ & \ \\
$H$ & $|H|$ & $n$ & $\tilde{n}$ & multi & subgroup & source & esti & condi & reason & estimate & additional & when & condi & source \\
& $(:=b)$ & $(:=a_2)$ & \ & plier & \ & \ & mate & tions & \ & $(:=a_1)$ & reasons & known & tions & \ \\
& \ & \ & \ & \ & \ & \ & \ & \ & \ & \ & \ & to me & \ & \ \\
\hline

& \ & \ & \ & \ & \ & \ & \ & \ & \ & \ & \ & \ & \ & \ \\
ON & $4.61\cdot 10^{11}$ & 18 & 11 & 3 & $A_2(7)$ & $^{\textnormal{l1})}$p.422 & 48 & $p\neq 7$ & d) & 18 & l2) & 10944 & $p=1,$ linear & $^{\textnormal{l1})}$p.461 \\
& \ & \ & \ & \ & $31\cdot 15$ & $^{\textnormal{c})}$p.225 & 15 & $p\neq 31$ & f) & \ & \ & 342 & $p=1,$ proj & $^{\textnormal{l1})}$p.468 \\
\hline

& \ & \ & \ & \ & \ & \ & \ & \ & \ & \ & \ & \ & \ & \ \\
$\cdot 1$ & $4.16\cdot 10^{18}$ & 49 & 24 & 2 & $2^{11}\cdot M_{24}$ & $^{\textnormal{c})}$p.217 & 24 & $p\neq 2$ & g1) & 24 & \ & 24 & $p\neq 2,$ proj & \ \\
& \ & \ & \ & \ & $3^6\cdot(2\cdot M_{12})$ & \ & 24 & $p\neq 3$ & g2) & \ & \ & 24 & $p=2,$ linear & \ \\
\hline

& \ & \ & \ & \ & \ & \ & \ & \ & \ & \ & \ & \ & \ & \ \\
$\cdot 2$ & $4.23\cdot 10^{13}$ & 24 & 14 & 1 & $M_{23}$ & $^{\textnormal{c})}$pp.216, & 22 & $p\neq 2,23$ & \ & 20 & \ & 23 & $p=1$ & \ \\
& \ & \ & \ & \ & $N(5^{1+9})$ & 217 & 20 & $p\neq 5$ & \ & \ & \ & 22 & $p=2$ & \ \\
\hline

& \ & \ & \ & \ & \ & \ & \ & \ & \ & \ & \ & \ & \ & \ \\
$\cdot 3$ & $4.96\cdot 10^{11}$ & 18 & 11 & 1 & $^{2}A_2(25)$ & $^{\textnormal{c})}$p.216 & 20 & $p\neq 5$ & d) & 20 & \ & 23 & $p=1$ & $^{\textnormal{G})}[34]$ \\
& \ & \ & \ & \ & $M_{23}$ & \ & 21 & $p\neq 2$ & a) & \ & \ & \ & \ & \ \\
\hline

& \ & \ & \ & \ & \ & \ & \ & \ & \ & \ & \ & \ & \ & \ \\
$M(22)$ & $6.46\cdot 10^{13}$ & 25 & 15 & 6 & $B_3(3)$ & $^{\textnormal{c})}$p.218 & 27 & $p\neq 3$ & d) & \ & \ & 78 & $p=1,$ linear & $^{\textnormal{c})}[67]$ \\
& \ & \ & \ & \ & $2^{10}\cdot M_{22}$ & \ & \ & $p\neq 2$ & \ & \ & \ & \ & \ & \ \\
\hline

& \ & \ & \ & \ & \ & \ & \ & \ & \ & \ & \ & \ & \ & \ \\
$M(23)$ & $4.09\cdot 10^{18}$ & 49 & 24 & 1 & $D_4(3)$ & $^{\textnormal{c})}$p.219 & 234 & $p\neq 3$ & d) & \ & \ & 782 & $p=1$ & $^{\textnormal{c})}[66]$ \\
& \ & \ & \ & \ & $2^{11}\cdot M_{23}$ & \ & \ & $p\neq 2$ & \ & \ & \ & \ & \ & \ \\
\hline

& \ & \ & \ & \ & \ & \ & \ & \ & \ & \ & \ & \ & \ & \ \\
$\mathcal{D}M(24)$ & $7.38\cdot 10^{22}$ & 83 & 37 & 3 & $D_4(3)$ & $^{\textnormal{c})}$p.219 & 234 & $p\neq 3$ & d) & 234 & \ & \ & \ & \ \\
& \ & \ & \ & \ & $2^{12}\cdot M_{24}$ & \ & 759 & $p\neq 2$ & p) & \ & \ & \ & \ & \ \\
\hline

& \ & \ & \ & \ & \ & \ & \ & \ & \ & \ & \ & \ & \ & \ \\
$F_5$ & $2.73\cdot 10^{14}$ & 27 & 16 & 1 & $^{2}A_2(64)$ & $^{\textnormal{c})}$p.226 & 56 & $p\neq 2$ & d) & 56 & \ & 133 & $p=1$ & $^{\textnormal{c})}[54]$ \\
& \ & \ & \ & \ & $N(5^{1+4})$ & \ & 100 & $p\neq 5$ & e) & \ & \ & \ & \ & \ \\
\hline

& \ & \ & \ & \ & \ & \ & \ & \ & \ & \ & \ & \ & \ & \ \\
$F_3$ & $9.07\cdot 10^{16}$ & 39 & 21 & 1 & $2^{(1+8)}\cdot \textnormal{Alt}_q$ & $^{\textnormal{c})}$p.225 & 128 & $p\neq 2$ & n) & 80 & \ & 248 & $p=1$ & $^{\textnormal{c})}[131]$ \\
& \ & \ & \ & \ & $3^5\cdot A_1(q)$ & $^{\textnormal{o})}$p.67 & 80 & $p\neq 3$ & q) & \ & \ & \ & \ & \ \\
\hline

& \ & \ & \ & \ & \ & \ & \ & \ & \ & \ & \ & \ & \ & \ \\
$F_2$ & $4.155\cdot 10^{33}$ & 260 & 89 & 2 & $^{2}E_6(4)$ & $^{\textnormal{c})}$p.219 & 1536 & $p\neq 2$ & d) & 594 & \ & 196 883 & $p=2,3,$ proj & m) \\
& \ & \ & \ & \ & $N(3^{(1+8)})$ & $^{\textnormal{o})}$p.68 & 594 & $p\neq 3$ & m) & \ & \ & 196 882 & $p=2,$ linear & \ \\
\hline

& \ & \ & \ & \ & \ & \ & \ & \ & \ & \ & \ & \ & \ & \ \\
$F_1$ & $8.08\cdot 10^{53}$ & 1472 & 343 & 1 & $2^{(1+24)}$ & $^{\textnormal{k})}$p.116 & 4096 & $p\neq 2$ & e) & 2132 & \ & 196 883 & \ & \ \\
& \ & \ & \ & \ & $3^8\cdot$ $^{2}D_4(9)$ & \ & 2132 & $p\neq 3$ & i) & \ & \ & \ & \ & \ \\
\hline

\end{tabular}
\end{center}

\end{landscape}

\newpage
\normalsize
\pagestyle{plain}

\subsection{}Reasons.\\

(a) For the groups $M_i, i=11, 12, 22, 23, 24,$ the exact values of degrees of their faithful irreducible (linear) representations in all characteristics are given in 
G. James, \cite{17}.\\

(b) From the Table T2.7 (groups of small degree) we can obtain the exact value of $\mathcal{C}d(H)$ if it is $\leq 5$. Otherwise we have, of course, $\mathcal{C}d(H)\geq 6.$\\

(c) S. A. Syskin, \cite{Sy}.\\

(d) For a group of Lie type we obtain our estimates from V. Landazuri and G. Seits [\cite{18}, p. 419]. Note that the ones we need are, generally, listed as exceptions.\\

(e) For an extraspecial group ($E$ of order) $q^{1+2\,m}, m\geq 1,$ we have the estimate $q^m$ if $q\neq p$ (see. (8.1)). However if $N_H(q^{1+2\,m})=:N(q^{1+2\,m})$ acts transitively on the center $q$ of $q^{1+2\,m}$ then the estimate becomes $(q-1)\,q^m$ because $N(q^{1+2\,m})$ permutes central characters of $q^{1+2\,m}$ and our claim follows from Clifford theory. N. B. e) can be tricky to use if the Schur multiplier of $H$ is divisible by $q$: in this case one has to trace what happens with $q^{1+2\,m}$ after a central extension.\\

(f) If $H$ contains a subgroup $R$ which is the middle term of $1\rightarrow E\rightarrow R\rightarrow M\rightarrow 1$ with $E=q^m$ (i.e. $E\simeq (\mathbb{Z}/q)^m$) with $p\neq q$ and $q$ a prime which is prime to the Schur multiplier then by Clifford theory the faithful irreducible projective representations of $R$ have degree not less than the minimal length of an orbit of $M$ on the characters of $E$, i.e. on $E^v:=Hom(E, \mathbb{F}_q)$.

Particular cases are: 

1) $q\cdot r(:=(\mathbb{Z}/q)\times(\mathbb{Z}/r))$ then the orbits have length $r$,

2) $M$ acts transitively on $E^v$ (for example $M\supseteq SL_m(\mathbb{F}_q)$), and 

3) $M$ is isomorphic to a Mathieu group $M_i$ in which case non-trivial orbits have length $\geq i$.\\

(g1) We apply f) and h). Since $M_{24}$ does not have subgroups of index $\leq 23$, each orbit of $M_{24}$ on $2^{11}$ has length 1 or $\geq 24$. since by a) $M_{24}$ is irreducible on $2^{11}$ it follows by j) that for $\cdot 1$ (resp. for $\cdot 0$) orbits of length 1 all lie in the center whence our claim.\\

(g2) For $3^6\cdot (2\cdot M_{12})$ we again have from Table T2.7 that $2\cdot M_{12}$ is irreducible on $E:=3^6=\mathbb{F}_3^6.$ Then $2\,M_{12}$ acts on $\mathbb{P}(E)\simeq \mathbb{P}^5(\mathbb{F}_3)$ and, one easily sees, has only orbits of length $\geq 12$. Since the center of $2\,M_{12}$ inverts the elements of $3^6$ it follows from the above the $2\,M_{12}$ has on $3^6$ only orbits of length $\geq 2\cdot 12=24$.\\

(h) We need several times the following argument: Suppose $q^m\cdot M$ is a subgroup of $H$ with $M$ irreducible on $q^m$ and $m$ odd. Suppose $\varphi: \tilde{H}\rightarrow H$ is a cover with central kernel $\mathbb{Z}/q$. Then the preimage of $q^m$ in $H$ is isomorphic to $q^{m+1}$ (although the representation of $\varphi^{-1}(M)$ on $q^{m+1}$ may not be completely reducible). Indeed, if the preimage $S:=\varphi^{-1}(q^m)$ is not commutative then the commutator map defines an alternating bilinear form on $q^m$ with values in $\ker{\varphi}$. Since $m$ is odd, the form must be degenerate and since it is (evidently) $M$-invariant and $M$ is irreducible, this form is zero. Thus $S$ is commutative. The map $x\rightarrow x^v$ on $S$ is also invariant under $\varphi^{-1}(M)$ and, as above, must be trivial, whence our claim.\\

(i) We consider $E\times \Omega_8^{-}(\mathbb{F}_3)$ where $E\simeq \mathbb{F}_3^8$ and $\Omega_8^{-}(\mathbb{F}_3)=\mathcal{D}(SO_8^{-}(\mathbb{F}_3)).$ We have $E\simeq E$ as $\Omega_8^{-}(\mathbb(F)_3$-module since it possesses an invariant non-degenerate quadratic and, therefore, bilinear form. One easily derives from Witt's theorem that $\Omega_8^{-}(\mathbb{F}_3)$ acts transitively on the vectors of the same length; possible length are 0, 1, -1. The stabilizers of the corresponding vectors are: the commutator of a parabolic subgroup $\mathbb{F}_3^8\cdot \Omega_6^{-}(\mathbb{F}_3)$, and (for both lengths $\neq 1$) $\Omega_7(\mathbb{F}_3)$. This gives lengths of orbits of $\Omega_8^{-}(\mathbb{F}_3)$ as 1, 2132 (isotropic), 2214 (for both anisotropic)  whence our estimate in view of f).\\

(j) If $H$ is an irreducible subgroup of $GL_u(k)$ with $H$ perfect with center of order $m$ then $m/n$. This is evident.\\

(k) R. L. Griess Jr., \cite{Gr}.\\

(l1) \cite{Nan}.\\

(l2) We have to show that $\mathcal{C}d(H)\geq 18$. We already have that $\mathcal{C}d(H)\geq 15$ and by use of $A_2(7)$ it is sufficient to consider the case $p=7$. So let $\varphi: \tilde{H}\rightarrow GL_n(k)$ be a centrally faithful irreducible representation of the universal cover $\tilde{H}$ of $H$. Let $C(\simeq \mathbb{Z}/3)$ be the center of $\tilde{H}$. $\tilde{H}$ contains (see D. Gorenstein and R. Lyons [\cite{14}, p. 61]) subgroup $E\simeq 3^{1+4}$ whose center is $C$. If $\varphi(C)\neq 1$ then it follows from (8.1) and from complete reducibility of $\varphi|E$ in characteristic 7 that $9|n$. Since $n\geq 15$ it implies that $n\geq 18$.

Thus it suffices to consider the case when $\ker{\varphi}=C$. In this case the Sylow 3-subgroup $\varphi(E)$ of $\delta(\tilde{H})\simeq H$ is isomorphic to $(\mathbb{Z}/3)^4$ and is, see D. Gorenstein, R. Lyons [\cite{14}, p. 61], a self dual $N_H(\varphi(E))/Z_H(\varphi(E))$-module. The structure of $N_H(\varphi(E))$ is given in l1) p. 422. One sees from that that $N_H(\varphi(E))$ is transitive on $\varphi(E)=\Id_n$ whence $n\geq 80$ in our case.\\

(m) R. L. Griess Jr., personal communication.\\

(n) (R. L. Griess) Since $2^{1+8}\cdot \textnormal{Alt}_9$ occuring in $F_3$ is the unique twisted holomorph, the estimate is $2^4$ (from $2^{1+8}$) times the degree of the smallest irreducible projective (and non-linear) representation of $A_9$. This latter is $\geq 8$ by (3.1)(b). Thus the minimal degree in characteristic $\neq 2$ is $\geq 2^4\cdot 8=128$.\\

(o) D. Gorenstein, R. Lyons, ?.\\

(p) By p. 245 of \cite{Con} we see that $\mathcal{D}M(24)$ contains a subgroup $A\cdot M_{24}, A\simeq 2^{12},$ (non-split) with the action of $M_{24}$ on $2^{12}$ being dual to that of $M_{24}$ on the Golay code (which is just another group $2^{12}$). The orbits of $M_{24}$ on the Golay code are of length 1, 759, 2576, 759, 1. Thus, if char $k\neq 2$, the shortest non-trivial orbit of $M_{24}$ on $\Hom{(A,k^{*})}$ is of length 759 whence our estimate by f).\\

(q) By o) p. 67 $F_3$ contains a subgroup $R\simeq 3^4\cdot SL_2(\mathbb{F}_9)$ (which is contained in $N(3\,C)$ in the notation of that paper). The action of $S:=SL_2(\mathbb{F}_9)$ on $E:=3^4$ can not be trivial. Thus we are dealing with a representation of $S$ over $\mathbb{F}_3$. This representation can not be reducible, for all homomorphisms of $SL_2(\mathbb{F}_9)$ into $SL_2(\mathbb{F}_3)$ or $SL_3(\mathbb{F}_3)$ are trivial. Thus it is irreducible. Let $\varphi$ be its extension to $\bar{\mathbb{F}}_3$. If $\varphi$ is reducible then it is a sum of two conjugate 2-dimensional representations of $SL_2(\mathbb{F}_9)$, or in the other words our representation is the natural 2-dimensional representation of $SL_2(\mathbb{F}_9)$ on $\mathbb{F}^2_9(\simeq E)$. In this case $S$ acts transitively on $E-\{1\}$ (and also on $E-\{1\}$ where $E$ is the set of characters of $E$). Thus in this case, by f), we have estimate $9^2-1=80$ on dimensions of representations of $F_3$.

If $\varphi$ is irreducible then $\varphi\simeq\varphi_1\otimes(\varphi_1\ \mathbb{F}r)$ where $\varphi_1$ is the natural 2-dimensional representation of $SL_2(\mathbb{F}_9)$. However in this case the center $C\simeq \mathbb{Z}/2$ of $S$ acts trivially on $E$. But this is impossible because the centralizer of the only (conjugacy class) element of order 3 in $F_3$ does not contains $3\times 3^4\cdot SL_2(\mathbb{F}_9)$. \\

(r1) See comment e) to Table T6.3.\\

(r2) Since the primitive cube root of 1 is not in $\mathbb{F}_5$, the element of order 3 fixes no line of $5^2\simeq \mathbb{F}^2_5$. Since the normalizer of every line permutes transitively the non-zero elements of this line (by o) p. 56) we see that $4\times\Sym_3$ has two orbits of length 12 and one of length 1 on $5^2$ (and on its dual) whence by f) the estimate.\\

\subsection{}Once the table is complete we see that $|H|\leq f(n)$ if the "adjusted estimate" is $\geq \min{n}$ for $H$. If this is not the case then $|H|\leq y_{a_1, a_2, H}(n)$ where $a_2:=\min{n}, a_1:=$"adjusted estimate" or an estimate for one or both of the $H_i$, combined with the corresponding restrictions on $p$. Now inspection of Table T7.2 completes the proof of (7.1).

\newpage
\section{\textbf{Estimate for extraspecial groups.}} 

Let $k$ be an algebraically closed field of characteristic exponent $p:=p(k)$. If $q\neq p$ is a prime then every extraspecial $q$-group $E$ of order $q^{1+2\,a}$ has faithful irreducible representations over $k$ and they all have degree $q^a$ (see D. Gorenstein [\cite{13}, ?]). Let $\Aut_c{E}$ denote the group of automorphisms of $E$ which are trivial on the center of $E$. We have $\Aut_c{E}\supseteq \textnormal{Inn}\,E$ and $\textnormal{Inn}\,E\simeq (\mathbb{Z}/q)^{2\,a}$. It is known (see e.g. B. Huppert [III.13.7 and III.13.8]) that $\Aut_c{E}/\textnormal{Inn}\,E\simeq Sp_{2\,a}(\mathbb{F}_q)$ if $q\neq 2$ and $\simeq O_{2\,a}(\mathbb{F}_2)$ if $q=2$; here $\epsilon = 0$ or 1 and $O_{2\,a}^{\epsilon}$ is the orthogonal group of a quadratic form on $\mathbb{F}_2^{2\,a}$ of Arf invariant $\epsilon$. Recall that $|Sp_{2\,a}(\mathbb{F}_q)|\leq q^{2\,a^2+a}$ (see (4.4.1)) and $|O_{2\,a}^{\epsilon}(\mathbb{F}_2)|=2\cdot |SO_{2\,a}^{\epsilon}(\mathbb{F}_2)|\leq 2\cdot 2^{2\,a^2-a}$ if $a> 1$ (see (4.4.1)) and take into account lower-dimensional isomorphisms, (see R. Steinberg [\cite{?}, Theorem 37]). We have $|O_2^{\epsilon}(\mathbb{F}_2)|\leq 6$ (see J. Dieudonne [\cite{10}, II.10.3]).

For $d\in \mathbb{N}$ set 
$$
N(d, q) :=
\begin{cases}
d^{2\,\log_3{d}+3} & \text{if $q\neq 2$} \\
2\,d^{2\,\log{d}+1} & \text{if $q=2, d>2$}  \\
24 & \text{if $q=2, d\leq 2$}\\
\end{cases}
$$ 

\prop{}
Let $\psi: E\rightarrow GL_d(k)$ be a faithful irreducible representation. Then \\

(a) $d=q^a$,\\

(b) $|\Aut_c{E}|\leq N(d,q)$,\\

(c) $|\Aut{E}|\leq (q-1)\,N(d,q)$,\\

(d) $|N_{GL_d(k)}(\psi(E))/Z_{GL_d(k)}(\psi(E))|\leq |\Aut_c{E}|$.\\
\eprop

\prf{} 
We have already given references for (a). Thus $a=\log_q{d}$. Note that, since $\psi$ is irreducible the center $C$ of $E$ acts on $k^d$ via a fixed character (with values in $\mu_q\subset k^{*}$). Therefore for $n\in N_{GL_d(k)}(\psi(E))$ we must have that $n\,\psi(C)n^{-1}=\psi(C)$ and then that $n\,\psi(c)n^{-1}=\psi(c)$ for any $c\in C$. Thus $N_{GL_d(k)}(\psi(E))$ induces on $E$ automorphisms from $\Aut_c{E}$. Thus (d) follows from (b). Since $\Aut{E}/\Aut_c{E}\simeq \Aut{C}\simeq \mathbb{Z}/{(q-1)}$, (c) also follows from (b). If $q\neq 2$ then $q\geq 3$ so that $a\leq \log_3{d}$. Therefore 
$$
|\Aut_c{E}|=q^{2\,a}\,|Sp_{2\,a}(\mathbb{F}_q)|\leq q^{2\,a}\,q^{2\,a^2+a}=(q^a)^{2\,a+3}=d^{2\,a+3}\leq d^{2\,\log_3{d}+3}=N(d,q).
$$
If $q=2$ then $a=\log{d}$ whence if $d> 2$ we have 
$$
|\Aut_c{E}|\leq 2^{2\,a}\cdot 2\cdot 2^{2\,a^2-a}=2\cdot 2^{2\,a^2+a}=2\cdot (2^a)^{2\,a+1}=2\,d^{2\,\log{d}+1}=N(d, 2) \textnormal{for}\ d> 2. 
$$ 
Finally if $q=2, a=1,$ then $d=2$ and $|\Aut_c{E}|\leq 2^2\cdot 6=24=N(2, 2).$
\epr

\cor{}
In the assumptions of (8.1) we have\\

\tiny
\begin{center}
\begin{tabular}{|c|c|c|c|c|c|c|c|c|c|c|c|c|}
\hline

& \ & \ & \ & \ & \ & \ & \ & \ & \ & \ & \ & \ \\
$n$ & $2$ & $3$ & $4$ & $5$ & $7$ & $8$ & $9$ & $11$ & $13$ & $16$ & $17$ & $19$ \\
\hline

& \ & \ & \ & \ & \ & \ & \ & \ & \ & \ & \ & \ \\
$|\Aut_c{E}|\leq$ & $24$ & $216$ & $1920$ & $3000$ & $16,464$ & $3.3\cdot 10^6$ & $4.2\cdot 10^8$ & $1.6\cdot 10^5$ & $3.7\cdot 10^5$ & $1\cdot 10^{11}$ & $1.4\cdot 10^6$ & $2.6\cdot 10^6$ \\
\hline

& \ & \ & \ & \ & \ & \ & \ & \ & \ & \ & \ & \ \\
$|\Aut{E}|\leq$ & $24$ & $432$ & $1920$ & $12000$ & $98,784$ & $3.3\cdot 10^6$ & $1.26\cdot 10^7$ & $1.6\cdot 10^6$ & $4.4\cdot 10^6$ & $10^{11}$ & $2.2\cdot 10^7$ & $4.7\cdot 10^7$ \\
\hline

& \ & \ & \ & \ & \ & \ & \ & \ & \ & \ & \ & \ \\
$2^{2\,\log{n}+1}$ & $16$ & $195$ & $2048$ & $17616$ & $778,230$ & $4.2\cdot 10^6$ & $2\cdot 10^7$ & $2.5\cdot 10^8$ & $4.5\cdot 10^9$ & $1.3\cdot 10^{11}$ & $3.9\cdot 10^{11}$ & $2.8\cdot 10^{12}$ \\
\hline

\end{tabular}
\end{center}
\ \\

\normalsize
Moreover, $|\Aut{E}|\leq 2\,n^{2\,\log{n}+1}$ for $n\geq 4$.
\ecor

\prf{}
If $n=q^a$ is a power of an odd prime then $|\Aut_c{E}|=n^2\,|C_a(q)|$ and $|\Aut{E}|=(q-1)\,|\Aut{E}|$ whence by direct computation of $|C_a(q)|$ the expressions for the above $n$. If $n=2^a$ then $|\Aut_c{E}|\leq n^2\cdot \max_{\pm}{|O_{2\,a}^{\pm}(2)|}$ and we use then expressions for $|O_{2\,a}^{\pm}(2)|$ from, say, J. Dieudonne [\cite{10}, $\S$II.10].

The estimate $n^{2\,\log_3{n}+4} < 2\,n^{2\,\log{n}+1}$ holds for $n>12$ and since the estimate $2\,n^{2\,\log{n}+1}$ applies by the table above to $4\leq n\leq 12, n$ a prime power, we obtain the concluding statement of (8.2).   
\epr

\cor{}
In the assumptions of (8.1) assume also that $q\neq 2$. Then
$$
|\Aut{E}|\leq 2\,d^{2\,\log_3{d}+1}.
$$
\ecor

\prf{}
The claim holds for $q=3$. Assume $q>3$. Let $d=q^a$ (by (8.1)(a)). Then from the proof of (8.1) we have $|\Aut{E}|=(q-1)\,|\Aut_c{E}|$ and $|\Aut_c{E}|\leq d^{2\,\log_q{d}+3}.$ We have thus to check whether $F(d):=2\,d^{2\,\log_3{d}+3}/{(q-1)\,d^{2\,\log_q{d}+3}}$ is $\geq 1$. We have 
$$
\ln{F(d)}=\ln{2}-\ln{(q-1)}+2\,(1/\ln{3}-1/\ln{q})\,\ln^2{d}\geq
$$
$$
\geq \ln{2}-\ln{q}+2\,(1/\ln{3}-1/\ln{q})\,a^2\,\ln^2{q}=
$$
$$
=\ln{2}+\frac{2\,a^2\,\ln^2{q}}{\ln{3}}-(2\,a^2+1)\,\ln{q}\geq
$$
$$
\geq 0.69+ 1.82\,a^2\,\ln^2{q}-(2\,a^2+1)\,\ln{q}.
$$

Thus 
$$
\ln{F(d)}> 1.82\,a^2\,\ln{q}\,\left(\ln{q}-\left(2+\frac{1}{a^2}\right)\,0.55\right)\geq 1.82\,a^2\,\ln{q}\,(\ln{q}-1.65).
$$

This latter expression is $> 0$ if $q\geq 6$. Thus our claim follows for $q\geq 7$. If $q=5$ then $0.69+1.82\,a^2\,\ln^2{q}-(2\,a^2+1)\,\ln{q}\geq 1.49\,a^2-0.9> 0,$ whence (8.3).
\epr

\newpage
\section{\textbf{Estimates for direct products of finite centrally simple groups.}}

Let $k$ be an algebraically closed field of characteristic exponent $p=p(k)$. Let $G_1,\ldots, G_m$ be finite perfect centrally simple universal (i.e., with trivial Schur multiplier) groups. Let $\varphi_i:G_i\rightarrow GL_{n_i}(k), i=1,\ldots, m,$ be non-trivial irreducible representations. Set $G:=\prod_{1\leq i\leq m}{G_i}, \varphi:=\otimes_{1\leq i\leq m}{\varphi_i}, n:=\prod_{1\leq i\leq m}{n_i}.$ Let $C_i$ be the center of $G_i, i=1,\ldots, m.$\\ 

We introduce the following subsets of the set of the $i, 1\leq i\leq m$.\\

$I_{\textnormal{Alt}}$, the set of $i$ such that $\varphi_i(G_i)$ is isomorphic to $\textnormal{Alt}_{a_i}, a_i\geq 10$;\\

$I_{\textnormal{Lie}, p}$, the set of $i$ such that $G_i$ is isomorphic to a group of Lie $p$-type (see Table T6.3 for exceptional isomorphisms);\\

$I_{\textnormal{Lie}, p'}$, the set of $i$ such that $G_i$ is isomorphic to a group of Lie $p'$-type but not isomorphic to any group of Lie $p$-type (see Table T6.3 for exceptional isomorphisms) and neither isomorphic to $^{2}A_3(9)$ or $D_4(2)$.\\

$I_{\textnormal{extrz-spor}}$, the set of $i$ such that $G_i$ is centrally isomorphic to Suz,$\cdot 1, \cdot 2$, if $p\neq 3,$ to $^{2}A_3(9)$, and if $p\neq 2$, to $D_4(2)$.\\

$I_{\textnormal{rest}}$, the set of the remaining indices; i.e. $\varphi_i(G_i)$ for $i\in I_{\textnormal{rest}}$ is isomorphic to $\tilde{\textnormal{Alt}}_{a_i}$ for $a_i\geq 10$, or centrally isomorphic to $\textnormal{Alt}_7, \textnormal{Alt}_9,$ or to one of the 23 sporadic groups not included in $I_{\textnormal{extra-spor}}.$\\

These sets are pairwise disjoint.\\

Set $H:=\prod_{i\notin I_{\textnormal{Lie}, p}}{G_i}, L:=\prod_{i\in I_{\textnormal{Lie}, p}}{G_i}, A:=\prod_{i\in I_{\textnormal{Alt}}}{G_i}, C=\textnormal{center of } H,$ and for $R\subseteq GL_n(k)$ set
$$
NZ(R):=N_{GL_n(k)}(R)/{Z_{GL_n(k)}(R)}
$$
$$
\overline{NZ}(R):=N_{GL_n(k)}(R)/{R\cdot Z_{GL_n(k)}(R)}.
$$

Thus $NZ(R)$ "describes" the part of the automorphism group of $R$ "realized" in $N_{GL_n(k)}(R)$ and $\overline{NZ}(R)$ the corresponding past of outer automorphism group. 

Set $f(x):= (2\,x+1)^{2\,\log_3{(2\,x+1)}+1}.$

\th{}
(a)
$$
|H/C|\leq 
\begin{cases}
f(n) & \text{if $n<10, n\neq 6, 8$} \\
1.58\,f(6) & \text{if $n=6$} \\
4.62\,f(8) & \text{if $n=8$} \\
(n+2)! & \text{if $n\geq 10$ }
\end{cases}
$$
$$
|\Aut{H}|\leq
\begin{cases}
n\,f(n) & \text{if $n=2, 3, 5, 7, 9, 10, 11$} \\
4.1\,f(4) & \text{if $n=4$} \\
12.61\,f(6) & \text{if $n=6$} \\
27.69\,f(8) & \text{if $n=8$} \\
231\,f(12) & \text{if $n=12$} \\
(n+2)! & \text{if $n>12$} \\
\end{cases}
$$ 
\ \ \ \ \ \ \ \ \ \ \ \ \ \ \ \ \ \ \ \ \ \ \ \ \ \ \ \ $
|\textnormal{Out} H|\leq [\log{n}]!\,n^2 \ \ \ \text{if $n\geq 2$}\\
$

(b) if $I_{\textnormal{Alt}}=\varnothing$ and \\

\ \ \ \ (i) $|I_{\textnormal{extra-spor}}|\neq 1$,\\

or (ii) $\prod_{1\leq i\leq m, i\notin I_{\textnormal{extra-spor}}}{n_i}\geq 4,$\\

or (iii) $m\geq 3,$\\

then $|H/C|\leq f(n),$\\

\ \ \ \ \ \ \ $|\Aut{H}|\leq 
\begin{cases}
n\,f(n) & \text{if $n\neq 4$} \\
4.1\,f(4) & \text{if $n=4$} \\
\end{cases}
$

\ \ \ \ \ \ \ $|\textnormal{Out}\,H|\leq [\log{n}]!\,n^2$\\

(c) if $m\leq 2, I_{\textnormal{Alt}}=\varnothing, |I_{\textnormal{extra-spor}}|= 1$ (assume $I_{\textnormal{extra-spor}}=\{1\}$), and $n/{n_1}=1, 2,$ or 3 then 
$$
|H/\textnormal{center}|\leq
\begin{cases}
F(G_1, n) & \text{if $m=1$} \\
60\,F(\cdot 1, n) & \text{if $m=2, G_1\simeq \cdot 0, 48\leq n \leq 72$} \\
f(n) & \text{otherwise} \\
\end{cases}
$$
$$
|\Aut{H}|\leq
\begin{cases}
F(G_1, n) & \text{if $m=1$} \\
3.05\cdot 48\,f(48) & \text{if $m=2, G_1/{C_1}\simeq \cdot 1, n=48$} \\
1.43\cdot 50\,f(50) & \text{if $m=2, G_1/{C_1}\simeq \cdot 1, n=50$} \\
n\,f(n) & \text{in the remaining cases} \\
\end{cases}
$$
\ \ \ \ \ \ \ \ \ \ \ \ \ \ \ \ \ \ \ \ \ \ $|\textnormal{Out}\,H|\leq 32$\\

where $F(G_1, n)=F(G_1/\textnormal{center}, n)$ is defined just before (A6).\\

The normalizer of $\varphi(G)$ in $GL_n(k)$ permutes the (isomorphic linear) groups $\varphi(G_i)$. Let $\bar{\varphi}: \overline{NZ}(\varphi(G))\rightarrow \Sym_m$ be the corresponding homomorphism.
\eth

\th{}\ \\

(a) $\ker{\bar{\varphi}}$ is solvable with $\mathcal{D}^3(\ker{\bar{\varphi}})=\{1\}$;\\

(b) $|\ker{\bar{\varphi}}|\leq n^2$ and, therefore,\\

(c) $|\overline{NZ}(\varphi(G))|\leq [\log{n}]!\cdot n^2$.\\

Remark. Of course, if $G$ does not contain $D_4(m)$ then $\mathcal{D}^2(\ker{\bar{\varphi}})=\{1\}$.
\eth

\lem{} 
If $|G_i/{C_i}|\leq f(n_i)$ for $i=1,\ldots, m$ and $I_{\textnormal{Lie}, p}=\varnothing$ then\\

(a) $|G/\textnormal{center}|\leq f(n)$ \\

(b) $|\Aut{G}|\leq
\begin{cases}
n\,f(n) & \text{if $n\neq 4$} \\
4.1\,f(4) & \text{if $n=4$}\\
\end{cases}$
\ \\

(c) $|\textnormal{Out}\,G|\leq [\log{n}]!\,n^2$\\

(d) $|\prod{\textnormal{Out}\,G_i}|\leq n^2$.

\elem

\prf{}
We have $|G/\textnormal{center}|=\prod{|G_i/{C_i}|}\leq \prod{f(n_i)}\leq f\left(\prod{n_i}\right)=f(n)$ (the last inequality by (A1)(a)). This proves (a).

Now we use (6.1)(b), (c). Note that the exceptions in (6.1)(b) are all, but $A_2(4)$ for $n=4$, excluded: $\textnormal{Alt}_6$ and $A_1(8)$ do not have representations of dimension 2 unless they are of Lie $p$-type (see column 7 in Table T6.3) and $^{2}A_3(9)$ is excluded for $n=6$ by $|^{3}A_3(9)|> f(6)$, similarly $D_4(2).$ Thus we have $|\textnormal{Out}\,G_i|\leq n^2_i$ and $|\Aut{G_i}|\leq n_i\,f(n_i)$ unless $n_i=4, G_i\simeq \tilde{A}_4(2)$ when $|\Aut{G_i}|=4.1\,f(4)$, if $G_i$ is of Lie $p'$-type. For the remaining centrally simple groups: alternating of degree $\neq 6$ (recall $\textnormal{Alt}_6\simeq \mathcal{D}B_2(2)$) and sporadic, we have $|\textnormal{Out}\,G_i|\leq 2$. Thus $|\textnormal{Out}\,G_i|\leq n^2_i.$ This implies (d).

To prove (b) and (c) let us split the set $J=\{1, \ldots, m\}$ of indices into the subsets $J_1,\ldots, J_r$ such that $G_i\simeq G_j$ if and only if $i, j\in J_s$ for some $s$. Then $\Aut{G}=\Aut{(\prod{G_i})}=\prod_{1\leq i\leq r}{\left(\left(\prod_{j\in J_i}{\Aut{G_j}}\right)\rtimes \Sym{J_i}\right)}$ where $\Sym{J_i}$ permutes the isomorphic $G_j, j\in J_i.$ Let $t_i:=|J_i|.$ Let $\tilde{n}_i:=\min{\{n_j, j\in J_i\}}$. Assume for definiteness that if some $G_j\simeq \tilde{A}_2(4)$ then $j\in J_1.$ Write $|\Aut{G_j}|\leq a_i\,\tilde{n}_i\,f(\tilde{n}_i)$ where $a_i=1$ unless $i=1, \tilde{n}_1=4,$ and $G_j\simeq \tilde{A}_2(4)$ for $j\in J_1$, and $a_i=1.025$ otherwise, for $i=1,\ldots, r$. If $n_j=2$ then $G_j\simeq \tilde{\textnormal{Alt}}_5$ (see Table T2.7), i.e., $|G_j/{C_j}|=60, |\Aut{G_j}|=120.$ Thus in notation of (A1) $|\Aut{G_j}|\leq \tilde{n}_i\,\tilde{f}(\tilde{n}_i)$. Then by (A1) we have 
$$
|\Aut{G}|=\prod_{1\leq i\leq r}{\left(\prod_{j\in J_i}{|\Aut{G_j}|}\right)\,\left(t_i!\right)}\leq \left(a_1^{t_1}\,\tilde{n}_1^{t_1}\,\tilde{f}(\tilde{n}_1)^{t_1}\,t_1!\right)\cdot \prod_{i\geq 2}{\tilde{n}_i^{t_i}\,\tilde{f}(\tilde{n}_i)^{t_i}\,t_i!}\leq
$$
$$
\leq a_1^{t_1}\,\prod{\tilde{n}_i^{t_i}\,\left(\prod{\tilde{f}(n_i)^{t_i}\,(t_i!)}\right)}\leq a_1^{t_1}\,n\cdot \prod{\tilde{f}(n_i^{t_i})}\leq a_1^{t_1}\,n\,\tilde{f}\left(\prod{\tilde{n}_i^{t_i}}\right)\leq a_1^{t_1}\,n\,\tilde{f}(n)\leq a_1^{t_1}\,n\,f(n).
$$ 
This proves (b) if $a_1=1.$ If $a_1=1.025$ and $t_1\geq 2$ we use (A1)(d) to get that $a_1^{t_1}\,f(\tilde{n}_1)^{t_1}\leq f(\tilde{n}_1^{t_1})$ and then the argument can be continued (from the middle) as above (but $a_1$-factor will be lost). If $a_1=1.025, t_1=1, m>1$ then the above $a_1\,f(n)\,\prod_{i\geq 2}{\tilde{f}(\tilde{n}_t^{t_i})}\leq a_1\,f(4)\,\tilde{f}(n/{n_1})\leq f(n)$ by (A1)(c). This concludes the proof of (b). 

Part (c) now follows as above:
$$
|\textnormal{Out}\,G|=\prod_{i}{\left(\prod_{j\in J_i}{|\textnormal{Out}\,G_j|}\right)\,t_i!}=\left(\prod_{s}{\textnormal{Out}\,G_s}\right)\,\prod{t_i!}\leq \prod{n_s^2}\cdot\left(\left(\sum{t_i}\right)!\right)\leq n^2\,\left([\log{n}]!\right).
$$
\epr

Set $g(x):=\Gamma(x+3);$ recall that $\Gamma(n+3)=(n+2)!$ for $n\in\mathbb{N}.$

\lem{}
If $G_i\simeq \tilde{\textnormal{Alt}}_{m_i}, m_i\geq 10,$ for $i=1,\ldots, m$ then\\

(a) $|\Aut{G}|\leq (n+2)!$\\

(b) $|\textnormal{Out}\,G|\leq n\cdot ([\log{n}]!)$ 
\elem

\prf{}
By (3.1)(a) we have $m_i\leq n_i+2,$ so that $|G_i/\textnormal{center}|\leq 1/2\,(n_i+2)!$. Since $\textnormal{Alt}_r=\Sym_r$ for $r\geq 7$ (see D. Gorenstein [\cite{13}, p. 304]) we have $|\Aut{G_i}|\leq (n_i+2)!$ Now the same argument as in the proof of (9.3) yields $|\Aut{G}|\leq \prod_{j}{\left(\Aut{G_j}\right)}\cdot \prod_{i}{t_i!}$ where $\sum_{i}{t_i}\leq \log{n}$. Then (A2) gives (a). To get (b) note that $|\textnormal{Out}\,G|=\prod_{j}{\textnormal{Out}\,G_j}\cdot \prod_{i}{t_i!}\leq 2^{\log{n}}\cdot [\log{n}]!$ as claimed. 
\epr

\lem{}
For a group $G_i, i\in I_{\textnormal{rest}},$ we have\\

(a) $|G_i/{C_i}|\leq f(n_i),$ \\

(b) $|\Aut{G_i}|\leq 2\,f(n_i),$ \\

(c) $|\textnormal{Out}\,G_i|\leq 2.$ \\

\elem

\prf{}
When $G_i$ is centrally isomorphic to a sporadic group (a) is contained in (7.1). If $G_i$ is centrally isomorphic to $\textnormal{Alt}_7$ and $\textnormal{Alt}_9$ we see from Table T2.7 (groups of small degree) that $\textnormal{Alt}_7$ has no non-trivial projective representations of dimension 2 and one easily sees that $|\textnormal{Alt}_7|=2520< f(3)$. For $\textnormal{Alt}_9$ we see from the same Table T2.7 that it does not have non-trivial projective representations of dimension $\leq 4$ and then $|\textnormal{Alt}_9|=781,440< f(5)$. 

If $G_i$ is isomorphic to $\tilde{\textnormal{Alt}}_a, a\geq 10,$ then (a) is contained in (3.4). 

Finally, (c) is known, see D. Gorenstein [\cite{13}, p. 304]. 
\epr

\subsection{} \textit{Proof (9.1)(a) and (9.2).} (9.1)(a) follows directly from (9.3), (9.4) and (A3) and (A6) since the estimates for one $G_i$ hold by (3.1), (6.1), (7.1), and (9.5). To prove (9.2) note that, similarly to the proofs of (9.3)(b), (c) and (9.4), we have that $\ker{\bar{\varphi}}$ is a subgroup of $\prod{\overline{NZ}(\varphi_i(G_i))}.$ Then: by (6.1)(c) $|\overline{NZ}(\varphi_i(G_i))|\leq n_i^2$ if $G_i$ is of Lie $p'$-type; by (5.1) $|\overline{NZ}(\varphi_i(G_i))|\leq 2\,n_i\,\log{n_i}\leq n_i^2$ if $G_i$ is of Lie $p$-type;  and by D. Gorenstein [\cite{13}, p. 304] we have $|\overline{NZ}(\varphi_i(G_i))|\leq 2\leq n_i^2$ in the remaining cases. Thus $|\overline{NZ}(\varphi(G))|\leq \left(\prod{n_i^2}\right)\cdot \left(m!\right)=n^2\cdot \left(m!\right)\leq n^2\,\left([\log{n}]!\right)$. This proves (9.2)(c); and (b) was also implicitly proved above. To prove (9.2)(a) recall that $\mathcal{D}^3(\textnormal{Out}\,\tilde{G})$ for every simple group $\tilde{G}$ of Lie type and $|\textnormal{Out}\,\tilde{G}|\leq 4$ in the remaining cases.$\Box$

\subsection{} \textit{Proof of (9.1)(b).} If $|I_{\textnormal{extra-spor}}|>1$ (and $I_{\textnormal{Alt}}=\varnothing$) then (A6), (A1), and (9.5) imply that $|H/C|\leq f(n)$. The same argument as in the proof of (9.3) together with the fact that $|\textnormal{Out}\,G_i|\leq 2 \leq n_i^2$ for sporadic groups, $|\textnormal{Out}$ $^{2}A_3(9)|=8< 36\leq n_i^2$ if $G_i\simeq$ $^{2}\tilde{A}_3(9),$ and $|\textnormal{Out}\,D_4(2)|=6\leq 64\leq n_1^2$ if $G_i\simeq \tilde{D}_4(2)$ gives $|\textnormal{Out}\,H|\leq ([\log{n}]!)\,n^2$. The estimates on $|\Aut{(\ )}|$ in (9.3), (9.5) together with evident estimates $|\Aut{G_i}|\leq F(G_i, n_i)$ for $i\in I_{\textnormal{extra-spor}}$, combined with (A1) and (A6) give that $|\Aut{H}|\leq n\,f(n)$. 

If $I_{\textnormal{extra-spor}}=\varnothing$ then the above proof works again but the reference to (A6) is not needed anymore.

If $|I_{\textnormal{extra-spor}}|=1$ and (9.1)(b)(iii) holds then so does (9.1)(b)(ii). So assume that $I_{\textnormal{extra-spor}}=\{1\}, n/{n_1}\geq 4.$ Then $\Aut{H}=\Aut{(H/{G_1})}\times \Aut{G_1}$. We have by above 
$$|\Aut{(H/{G_1})}|\leq 
\begin{cases}
(n/{n_1})\,f(n/{n_1}) & \text{if $n/{n_1}\neq 4$} \\
4.1\,f(4) & \text{if $n=4\,n_1$}.\\
\end{cases}
$$

Now $|\Aut{G_1}|\leq F(G_1, n)$ and the claim follows from (A6).$\Box$

\subsection{} The proof of (9.1)(c) is also similar. If $m=1$ then we are done in view of (7.1) and of definition of $F(G, n).$ If $m=2$ then $n=n/{n_1}=2$ or 3. If $2\notin I_{\textnormal{Lie,}p}$ then from the Table T2.7 (groups of small degree) $G_2$ is isomorphic to $\tilde{\textnormal{Alt}}_a$ for $a=5, 6,$ or 7, or $\tilde{A}_1(7).$ First, $|\textnormal{Out}\,G_2|\leq 4$ (and $=4$ if $G_2\simeq \tilde{\textnormal{Alt}}_6$) and $|\textnormal{Out}\,G_1|\leq 8$ (and $=8$ if $G_1\simeq$ $^{2}\tilde{A}_3(9)$), whence the estimate on $\textnormal{Out}\,H$.

Now if $n_2=2$ then $|G_2/{C_2}|=60$ and $|H/C|\leq 60\,F(G_1, n_1).$ If $n_2=3$ then $|G_2/{C_2}|=2,520$ and $|H/C|\leq 2,520\,F(G_1, n_1)$.

We now refer to Table TA6 and use notation therefrom. We have $F(G_1, n_1)\geq F(G_1, a_1)=|G_1/{C_1}|$. We see from the Table that $60\,F(G_1, a_1)\leq F(G_1, 2\,a_1)=f(2\,a_1)$ and $2,520\cdot$ $F(G_1, a_1)$ $\leq F(G_1, 3\,a_1)=f(3\,a_1)$ unless $G_1\simeq \cdot 0$. If $G_1/{C_1}\simeq\cdot 1$ we establish by direct calculations that $2,520\,F(G_1, n_1)\leq F(G_1, 3\,n_1)=f(3\,n_1)$ if $n_1\geq 25$ and $60\,F(G_1, n_1)\leq F(G_1, 2\,n_1)\leq f(G_1, 2\,n_1)=f(2\,n_1)$ if $n\geq 31.$ The remaining part of the estimate on $|H/C|$ is then verified by direct calculations again.

If, finally, $2\in I_{\textnormal{Lie,}p}$ then we have to check weaker inequalities $F(G_1, n_1\,n_2)\geq |\Aut{G_1}|$ which, of course, hold if the ones above hold.

The estimate on $|\Aut{H}|$ holds if $m=1$ holds by the definition of $F(G, n)$ and if $m=2$ by (A6) except when $m=2, G_1\simeq\cdot 0$.

We assume now that $m=2, G_1/{C_1}\simeq\cdot 1$. Then $\Aut{G_1}=G_1/{C_1}.$ If $n_2=2$ we see that $120\cdot |\cdot 1|\leq 2\,n_1\,f(2\,n_1)$ if $n_1\geq 26;$ for $n_1=24, 25$ we have values given by (9.1)(c). If $n_2=3$ then $5040\cdot |\cdot 1|\leq 3\,n_1\,f(3\,n_1)$ if $n_1\geq 24$.

That $|\Aut{H}|\leq n\,f(n)$ holds in the remaining cases is easily seen from Table TA6. $\Box$

\newpage
\section{\textbf{Estimates for direct products of centrally simple and extraspecial groups.}}  

Let $k,m,$ the $G_i, \varphi_i, n_i, i=1,\ldots, m, \varphi, H,$ and $L$ have the same meaning as in $\S 9$. Recall $f(n)=(2\,n+1)^{2\,\log_3{(2\,n+1)}+1}.$ Set
$$
f_H(n)=
\begin{cases}
F(G_1, n) & \text{if $H=G_1, G_1/{C_1}\simeq$$^{2}\bar{A}_3(9), D_4(2),$ Suz, $\cdot 1,$ or $\cdot 2$} \\
1.025\,f(4) & \text{if $m=1, H/C\simeq \bar{A}_2(4)$} \\
3.05\,f(48) & \text{if $m=2, G_1/{C_1}\simeq \cdot 1, n=48$} \\
1.43\,f(50) & \text{if $m=2, G_1/{C_1}\simeq \cdot 1, n=50$} \\
f(n) & \text{in the remaining cases.}
\end{cases}
$$  

Let $E_1, \ldots, E_r$ be extraspecial groups of orders $p_1^{1+2\,a_1},\ldots ,p_r^{1+2\,a_r}$ where $p_1< p_2 <\ldots < p_r, p_i\neq p$ for $i=1,\ldots ,r,$ and $a_i\geq 1$ for $i=1,\ldots ,r.$ Let $\varphi_i: E_i\rightarrow GL_{d_i}(k)$ be faithful irreducible representations. By (8.1) we have $d_i=p_i^{a_i}.$ Set $\psi:=\otimes_{1\leq i\leq r}{\psi_i}, d:=\prod_{1\leq i\leq r}{d_i}, E:=\prod_{1\leq i\leq r}{E_i}$ and, if $p_1=2,$ set also $\bar{E}:=\prod_{2\leq i\leq r}{E_i}, \bar{\psi}:=\otimes_{2\leq i\leq r}{\psi_i}, \bar{d}:=\prod_{2\leq i\leq r}{d_i}$.

Set $G:=L\times H\times E, \omega:=(\otimes{\varphi_i})\otimes \psi, n:=(\prod{n_i})\,d, N:=|\Aut{(H\times E)}|, d:=(\log{3}-1)/2\leq 0.2925, \beta:=\log_{24}{|\cdot 1|}-2\,\log{24}\leq 4.32.$

\th{}\ \\

(a) $
|\Aut{(H\times E)}|\leq
\begin{cases}
n\,f(n) & \text{if $n=2, 3, 5, 7, 9, 10, 11$} \\
4.1\,f(4) & \text{if $n=4$} \\
12.61\,f(6) & \text{if $n=6$} \\
27.69\,f(8) & \text{if $n=8$} \\
231\,f(12) & \text{if $n=12$} \\
(n+2)! & \text{if $n> 12$} \\
\end{cases}
$
\ \\

(b) If $I_{\textnormal{Alt}}=\varnothing$ and $|E/\bar{E}|=2^{1+2\,a}$ then $|\Aut{(H\times E)}|\leq 2^{2\,a^2+1}\,n\,F_H(n/{2^a})$. 

\eth   

\cor{}

If $I_{\textnormal{Alt}}=\varnothing, E\neq \{1\}$ then\\

(a) $N\leq n\,f(n)$ if any one of the following holds\\

\ \ \ \ \ \ (i) $p_1\neq 2$,\\

\ \ \ \ \ \ (ii) $p_1=2, 1\leq a_1\leq 5$,\\

\ \ \ \ \ \ (iii) $p_1=2, a_1\geq 2, n\geq 2^{(\alpha+1)\,a_1}$\\

(b) $N\leq 2\,n^{2\,\log{n}+1}$ otherwise. 

\ecor

\cor{}

If $I_{\textnormal{Alt}}=\varnothing$ then \\

(a) $|\Aut{(H\times E)}|\leq \max{\{n\,F_H(n), 2\,n^{2\,\log{n}+1}\}}$, \\

(b) $|\Aut{(H\times E)}|\leq n^{2\,\log{n}+5}$ for $n\geq 2$, \\

(c) $|\Aut{(H\times E)}|\leq 2\,n^{2\,\log{n}+1}$ if $n\geq 39$.

\ecor

\subsection{} \textit{Proof of (10.1).} Note that since the $E_i$ are not isomorhpic and none of them is isomorphic to a quotient of $H$ we have 

\subsubsection{} $|\Aut{(H\times E)}|=|\Aut{H}|\cdot \prod_{1\leq i\leq r}{|\Aut{E_i}|}.$

Since one easily sees that $2\,n^{2\,\log{n}+1}< (n+2)!$ for $n> 12$ and $\leq f(n)$ for $n< 14$ one gets (a) from (A1)(a), (A2)(a), (8.2), and (9.1)(a).

To prove (b) note first that $n^{2\,\log_3{n}+3}\leq f(n)$ by (A7)(a). Therefore by (9.1)(b), (c)
$$
|\Aut{(H\times \bar{E})}|\leq (n/{2^a})\,F_H(n/{2^a}).
$$

Assuming that $a\neq 0$ (otherwise we done) we have $p_1=2$ and $|\Aut{E_1}|=2\cdot(2^a)^{2\,a+1}=2^{2\,a^2+a+1}$ whence (b) follows, in view of (10.4.1).$\Box$

\subsection{}\textit{Proof of (10.2).} We have $I_{\textnormal{Alt}}=\varnothing, E\neq \{1\}$.

\subsubsection{}Assume first that 
$$
|\Aut{H}|\leq \left(\prod{n_i}\right)\,f\left(\prod{n_i}\right).
$$

If $p_1\neq 2$ then by (10.4.1), (8.3), and (A7)(a)
$$
|\Aut{(H\times E)}|\leq \left(\prod{n_i}\right)\,f\left(\prod{n_i}\right)\cdot\prod{2\,d_j^{2\,\log_3{d_j}+3}}\leq
$$
$$
\leq\left(\prod{n_i}\right)\,f\left(\prod{n_i}\right)\,\prod{f(d_j)}\leq n\,f\left(\prod{n_i}\,\prod{d_j}\right)\leq n\,f(n), 
$$

whence (10.2)(a)(i).

If $p_1=2$ but $1\leq a_1\leq 5$ then $24< f(2)$ and $2\,d_1^{2\,\log_3{d_1}+1}\leq d_1\,f(d_1)$ (direct verification). Therefor (10.4.1), (10.1)(b) and (A1) give 
$$
|\Aut{(H\times E)}|\leq n\,f(d_1)\,f(n/{d_1})\leq n\,f(n)
$$

whence (10.2)(a)(ii).

Now let $p_1=2, a_1\geq 2, n\geq 2^{(\alpha+1)\,a_1}.$ Then by (A5) $2\cdot 2^{2\,a_1^2+a_1}\,f(n/{2^{a_1}})\leq f(n)$ whence by (10.4.1) and (10.1)(b) we get (10.2)(a)(iii).

For (b) we note that the "otherwise" condition implies $a_1\geq 6$. Now (b) follows from (A7)(c).

\subsubsection{}Now we consider the remaining case when $|\Aut{H}|\geq \left(\prod{n_i}\right)\,f\left(\prod{n_i}\right)$. It follows then from (9.1)(b), (c) that either $m=1$ and $H=G_1$ is centrally isomorphic to Suz, $\cdot 1, \cdot 2$, $^{2}A_3(9), D_4(2)$ or$A_2(4)$ or $m=2, n=48$ or 50, and $G_1/{C_1}$ or $G_2/{C_2}$ is isomorphic to $\cdot 1$.

Assume first that $\bar{E}\neq \{1\}$. Then by (8.3), (A7)(a): $|\Aut{\bar{E}}|\leq \prod{2\,d_i^{2\,\log_3{d_i}+3}}\leq f(\bar{d})$ whence $|\Aut{(H\times \bar{E})}|\leq |\Aut{H}|\cdot f(\bar{d})$. Using one of (A1)(c) or (A6)(c), (d) we get, therefore that $|\Aut{(H\times \bar{E})}|\leq \bar{d}\left(\prod{n_i}\right)\,f\left(\bar{d}\,\prod{n_i}\right)$ unless $m=1, G_1/{C_1}\simeq\cdot 1$ and $\bar{d}=3$. In this latter case we have $|\Aut{\bar{E}}|=432$ (by (8.2)) and $\Aut{\cdot 1}=\cdot 1$. Then
$$
|\Aut{(H\times \bar{E})}|=432\cdot 4.16\cdot 10^{18}\sim 1.8\cdot 10^{21}\leq (3.24)\,f(3.24)\sim 4\cdot 10^{23}
$$ 
and since $n\,f(n)\geq 72\,f(72)$ if $n\geq 72$ it follows that in all cases (if $\bar{E}\neq \{1\}$) $|\Aut{(H\times \bar{E})}|\leq \bar{d}\left(\prod{n_i}\right)\,f\left(\bar{d}\,\prod{n_i}\right)$.

Now (10.2) follows as in (10.5.1).

\subsubsection{}Now we can assume that $\bar{E}=\{1\}$ so that $E=E_1\simeq 2^{1+2\,a}, a:=a_1\geq 1$. If $a\leq 5$ then $|E|\leq f(d)$ and the argument of (10.5.2) works except when $d=2$. When $d=2, m=1, G_1\simeq A_2(4)$ we still get our claim by (A1)(c). In the remaining cases we have to check that $24\cdot |\Aut{H}|\leq 2\cdot\left(\prod{n_i}\right)\,f\left(2\,\prod{n_i}\right)$. If $m=1$ this is readily seen (except $\cdot 1$) from Table TA6 (compare lines $F(G, 2\,a_1)$ and $24\,|\Aut{G}|$). If $m=1$ and $G_1\simeq \cdot 1$ then one has $|\Aut{(H\times E)}|=10^{20}$ and $48\,f(48)\geq 1.63\cdot 10^{20}$ whence the desired estimate holds in this case as well. The cases $m=2$ and $n=48$ (resp. 50) and $|\Aut{H}|\leq 3.05\cdot 48\,f(48)$ (resp. $\leq 1.43\cdot 50\,f(50)$) are verified directly.

\subsubsection{}We now assume that $a\geq 6$ (and the assumptions of (10.5.2) and (10.5.3) hold). Suppose first that $|\Aut{H}|\leq d\cdot \prod{n_i}\,f\left(\prod{n_i}\right)$ and $n\geq 2^{(1+\alpha)\,a}.$ Then (by (10.4.1) and (A5))
$$
|\Aut{(H\times E)}|\leq 2\,d^{2\,\log{d}+1}\,\left(d\,\prod{n_i}\right)\,f\left(\prod{n_i}\right)\leq \left(d\,\prod{n_i}\right)\,f\left(d\,\prod{n_i}\right)\leq n\,f(n).
$$ 
This gives (10.2)(a)(iii). The condition $|\Aut{H}|\leq d\,\left(\prod{n_i}\right)\,f\left(\prod{n_i}\right), d\geq 64,$ evidently (from expressions in (9.1)(c)) holds if $m=2$. If $m=1, G_1/{C_1}\simeq A_2(4)$ (resp. $^{2}\bar{A}_3(9), D_4(2)$, Suz, $\cdot 1, \cdot 2$) then $n_1\geq 4$ (resp. 6, 8, 12, 24, 20 (by Table T7.2)) and (10.2)(a)(iii) follows from the comparison of lines $64\,a_1\,f(a_1)$ and $|\Aut{G}|$ in Table TA6, except, of course, when $G_1/{C_1}\simeq\cdot 1$.  In this case it holds if $d\geq 4096=2^{12}.$ Thus it is sufficient to check (10.2)(a)(iii) for $m=1, G_1/{C_1}\simeq \cdot 1, d=2^a, 6\leq a\leq 11,$ directly. We have \\

\begin{tabular}{|c|c|c|c|c|c|c|}
\hline

& \ & \ & \ & \ & \ & \ \\
$a$ & $6$ & $7$ & $8$ & $9$ & $10$ & $11$ \\ 
\hline

& \ & \ & \ & \ & \ & \ \\
$|\Aut{E}|=2^{2\,a^2+a+1}$ & $6\cdot 10^{23}$ & $8.1\cdot 10^{31}$ & $1.7\cdot 10^{41}$ & $6\cdot 10^{51}$ & $3.2\cdot 10^{63}$ & $2.9\cdot 10^{76}$ \\
\hline

& \ & \ & \ & \ & \ & \ \\
$|\cdot 1|\cdot|\Aut{E}|$ & $2.5\cdot 10^{42}$ & $3.4\cdot 10^{50}$ & $7.1\cdot 10^{59}$ & $2.5\cdot 10^{70}$ & $1.4\cdot 10^{82}$ & $1.2\cdot 10^{95}$ \\
\hline

& \ & \ & \ & \ & \ & \ \\
$f(2^a\cdot 24)$ & $2.9\cdot 10^{54}$ & $9\cdot 10^{63}$ & $1.5\cdot 10^{74}$ & $1.5\cdot 10^{85}$ & $9\cdot 10^{96}$ & $3\cdot 10^{109}$ \\
\hline

\end{tabular}
\ \\

Thus (10.2)(a)(iii) holds in all cases.

\subsubsection{}It remains to prove (10.2)(b) in the case when the conditions of (10.5.2), (10.5.3) and (10.5.2) hold. So $|E|=2^{1+2\,a}, a\geq 6$. We have to check 
$$
|\Aut{E}|\cdot |\Aut{H}|\leq r^{2\,\log{r}+1},
$$
where $r=2^a\cdot s, s=\prod{n_i}.$ Substituting $|\Aut{E}|\leq 2^{2\,a^2+a+1}$ and taking $\log$ we see that it suffices to check
$$
\log_{10}{|\Aut{H}|}\leq \left(\frac{2}{\log_{10}{2}}\right)\,\log_{10}^2{s}+\left(4\,a+1\right)\,\log_{10}{s}.
$$

Noting that $a\geq 6,$ i.e., $4\,a+1\geq 25$ we easily see that $(4\,a+1)\,\log_{10}{s}$ supplies a power of 10 sufficient to overcome $|\Aut{H}|$ (here, of course, $s\geq 4$ (resp. 6, 8, 12, 24, 20) if $m=1$ and $G_1/{C_1}\simeq \bar{A}_2(4)$ (resp. $^{2}\bar{A}_3(9), \bar{D}_4(2)$, Suz, $\cdot 1, \cdot 2$) and $s=48$ or 50 if $m=2$ and $G_1/{C_1}$ or $G_2/{C_2}\simeq \cdot 1$). This proves (10.2)(b).$\Box$

\subsection{}\textit{Proof of (10.3).} If $I_{\textnormal{Alt}}=\varnothing, E=\{1\}$ then (10.3)(a) follows from (9.1)(b), (c) (and the difinition of $F_H(n)$). If $E\neq \{1\}$ it follows from (10.2).

We have $F(G, n)\leq 2\cdot n^{2\,\log{n}+4.32}$ for all sporadic $G$ (since $\min{\tilde{n}}\leq$ (adjusted estimate) in Table T7.2). We have $5-4.32=0.68$ and $12^{0.68}=5.42 > 2$. Thus $|\Aut{G}|\leq n^{2\,\log{n}+5}$ for a sporadic $G$. We also have $f(n)\leq n^{2\,\log{n}+5}$ if $n\geq 4$ by (A7)(b), $4.1\cdot f(4) < 4^{2\,\log{4}+5}, 12.61\cdot f(6) < 6^{2\,\log{6}+5}, 27.69\cdot f(8) < 8^{2\,\log{8}+5}$ whence (10.3)(b) holds (by (10.1)(a)) for all $4\leq n\leq 12$. If $n=2$ or 3 then by Table T2.7 (groups of small degree) $|\Aut{(H\times E)}|\leq 120$ (resp. 5040) which $\leq 2^7$ (resp. $3^{8.17}$), whence (b) holds also for $n=2$ and 3. If $E=\{1\}$ it leaves (in view of (9.1)(b), (c)) only the case when $m=2, n_1\,n_2=48$ or 50 and $|\Aut{H}|/{n_1\,n_2\,f(n_1\,n_2)}\leq 3.05$. One verifies that then $|\Aut{H}|\leq (n_1\,n_2)^{2\,\log{(n_1\,n_2)}+5}$.  

Thus $E\neq \{1\}$. In this case our claim follows from (10.2) in view of (A7)(c) and the evident inequality $2\,n^{2\,\log{n}+1} < n^{2\,\log{n}+5}$ for $n\geq 2.$

To prove (10.3)(c) note that by (A7)(c) $n\,f(n)\leq 2\,n^{2\,\log{n}+1}$ for $n\geq 37$ and, in addition, one verifies directly that $|\Aut{G}|\leq 2\cdot39^{2\,\log{39}+1}$ for any sporadic $G$. Further one verifies that $3.05\cdot 48\cdot f(48) < 2\cdot 48^{2\,\log{48}+1}$ and $1.43\cdot 50\cdot f(50) < 2\cdot 50^{2\,\log{50}+1}$. These remarks together with (9.1)(b), (c) and (10.2) imply (10.3). $\Box$ 

\newpage
\section{\textbf{Estimates for finite quasi-primitive group.}}

Let $k$ be an algebraically closed field of characteristic exponent $p=p(k)$, $M$ a finite subgroup of $GL_n(k)$, and $C$ the center of $M$. Recall (see R. Brauer [\cite{5}, p. 64, where $q$ should be $K$]) that $M$ is quasi-primitive if it is irreducible and if for every normal subgroup $N$ of $M$, any two irreducible constituents are equivalent, Of course, any irreducible representation of a centrally simple group is quasi-primitive. Let $\bar{S}$ be the socle of $M/C$ and $S$ its preimage in $M$. Recall that $f(n)=(2\,n+1)^{2\,\log_3{(2\,n+1)}+1}$. 

\th{}
Suppose that $M$ is quasi-primitive. Then $M$ contains\\ 

(i) a normal subgroup $A$ isomorphic to a direct product of alternating groups $\textnormal{Alt}_{m_i}, m_i\geq 10$;\\

(ii) a normal perfect subgroup $L$ centrally isomorphic to a direct product of finite simple groups of Lie $p$-type;\\

(iii) a normal subgroup $E$ isomorphic to a direct product of extra-special groups whose orders are powers of distinct primes $q$ with $q|n, q\neq p$;\\

such that\\

(a) $
|M/{CL}|\leq
\begin{cases}
n\,f(n) & \text{if $n=2, 3, 5, 7, 9, 10, 11$} \\
4.1\,f(4) & \text{if $n=4$} \\
12.61\,f(6) & \text{if $n=6$} \\
27.69\,f(8) & \text{if $n=8$} \\
231\,f(12) & \text{if $n=12$} \\
(n+2)! & \text{if $n> 12$} \\
\end{cases} 
$
\ \\

(b) \begin{tabular}{cccc}
$|M/{ACL}|$ & $\leq$ & $n^{2\,\log{n}+5}$ & if $n\geq 2$ \\
and & $\leq$ & $2\,n^{2\,\log{n}+1}$ & if $n\geq 39$\\
\end{tabular}
\ \\

(c) $|M/S|\leq [\log{n}]!\,n^2\cdot |N_M(E)/{Z_M(E)}\cdot E|.$
\eth

\lem{}
$S$ is a central product of $C$ with centrally simple groups $G_1,\ldots, G_m$ and extraspecial groups $E_1,\ldots, E_r$ with $|E_i|=p_i^{1+2\,a_i}, a_i > 0, p_i\neq p, p_i|n, p_i\neq p_j$ if $i\neq j, i, j=1, \ldots, r.$
\elem

\prf{}
$\bar{S}$ is a direct product of simple groups. Write $\bar{S}=\bar{G}_1\times \cdots\times \bar{G}_m\times \bar{E}_1\times \cdots \times \bar{E}_r$ where the $\bar{G}_i$ are simple non-commutative and $\bar{E}_i$ are elementary abelian, $|\bar{E}_i|=p_i^{b_i}, p_i\neq p_j$ if $i\neq j.$ Let $\tilde{G}_i, i=1, \ldots, m,$ be the preimage of $\bar{G}_1$ in $M$. Then $G_i:=[\tilde{G}_i, \tilde{G}_i]$ is centrally simple. 

Let $\tilde{E}_i, i=1, \ldots, r,$ be the preimage of $\bar{E}_i$ in $M$. Note that since $M$ is primitive every normal commutative subgroup of $M$ is central. The pairing $\bar{E}_i\times\bar{E}_i\rightarrow C$ given by $[\bar{x}, \bar{y}]=[x, y]$ where $x, y$ are preimages of $x$ and $y$ in $\tilde{E}_i$ does not depend on the choice of $x$ and $y$. Let $\bar{F}_i:= \left\{x\in\bar{E}_i| [x, E_i]=\{1\}\right\}$ and $F_i$ the preimage of $\bar{F}_i$ in $E_i$. Then $F_i$ is commutative. As it is a characteristic subgroup of $\tilde{E}_i$ and, therefore, of $M$, it is normal. Therefore $F_i\subseteq C,$ i.e., $\bar{F}_i=\{1\}.$ Thus our pairing is non-degenerate. It follows now from D. Gorenstein [\cite{13}, ?] that $\tilde{E}_i=E_i\cdot C, i=1, \ldots, r,$ where the $E_i$ are extraspecial, $|E_i|=p_i^{1+2\,a_i}, p_i\neq p_j$ if $i\neq j$. The center of $E_i$ is of order $p_i$ and is contained in scalar matrices of degree $n$. Therefore $p_i|n, i=1, \ldots, r.$ Finally, $p_i\neq p, i=1, \ldots, r,$ since the $p$-subgroup of $C$ is trivial ($k$ contains only trivial $p$-th roots of 1).
\epr
\subsection{}\textit{Proof of (11.1).} Let $V$ be an irreducible component of the action of $S$ on $k^n$ and $V$ and $v:=\dim{V}$. Let, further, $D$ be the product of the $G_i$ which are isomorphic neither to alternating groups $\textnormal{Alt}_a, a\geq 10$, nor to groups of Lie $p$-type. Let $U$ (resp. $W$) be an irreducible component of the action of $ADE$ (resp, $DE$) on $V$ (resp. U). Let $w:=\dim{W}, u:=\dim{U}.$ We have by (10.1)(a) that $\Aut{ADE}$ is bounded as claimed by (11.1)(a) with $n$ replaced by $u$. Let $F(n)$ denote the right-hand side of (11.1)(a). Thus $|\Aut{ADE}|\leq F(u)$.

Write $V=U\otimes\bar{U}$ where $\bar{U}$ is an irreducible representation of $L$. Let $N:=N_{GL_n(k)}(L)/L\cdot Z_{GL_n(k)}(L)$. Then $N$ can be identified with a subgroup of $\textnormal{Out}\,L$. Since $k^n$ is a multiple of $\bar{U}$ as the $L$-module (as $k^n$ is primitive) it follows that $N$, acting on the irreducible representations of $L$ by $\varphi^n(l):=\varphi(\tilde{n}\,l\,\tilde{n}^{-1})$, where $\tilde{n}$ is a lift of $n$ to $N_{GL_n(k)}(L)$, preserves the equivalence class of $\varphi$. Then (5.5) together with the argument of the proof of (9.2) shows therefore that $|N|\leq [\log{\bar{u}}]!\,\bar{u}^2$.

The action of $M$ on $S$ induces automorphisms of $A, D, E, L.$ Let us denote by $\omega_A$ and $\bar{\omega}_A$ (resp. $\omega_D$, etc) the corresponding maps $\omega_A: M\rightarrow \Aut{A}$ and $\bar{\omega}_A:= M\rightarrow \textnormal{Out}\,A$ etc. We have $\bar{\omega}_L(M)\subseteq N$ whence $|\bar{\omega}_L(M)|\leq [\log{\bar{u}}]!\,\bar{u}^2$. Clearly $\ker{\omega_L}\bigcap \omega_{ADE}\subseteq Z_M(S)=C$. Therefore $|M/{CL}|\leq |\bar{\omega}_L(M)|\cdot |\omega_{ADE}(M)|\leq |\bar{\omega}_L(M)|\cdot |\Aut{ADE}|\leq [\log{\bar{u}}]!\,\bar{u^2}\cdot F(u)$. One checks easily that $\bar{u}\,u\leq n$ implies $|M/{CL}|\leq F(n),$ thus proving (11.1)(a).

For (b) we have by (10.3)(b) that $|\Aut{DE}|\leq w^{2\,\log{w}+5}$. We have also $|\bar{\omega}_A(M)|\leq [\log{u/w}]!\,2^{\log{u/w}}=[\log{u/w}]!\,u/w$ and $\ker{\omega_A}\bigcap \ker{\omega_L}\bigcap \ker{w_{DE}}=\left\{1\right\}$ whence as above 
$$
|M/{AC}|\leq |\bar{\omega}_A(M)|\cdot |\omega_L(M)|\cdot |\Aut{DE}|\leq [\log{\bar{u}}]!\,\bar{u}^2\,[\log{u/w}]!\,u/w\,w^{2\,\log{w}+5}.
$$

From $\bar{u}\cdot (u/w)\cdot w\leq n$ it follows that this $\leq n^{2\,\log{n}+5}$ as desired.

For the second estimate in (b) we have as above
$$
|M/{ACL}|\leq [\log{\bar{u}}]!\,\bar{u}^2\cdot [\log{u/w}]!\,u/w\cdot |\Aut{DE}|.
$$

If $w\geq 39$ we have by (10.3)(c) that $|M/{ACL}|\leq 2\,w^{2\,\log{w}+1}$ which together with the above and the inequality $\bar{u}\cdot (u/w)\cdot w\leq n$ implies the desired inequality. Assume that $s:=\bar{u}\cdot (u/w)=2$ (resp. 3) and $38 \geq w \geq [39/s]$. Then $|M/{ACL}|\leq s^2\,|\Aut{DE}|$ and one verifies using (10.2)(a) and (9.1)(c) that our claim holds in this case as well. If $s\geq 4$ then we can use $[\log{\bar{u}}]!\,\bar{u}^2\cdot [\log{u/w}]!\,u/w \leq f(s)$ and invoke (A6) to conclude the proof.

Finally (c) is an evident corollary of the estimate in (9.1)(c) on $\textnormal{Out}\,H,$ the estimate we established on $(\textnormal{Out}\,L)_{\bar{U}}$ and of an evident estimate on the contribution of $N_M(E)$ to $\textnormal{Out}\,E. \Box$ 

\newpage
\section{\textbf{Estimates for irreducible finite groups.}} 

Let now $k$ be an algebraically closed field of characteristic exponent $p=p(k)$.

\th{}
Let $H$ be an irreducible finite subgroup of $GL_n(k)$. Then $H$ contains \\

(i) a commutative normal diagonalizable subgroup $B$;\\

(ii) a normal perfect subgroup $L$ centrally isomorphic to a direct product of simple groups of Lie $p$-type\\

so that\\

(a) $|H/{BL}|\leq 
\begin{cases}
n^4\,(n+2)! & \text{if $n\leq 63$} \\
(n+2)! & \text{if $n > 63$} \\
\end{cases}
$
\ \\

Moreover\\

(b) for all $n\geq 2$
$$
|H/{BL}|\leq (n+2)!\cdot n^{4020/{((n-20)^2+1000)}}
$$
\ \\

(c) for $n\leq 63$ we have
$$
|H/{BL}|\leq (n+2)!\cdot n^{4\,\alpha_{\textnormal{irr}}}
$$
\ \\

where $\alpha_{\textnormal{irr}}$ is given by Table T12.1.
\eth

\newpage
\pagestyle{empty} 

\begin{center}
Table T12.1.\\
\end{center}
\ \\

\tiny
\begin{center}
\begin{tabular}{ccccccccccccccc}

& \ & \ & \ & \ & \ & \ & \ & \ & \ & \ & \ & \ & \ & \ \\
& $\alpha_{\textnormal{irr}}$ & $\alpha_{\textnormal{all}}$ & \ & \ & \ & \ & \ & \ & \ & \ & \ & \ & $\alpha_{\textnormal{irr}}$ & $\alpha_{\textnormal{all}}$ \\ 

& \ & \ & \ & \ & \ & \ & \ & \ & \ & \ & \ & \ & \ & \ \\
$2$ & $.34$ & $.34$ & \ & \ & \ & \ & \ & \ & \ & \ & \ & $33$ & $.83$ & $.83$ \\

& \ & \ & \ & \ & \ & \ & \ & \ & \ & \ & \ & \ & \ & \ \\
$3$ & $.7$ & $.7$ & \ & \ & \ & \ & \ & \ & \ & \ & \ & $34$ & $.53$ & $.61$ \\

& \ & \ & \ & \ & \ & \ & \ & \ & \ & \ & \ & \ & \ & \ \\
$4$ & $.78$ & $.78$ & \ & \ & \ & \ & \ & \ & \ & \ & \ & $35$ & $0$ & $.6$ \\

& \ & \ & \ & \ & \ & \ & \ & \ & \ & \ & \ & \ & \ & \ \\
$5$ & $.37$ & $.53$ & \ & \ & \ & \ & \ & \ & \ & \ & \ & $36$ & $.77$ & $.77$ \\

& \ & \ & \ & \ & \ & \ & \ & \ & \ & \ & \ & \ & \ & \ \\
$6$ & $.81$ & $.81$ & \ & \ & \ & \ & \ & \ & \ & \ & \ & $37$ & $0$ & $.55$ \\

& \ & \ & \ & \ & \ & \ & \ & \ & \ & \ & \ & \ & \ & \ \\
$7$ & $.67$ & $.76$ & \ & \ & \ & \ & \ & \ & \ & \ & \ & $38$ & $.47$ & $.54$ \\

& \ & \ & \ & \ & \ & \ & \ & \ & \ & \ & \ & \ & \ & \ \\
$8$ & $.88$ & $.88$ & \ & \ & \ & \ & \ & \ & \ & \ & \ & $39$ & $.71$ & $.71$ \\

& \ & \ & \ & \ & \ & \ & \ & \ & \ & \ & \ & \ & \ & \ \\
$9$ & $.89$ & $.89$ & \ & \ & \ & \ & \ & \ & \ & \ & \ & $40$ & $.44$ & $.48$ \\

& \ & \ & \ & \ & \ & \ & \ & \ & \ & \ & \ & \ & \ & \ \\
$10$ & $.58$ & $.79$ & \ & \ & \ & \ & \ & \ & \ & \ & \ & $41$ & $0$ & $.47$ \\

& \ & \ & \ & \ & \ & \ & \ & \ & \ & \ & \ & \ & \ & \ \\
$11$ & $.1$ & $.81$ & \ & \ & \ & \ & \ & \ & \ & \ & \ & $42$ & $.64$ & $.64$ \\

& \ & \ & \ & \ & \ & \ & \ & \ & \ & \ & \ & \ & \ & \ \\
$12$ & $.94$ & $.94$ & \ & \ & \ & \ & \ & \ & \ & \ & \ & $43$ & $0$ & $.41$ \\

& \ & \ & \ & \ & \ & \ & \ & \ & \ & \ & \ & \ & \ & \ \\
$13$ & $0$ & $.81$ & \ & \ & \ & \ & \ & \ & \ & \ & \ & $44$ & $.37$ & $.4$ \\

& \ & \ & \ & \ & \ & \ & \ & \ & \ & \ & \ & \ & \ & \ \\
$14$ & $.62$ & $.77$ & \ & \ & \ & \ & \ & \ & \ & \ & \ & $45$ & $.57$ & $.57$ \\

& \ & \ & \ & \ & \ & \ & \ & \ & \ & \ & \ & \ & \ & \ \\
$15$ & $.97$ & $.97$ & \ & \ & \ & \ & \ & \ & \ & \ & \ & $46$ & $.34$ & $.34$ \\

& \ & \ & \ & \ & \ & \ & \ & \ & \ & \ & \ & \ & \ & \ \\
$16$ & $.93$ & $.93$ & \ & \ & \ & \ & \ & \ & \ & \ & \ & $47$ & $0$ & $.32$ \\

& \ & \ & \ & \ & \ & \ & \ & \ & \ & \ & \ & \ & \ & \ \\
$17$ & $0$ & $.77$ & \ & \ & \ & \ & \ & \ & \ & \ & \ & $48$ & $.49$ & $.49$ \\

& \ & \ & \ & \ & \ & \ & \ & \ & \ & \ & \ & \ & \ & \ \\
$18$ & $.98$ & $.98$ & \ & \ & \ & \ & \ & \ & \ & \ & \ & $49$ & $0$ & $.25$ \\

& \ & \ & \ & \ & \ & \ & \ & \ & \ & \ & \ & \ & \ & \ \\
$19$ & $0$ & $.81$ & \ & \ & \ & \ & \ & \ & \ & \ & \ & $50$ & $.26$ & $.26$ \\

& \ & \ & \ & \ & \ & \ & \ & \ & \ & \ & \ & \ & \ & \ \\
$20$ & $.89$ & $.89$ & \ & \ & \ & \ & \ & \ & \ & \ & \ & $51$ & $.41$ & $.41$ \\

& \ & \ & \ & \ & \ & \ & \ & \ & \ & \ & \ & \ & \ & \ \\
$21$ & $.97$ & $.97$ & \ & \ & \ & \ & \ & \ & \ & \ & \ & $52$ & $.22$ & $.22$ \\

& \ & \ & \ & \ & \ & \ & \ & \ & \ & \ & \ & \ & \ & \ \\
$22$ & $.63$ & $.78$ & \ & \ & \ & \ & \ & \ & \ & \ & \ & $53$ & $0$ & $.16$ \\

& \ & \ & \ & \ & \ & \ & \ & \ & \ & \ & \ & \ & \ & \ \\
$23$ & $0$ & $.75$ & \ & \ & \ & \ & \ & \ & \ & \ & \ & $54$ & $.32$ & $.32$ \\

& \ & \ & \ & \ & \ & \ & \ & \ & \ & \ & \ & \ & \ & \ \\
$24$ & $.95$ & $.95$ & \ & \ & \ & \ & \ & \ & \ & \ & \ & $55$ & $0$ & $.07$ \\

& \ & \ & \ & \ & \ & \ & \ & \ & \ & \ & \ & \ & \ & \ \\
$25$ & $0$ & $.75$ & \ & \ & \ & \ & \ & \ & \ & \ & \ & $56$ & $.13$ & $.13$ \\

& \ & \ & \ & \ & \ & \ & \ & \ & \ & \ & \ & \ & \ & \ \\
$26$ & $.61$ & $.73$ & \ & \ & \ & \ & \ & \ & \ & \ & \ & $57$ & $.23$ & $.23$ \\

& \ & \ & \ & \ & \ & \ & \ & \ & \ & \ & \ & \ & \ & \ \\
$27$ & $.92$ & $.92$ & \ & \ & \ & \ & \ & \ & \ & \ & \ & $58$ & $.09$ & $.09$ \\

& \ & \ & \ & \ & \ & \ & \ & \ & \ & \ & \ & \ & \ & \ \\
$28$ & $.74$ & $.74$ & \ & \ & \ & \ & \ & \ & \ & \ & \ & $59$ & $0$ & $0$ \\

& \ & \ & \ & \ & \ & \ & \ & \ & \ & \ & \ & \ & \ & \ \\
$29$ & $0$ & $.69$ & \ & \ & \ & \ & \ & \ & \ & \ & \ & $60$ & $.13$ & $.13$ \\

& \ & \ & \ & \ & \ & \ & \ & \ & \ & \ & \ & \ & \ & \ \\
$30$ & $.88$ & $.88$ & \ & \ & \ & \ & \ & \ & \ & \ & \ & $61$ & $0$ & $0$ \\

& \ & \ & \ & \ & \ & \ & \ & \ & \ & \ & \ & \ & \ & \ \\
$31$ & $0$ & $.67$ & \ & \ & \ & \ & \ & \ & \ & \ & \ & $62$ & $0$ & $0$ \\

& \ & \ & \ & \ & \ & \ & \ & \ & \ & \ & \ & \ & \ & \ \\
$32$ & $.64$ & $.65$ & \ & \ & \ & \ & \ & \ & \ & \ & \ & $63$ & $.04$ & $.04$ \\

\end{tabular}

\end{center}

\newpage
\normalsize
\pagestyle{plain}

\subsection{}\textit{Proof.} Set $V:=k^n$ and let $V=\otimes^{m}_{i=1}{V_i}$ be an imprimitivity system for $H$ on $V$. Let $H_i:=\left\{h\in H| h\,V_i= V_i\right\}$. Then $H_i$ is primitive on $V_i$. Let $L_i$ be the largest perfect normal subgroup of $H_i$ centrally isomorphic to a direct product of finite simple groups of Lie $p$-type and $C_i$ the center of $H_i$. Then

\subsubsection{}$|H_i/{L_i\,C_i}|$ satisfies (11.1)(a).\\

Let $\varphi: H\rightarrow \Sym_m$ be the homomorphism defined by $h(V_i)=V_{\varphi(h)i}$ and let $M:=\ker{\varphi}$. Then $MV_i=V_i$ for $i=1,\ldots,m.$ In particular, we have homomorphisms $\omega_i: M\rightarrow H_i.$ Set $M_i:=\omega_i(M), i=1,\ldots,m.$ Since $M$ is normal in $H$ we have that $M_i$ is normal in $H_i$. Therefore each perfect factor of $L_i$ either it is contained in $M_i$ or $M_i$ intersects it in its center. Let $L'_i$ be the largest perfect normal subgroup of $M_i$ which is centrally isomorphic to a direct product of finite simple groups of Lie $p$-type and let $C'_i$ be the center of $M_i$. Then since $M_i\,L_i\,C_i/{L_i\,C_i}\simeq M_i/(L_i\,C_i\bigcap M_i)$ and since by the above comments $C_i\,L_i\bigcap M_i=C'_i\,L'_i$ we get

\subsubsection{}$|M_i/{L_i\,C_i}|\leq |H_i/{L_i\,C_i}|$.\\

Let $L:=\bigcap^{m}_{i=1}{\omega_i^{-1}(L'_i)}, B:=\bigcap^{m}_{i=1}{\omega_i^{-1}(C_i)}$; clearly $L$ and $B$ are normal. Consider the evident homomorphisms $\omega:=\otimes^{m}_{i=1}{\omega_i}: M \rightarrow \prod^{m}_{i=1}{M_i}.$ Note that

\subsubsection{}$\ker{\omega}=\left\{1\right\}.$ \\

We have $\omega(B)\subseteq \prod^{m}_{i=1}{C'_i}.$ This and (12.2.3) give that

\subsubsection{}$B$ is a commutative normal subgroup of $M$.\\

We also have $\omega(L)\subseteq \prod^{m}_{i=1}{L'_i}.$ Moreover, since $\ker{\omega}=\left\{1\right\}$ and the projection of $\omega(L)$ one each $L'_i$ is the whole $L'_i$ it follows that

\subsubsection{}$L$ is perfect and centrally isomorphic to a direct product of finite simple groups of Lie $p$-type. Clearly, $\omega(L)=\omega(M)\bigcap\prod^{m}_{i=1}{L'_i}.$ Thus 
$$
\omega(M)\,\prod^{m}_{i=1}{L'_i\,C'_i}/{\prod^{m}_{i=1}{L'_i\,C'_i}}\simeq \omega(M)/{\omega(M)}\bigcap\prod^{m}_{i=1}{L'_i\,C()}=\omega(M)/{\omega(LB)}\simeq M/{LB}.
$$

Since evidently 
$$
\left|\omega(M)\,\prod^{m}_{i=1}{L'_i\,C'_i}/{\prod^{m}_{i=1}{L'_i\,C'_i}}\right|\leq |\prod^{m}_{i=1}{M_i}/{\prod^{m}_{i=1}{L'_i\,C'_i}}|
$$

we get from (12.2.2):
$$
\left|M/{LB}\right|\leq \prod^{m}_{i=1}{|H_i/{L_i\,C_i}|}.
$$

Since $H$ is irreducible on $V$ it follows from Clifford theory that $H_i\simeq H_j$ for $i,j=1,\ldots,m.$ Thus $\prod^{m}_{i=1}{|H_i/{L_i\,C_i}|}=|H_1/{L_1\,C_1}|^m$. Since $M=\ker{\left\{H\rightarrow \Sym_m\right\}}$ this implies

\subsubsection{}$|H/{LC}|\leq m!\,|H_1/{L_1\,C_1}|^m$.\\

If now $n/m(=\dim{V_1}) > 12$ then (A9)(i) implies that $|H/{BL}|\leq (n+2)!$ Thus it is sufficient to consider the cases when $n/m\leq 12.$ Set $r:=n/m.$ For $2\leq r\leq 12$ values we use the estimate $|H_1/{L_1\,C_1}| < t_r$ where $t_r$ is from Table T12.2. For $r=2, 3, 4, 5$ we get a good estimate because all possible simple groups $H_1$ are known from Table T2.7 (groups of small degree) and we get good estimates on their normalizers from W. Feit [\cite{?}, p. 76] for $r=2$ and 3 (maxima are respectively for $H_1\simeq \textnormal{Alt}_5$ and $\textnormal{Alt}_7$), from Zalessky [\cite{33}, p. 95] for $r=4$ and 5 (maxima in both cases are for $\Aut{B_2(3)}$). For $6\leq r\leq 11, r\neq 8,$ we take $t_r=r\,f(r)$ (use (6.1)(d) to justify). For $r=8$ (resp. 12) we take $t_8=|\Aut{D_4(2)}|$ (resp. $t_{12}=|\Aut{(\textnormal{Suz})}|$). Of course, it has to be verified that admitting central products for $H_1$ lowers the estimate, but this is straightforward in our range.

Set $F(m, r, t):=t^m\,(m!)/{(r\,m+2)!}$ For each $r\leq 12$ we find, using a computer, the first $m$ such that $F(m, r, t_r) < 1, F(m+1, r, t_r)\geq 1.$ We denote this $m$ by $m_r;$ it is given in Table T12.2.

It is then easy to see that $F(m, r, t_r)< 1$ for all $m > m_r$. Namely
$$
F(m+1, r, t_r)=F(m, r, t_r)\cdot t_r\cdot (m+1)/{(r\,m+3)\,(r\,m+4)\ldots(r\,m+r+2)}.
$$

It is evident that $F_1(m, r, t_r):=t_r(m+1)/(r\,m+3)\ldots(r\,m+r+3)$ decreases when $m$ increases. By the definition of $m_r$ we have $F_1(m_r, r, t_r) < 1.$ Thus the same holds for all $m\geq m_r.$

Note that maximum of $r\,m_r, 2\leq r\leq 12,$ is 63. Thus $|M| < (n+2)!$ if $n > 63.$ For each $n\leq 63$ the computer takes for an estimate on $|M|$ the maximum of $(n+2)!$ and all $t^{m}_r(m!)$ such that $r\leq 12, m\,r=n.$ One then checks, on the computer again, that these maxima satisfy the inequalities claimed in (12.1). 
\ \\

\begin{center}
Table T12.2.\\
\end{center}
\ \\
\scriptsize
\begin{center}
\begin{tabular}{|c|c|c|c|c|c|c|c|c|c|c|c|}
\hline

& \ & \ & \ & \ & \ & \ & \ & \ & \ & \ & \ \\
$r$ & $2$ & $3$ & $4$ & $5$ & $6$ & $7$ & $8$ & $9$ & $10$ & $11$ & $12$ \\
\hline

& \ & \ & \ & \ & \ & \ & \ & \ & \ & \ & \ \\
$t_r$ & $60$ & $2520$ & $51840$ & $51840$ & $1.24\cdot 10^7$ & $6.6\cdot 10^7$ & $1.05\cdot 10^9$ & $1.22\cdot 10^9$ & $4.47\cdot 10^9$ & $1.5\cdot 10^{10}$ & $9\cdot 10^{11}$ \\
\hline

& \ & \ & \ & \ & \ & \ & \ & \ & \ & \ & \ \\
$m_r$ & $30$ & $21$ & $13$ & $4$ & $7$ & $5$ & $3$ & $3$ & $2$ & $2$ & $1$ \\
\hline

\end{tabular}

\end{center}

\newpage
\normalsize
\section{\textbf{Estimates for arbitrary finite linear groups.}}

Let $k$ be an algebraically closed field of characteristic exponent $p=p(k)$.

\th{}
Let $G$ be a finite subgroup of $GL_n(k)$. Then $G$ contains\\

(i) a triangulizable normal subgroup $T, T\supseteq O_p(G),$\\

(ii) a normal subgroup $L$ such that $L\supseteq O_p(G)$ and $L/{O_p(G)}$ is perfect and centrally isomorphic to a direct product of finite simple groups of Lie $p$-type\\

so that\\

(a) $|G/{LT}|\leq
\begin{cases}
n^4\,(n+2)! & \text{if $n\leq 63$} \\
(n+2)! & \text{if $n > 63$}\\
\end{cases}
$ 
\ \\

(b) for all $n\geq 2$
$$
|G/{LT}|\leq (n+2)!\,n^{4020/{((n-20)^2+1000)}}
$$
\ \\

(c) for $n\leq 63$ we have
$$
|G/{LT}|\leq (n+2)!\,n^{4\,\alpha_{\textnormal{all}}}
$$
\ \\

where $\alpha_{\textnormal{all}}$ is given in Table T12.1.

\eth

\subsection{}This result implies of R. Brauer and W. Feit \cite{6}.

\textbf{Corollary.} Suppose that $p^a$ is the highest power of $p$ dividing $|G|$. Then $G$ contains a normal commutative diagonalizable subgroup $B$ such that
$$
\left|G/B\right|\leq p^{3\,a}\,|G/L|\leq 
\begin{cases}
p^{3\,a}\,n^4\,(n+2)! & \text{if $n\leq 63$} \\
p^{3\,a}\,(n+2)! & \text{if $n > 63$}\\
\end{cases}
$$

\prf
Let $p^c$ be the order of the Sylow $p$-subgroup of $L/O_p(G).$ Then one easily sees from (4.4.1) that $|L/{O_p(G)}|\leq p^{3\,c}$. Let $p^t=|O_p(G)|$. Write $D$ for a $p'$-complement to $O_p(H)$ in $T$. Then $D$ is commutative. The action of $D$ by conjugation on $R:=O_p(H)$ defines a linear action of $D$ in a $\mathbb{F}_p$-vector space $\bar{R}:=R/[R, R]\cdot R^p,$ i.e., a homomorphism $\omega:=D\rightarrow GL(\bar{R})$. Let $\bar{t}=\dim_{\mathbb{F}_p}{R}$ so that $GL(\bar{R})\simeq GL_{\bar{t}}(\mathbb{F}_p)$. It is evident that every commutative $p'$-subgroup of $GL_{\bar{t}}(\mathbb{F}_p)$ has order $\leq p^{\bar{t}}$. Let $B:=\ker{\omega}$. By the above $|B|\geq |D|/p^{\bar{t}}\geq |D|/{p^t}$. By D. Gorenstien \cite{13} $B$ (acting trivially on $\bar{R}$) acts trivially on $R$. Thus $B$ is the $p'$-component of the center of $T$. In particular, it is a characteristic subgroup of $T$ and, therefore, a normal subgroup of $G$. We have
$$
|G/B|=|G/L|\cdot |L/R|\cdot |D/B|\cdot |R|\leq |G/L|\cdot p^{3\,c}\cdot p^t\cdot p^t.
$$

Since $p^t\cdot p^c$ is the order of Sylow $p$-subgroup of $L$ we have $t+c\leq a$. Thus the above gives
$$
|G/B|\leq p^{3\,a}\,|G/L|
$$
as claimed.
\epr

\subsection{}\textit{Proof of (13.1).} Set $V:=k^n$. Let $V_1:=V\supseteq V_2\supseteq \ldots \supseteq V_m\supseteq V_{m+1}=:0$ be a sequence of $G$-submodules of $V$ such that $V_i \neq V_{i+1}, i=1,\ldots, m$ and $W_i:=V_i/{V_{i+1}}$ is irreducible for $G$ for $i=1,\ldots, m.$ Set $n_i:=\dim{W_i}$. The action of $G$ on $W_i$ defines a homomorphism $\omega_i: G \rightarrow GL(W_i)\simeq GL_{n_i}(k)$. Set $H_i:=\omega_i(G).$

Set $\omega:= \oplus{\omega_i}: G\rightarrow \prod{GL_{n_i}(k)}$. The kernel of $\omega$ is a unipotent subgroup of $G$. Since each $H_i$ is irreducible it has no unipotent normal subgroups and, therefore, $\ker{\omega}$ is the largest unipotent normal subgroup of $G$. Since $G$ is finite this latter is just $O_p(G).$ Thus 

\subsubsection{} $\ker{\omega}=O_p(G).$\\

Set $H:= G/{O_p(G)}$ and let $\bar{\omega}_i: H\rightarrow H_i$ be the map induced by $\omega_i, i=1,\ldots, m.$ Set $\bar{\omega}:= \oplus{\bar{\omega}_i}.$ We have 

\subsubsection{} $\ker{\bar{\omega}}=\{1\}.$\\

Let $L_i$ and $B_i$ be the subgroups claimed in (12.1) for $H_i$. Let $L':= \bigcap^{m}_{i=1}{\bar{\omega}^{-1}_i{(L_i)}}, B':= \bigcap^{m}_{i=1}{\bar{\omega}^{-1}_i{(B_i)}}$. As in the proof of (12.2.4) and (12.2.5) we get 

\subsubsection{}$B'$ is a commutative normal subgroup of $H$.\\

\subsubsection{}$L'$ is a perfect and centrally isomorphic to a direct product of finite simple groups of Lie $p$-type. 

We also have
$$
\bar{\omega}(H)\,\prod^m_{i=1}{L_i\,B_i}/\prod^m_{i=1}{L_i\,B_i}\simeq \bar{\omega}(H)/{\bar{\omega}(H)}\bigcap \prod^m_{i=1}{L_i\,B_i}\simeq \bar{\omega}(H)/{\omega(L'\,B')}\simeq H/{L'\,B'}.
$$

Since 
$$
|\bar{\omega}(H)\,\prod^{m}_{i=1}{L_i\,B_i}/{\prod^m_{i=1}{L_i\,B_i}}|\leq |\prod^m_{i=1}{H_i}/\prod^m_{i=1}{L_i\,B_i}|
$$

we have that 

\subsubsection{}$|H/{L'\,B'}|\leq \prod^m_{i=1}{|H_i/{L_i\,B_i}|}$. \\

If now all $n_i\geq 64$ then our claim follows from (A9)(ii). In the remaining cases we use estimates on $|H_i/{L_i\,B_i}|$ for $n_i\leq 63$ obtained by the computer as described in (12.2). For each pair $m_1, m_2, 2\leq m_1, m_2\leq 64$ we take for a new estimate for $m_1+m_2$ the maximum of the estimate obtained before for $m_1+m_2$ and of the product of these estimates for $m_1$ and $m_2$. We repeat this procedure until it stabilizes. It turns out that it gives new (compared with $\S 12$) values only for $n\leq 55$. Then our checks using the computer that the estimates claimed in (13.1) for $n\leq 126$ hold for $|H/{B'\,L'}|$. The case when some $n_i\geq 63$ and some $\geq 64$ is handled as follows. If $I$ is a subset of $1, \ldots, m$ and $\sum_{i\in I}{n_i}\leq 126$ then, as remarked, the computer establishes the required estimate.

In view of (A9)(ii) it remains to show that if $A$ is the estimate (hold by the computer) for $r\leq 64$ and if $d\geq 64$ is such that $r+d > 126$ then 
$$
A\,(d+2)!\leq (r+d+2)!
$$

We know (by a check on a computer) that this hols for $d=64$. By induction suppose it holds for some $d$. Then for $d:=d+1$ we have
$$
A\,(d+3)!=A\,(d+2)!\,(d+3)\leq (r+d+2)!\,(d+3) < (r+d+3)!
$$

whence our present claim:

\subsubsection{}The estimates of (13.1) hold for $H$. \\

Let now $L$ be the preimage of $L$ in $G$ and $T$ the preimage of $B$. Since $B$ is diagonalizable $T$ is triagulizable. This concluds our proof of (13.1).

\newpage
\section{\textbf{Extension to infinite linear groups.}}

Let $k$ be an algebraically closed field of characteristic exponent $p=p(k)$. For a subgroup $H$ of $GL_n(k)$ let $H^c$ denote its Zariski closure and set $H^{\circ}:= H\bigcap (H^c)^{\circ}.$ When $H^c$ is semi-simple it contains, by J. Tits [\cite{27}, Theorems 3 and 4], a smallest (automatically connected) normal subgroup such that $H/H \bigcap \mathcal{F}$ is periodic; we call this $\mathcal{F}$ the Tits subgroup of $H^c$.

\th{}
Let $G$ be a subgroup of $GL_n(k)$. Then there exist\\

(i) a normal triangulizable subgroup $T$ of $G$,\\

(ii) normal subgroups $F, P, L$ of $G$ with $T=F\cap P\cap L$ \\

such that \\

(a) $P^c/T^c$ and $F^c/T^c$ are connected, semi-simple and commute,\\

(b) $F^c/T^c$ is the smallest among normal subgroups $\mathcal{H}$ of $G^c/T^c$ such that $\mathcal{H}\cap G/T$ projects onto the image of $(G/T)^{\circ}$ in the Tits subgroup of $(G/T)^c$, \\

(c) $FP\supseteq \mathcal{D}G^{\circ}$ and $G^{\circ}/FP$ is finite commutative; in particular, $F^c\,P^c=(G^c)^{\circ},$\\

(d) $P/T$ is a direct product of infinite simple groups of Lie $p$-type, \\

(e) $L/T$ is a direct product of finite simple groups of Lie $p$-type, \\

(f) $|G/{PFL}|\leq
\begin{cases}
n^4\,(n+2)! & \text{if $n\leq 63$} \\
(n+2)! & \text{if $n\geq 64$} \\
\end{cases}$ 
\ \\

Moreover\\

(g) if $G$ is finitely generated then $P=F$.
\eth

\subsection{}\textit{Proof of (14.1).} First, let us show how (g) follows from (a)-(f). By (c) $FP$ is of finite index in $G$. Therefore $FP$ is finitely generated if $G$ is finitely generated. But then $P/T=PF/F$ is also finitely generated. This is not so if $P/T$ is infinite (by (d)). Hence $P/T=\{1\}$ if $G$ is finitely generated, whence (g).

\subsection{}Now consider the case when Zarisski-closure of $G$ is connected and almost simple and $G$ is periodic. By J. Tits [\cite{27}, Theorem 3 and 4 (iv)] we have then $p>1$ and $G^c$ is defined over $\bar{\mathbb{F}}_p$. Let us fix a (rational) irreducible representation $g: G^c\rightarrow GL_d.$ Since $G^c$ is connected $g$ is also primitive (for otherwise $G^c$ would contain a subgroup of finite index preserving a decomposition of $k^d$ into a direct sum).

Since $G$ is irreducible and primitive on $V:=k^d$ there exists a finitely generated (and, hence, finite) subgroup $G_1$ of $G$ which is also irreducible and primitive. Write $G=\bigcup^{\infty}_{i=1}{G_i}$ where $G_{i+1}\supseteq G_i$ and $G_{i+1}\neq G_i$ (for example, $G_{i+1}=\left\langle x_i, G_i\right\rangle$ where $x_i\in G-G_i$). Let $S_i$ be the preimage in $G_i$ of the socle of $G_i/\textnormal{center}$. Then $S_i$ is a central product of centrally simple perfect groups $H_{i, 1}, \ldots, H_{i, m_i}$ and of extraspecial groups $E_{i, 1}, \ldots, E_{i, r_i}$ of relatively prime power orders. It is clear that each $H_{i, j}$ is contained in some $H_{i+1, s}.$ Similarly, $E_{i+1, j}\subseteq E_{i, s}$ if $E_{i+1, j}$ and $E_{i, s}$ have \textit{not} relatively prime orders. Therefore 

\subsubsection{}There exists $c$ such that $m_i=m_c, r_i=r_c$ for $i\geq c$ and $E_{i, j}=E_{c, s}$ for $i\geq c$ and appropriate $j$ and $s$.

We can assume (after renumeration) that $c=1, H_{i+1, j}\supseteq H_{i, j}, E_{i+1, j}=E_{i, j}$ for $i\geq 1.$ Then set $H_j:=\bigcup^{\infty}_{i=1}{H_{i, j}}, E_j:=E_{1, j}.$ We have that the $H_j$ and $E_s$ commute. Therefore so do $H_j^c$ and the $E_s$. Clearly $\prod_{j}{H_j^c}\cdot \prod_{s}{E_s}$ is a normal subgroup of $G^c$. Since $G^c$ is connected and almost simple this implies that $m_1=1$ and $r_1=0$ (that is, there is only one $H_j$ and no $E_s$).

Set $H_i:=H_{i, 1}, P:=\bigcup^{\infty}_{i=1}{H_i}$. Then

\subsubsection{}$H$ is centrally simple.\\

By a recent (overlapping) results of V. Belyaev \cite{3}, A. Borovik \cite{4}, N. Chernikov \cite{8}, B. Hartley and B. Shute \cite{15}, S. Thomas \cite{26} (see, for definiteness S. Thomas [\cite{26}, Theorem 2]) we get

\subsubsection{}$P$ is centrally simple of Lie $p$-type over a subfield $K$ of $\bar{\mathbb{F}}_p$. 

We have that $P$ is normal in $G$ and therefore, $G$ acts by automorphisms on $P$. By R. Steinberg [\cite{?}, Theorems 30 and 36], any automorphism of $P$ is a product of a diagonal, graph, field, and inner one. However, since $G^c$ is connected and since graph automorphisms do not belong to $G^c$ and field automorphisms do not induce automorphisms of $G^c$ we get

\subsubsection{}$G/P$ consists of diagonal automorphisms.\\

This implies

\subsubsection{}$G/P$ is finite commutative; it is given in column $A_d$ of Table T4.4.\\

\subsection{}Assume now that $\mathcal{Y}:=G^c$ is connected and semi-simple. Then by J. Tits [\cite{27}, Theorem 3 and 4(i)] $\mathcal{Y}$ contains a connected normal subgroup $\tilde{\mathcal{F}}$ which is the smallest such subgroup with the condition that the image of $G$ in $\mathcal{Y/\tilde{F}}$ is periodic. Write $\mathcal{Y}_1,\ldots, \mathcal{Y}_m$ for almost simple quotients of $\mathcal{Y/\tilde{F}}$. Let $G_i$ be the image of $G$ in $\mathcal{Y}_i, i=1,\ldots, m.$ Then (14.3) applies to $G_i$ and we have by (14.3.5) and (14.3.3) that $\mathcal{D}G_i$ is centrally simple of Lie $p$-type. Set $P_i:=\mathcal{D}G_i.$ Let $\varphi_i: G\rightarrow G_i$ be the natural projection. Then $\tilde{P}:= \bigcap^m_{i=1}{\varphi_i^{-1}\,(P_i)}$ satisfies

\subsubsection{}$\varphi_i(\tilde{P})=P_i, \mathcal{D}G\subseteq \tilde{P}, |G/{\tilde{P}}| < \infty$. \\

Let $\varphi: \mathcal{Y}\rightarrow \tilde{\mathcal{F}}/{\textnormal{center}}$ be the natural projection and let $P:=\ker{\varphi}\bigcap \tilde{P}$. Then $\varphi_i(p)$ is normal in $P_i$. Therefore, since $P_i$ is centrally simple, $\varphi_i(P)$ is either in the center of $P_i$ or contains $\mathcal{D}P_i$. Since $P$ is a subgroup of $\mathcal{Y}_1\ldots \mathcal{Y}_m$ it follows that

\subsubsection{}$P$ is centrally isomorphic to a direct product of some of the $P_i$.\\

Next, $P$ is normal in $\tilde{P}$ and, by construction, $\tilde{P}$ induces only inner automorphism of $P$. Thus

\subsubsection{}$\tilde{P}=P\cdot F$ where $F=Z_{\tilde{p}}(P)$. \\

Let $\mathcal{F}:=F^c, \mathcal{P}:=P^c.$ Then in view of (14.4.1) and (14.4.3).

\subsubsection{}$\mathcal{F\cdot P}=G^c$ \\

Now, if $\mathcal{H}$ is a smallest factor of $\mathcal{Y}$ such that $\mathcal{H}\bigcap G$ projects onto the image $G$ in $\tilde{\mathcal{F}}$ then, clearly, $\mathcal{H}\subseteq \mathcal{F}.$ If $\mathcal{H}\neq \mathcal{F}$ then $\varphi(G)/{\varphi(\mathcal{H}\bigcap G)}$ is isomorphic 
to the projection of $G$ onto the complement to $\mathcal{H}$ in $\mathcal{F}$ and, in particular, is $\neq \{1\}.$ Thus 

\subsubsection{}$\mathcal{F}$ is the smallest among normal subgroups $\mathcal{H}$ of $G^c$ such that $\varphi(\mathcal{H}\bigcap G)=\varphi(G)$.\\

\subsubsection{}Remark. Note that since $G$ may induce non-linear, hence, diagonal, automorphisms of $P$ there is no decomposition similar to (14.4.3) for $G$.\\

\subsubsection{}Example. Let $H:=PGL_n(\bar{\mathbb{F}}_p[t])$. For $c\in \mathbb{F}_p$ we have a specialization homomorphism $\varphi_c: t\rightarrow c$ of $\bar{\mathbb{F}}_p[t]$ to $\bar{\mathbb{F}}_p$. It induces an epimorphism $\varphi_c: H\rightarrow PGL_n(\bar{\mathbb{F}}_p)$. Let $c_1,\ldots, c_m\in \mathbb{F}_p$. Then we define a homomorphism $\varphi:=\id\times \prod^m_{i=1}{\varphi_{c_i}}: H\rightarrow (PGL_n)^{m+1}=:\mathcal{Y}$ with $G:=\varphi(H)$ Zariski-dense in $\mathcal{Y}.$ The group $\tilde{\mathcal{F}}$ is then the first simple factor of $\mathcal{Y}.$ Write $\mathcal{Y}=\tilde{\mathcal{F}}\times \prod^m_{i=1}{\mathcal{Y}_i}$ where $\mathcal{Y}_i:=(\varphi_{c_i}(H))^c$. Let $\mathcal{H}$ be a direct factor of $\mathcal{Y}$. If $\mathcal{F}$ is not a factor of $\mathcal{Y}_1$ then $G\bigcap \mathcal{Y}_1=\{1\}$. If $\mathcal{H}=\tilde{\mathcal{F}}\times \prod_{i\in I}{\mathcal{Y}_i}$ for $I$ a subset of $\{1,\ldots, m\}$, then $G/G\bigcap \mathcal{H}\simeq \prod_{i\notin I}{\varphi_{c_i}(H)}$. Thus $\mathcal{Y}$ itself is the smallest normal subgroup $\mathcal{H}$ of $\mathcal{Y}$ such that the projection of $G/G\bigcap \mathcal{H}$ on $\tilde{\mathcal{F}}$ is equal to that of $G$. So $\mathcal{F}=\mathcal{Y}$.

We conclude this example by pointing out that it is not necessary that all simple factors of $\mathcal{F}$ are of the same type. For example, replacing $H$ by $\varphi^{-1}_{c_1}(PSO_n(\mathbb{F}_p))$ and then proceeding as above we get that $\mathcal{Y}_1\simeq PSO_n$.

\subsection{}Assume now that $G$ is primitive. Let $C$ be the center of $G$ (so that $T=C$ in our case). Let $\mathcal{Y}:=(G^c)^{\circ}$. Since $G$ is primitive $\mathcal{Y}$ is semi-simple. Set $\mathcal{Y}=\mathcal{Y}_1\ldots\mathcal{Y}_s$, an almost direct product of almost simple groups. The group $\tilde{N}:=N_{GL_n}(\mathcal{Y})/{Z_{GL_n}(\mathcal{Y})}$ consists of permutations of factors and of outer (that is, graph) automorphisms of the $\mathcal{Y}_i$. By (14.4) $G^{\circ}$ contains normal subgroups $P$ and $F$ such that (14.1)(a)-(d) hold with $T=C$. Thus (since $P^c F^c=(G^c)^{\circ}$ and $G^{\circ}/{PF}$ consists of diagonal automorphisms of $P$) it follows that for $N:=N_{GL_n(k)}(PF)/{Z_{GL_n(k)}(PF)}$ we have $|N|\leq |\tilde{N}|\cdot |\text{Outer diagonal automorphism group of $P$}|.$ As in (4.5.4) this implies

\subsubsection{}$|N|\leq n^2\,\log{n}$\\

It is, therefore, now remains to study $Z:=Z_G(PF)$. This latter group is finite and, since $G$ is primitive, it is completely reducible. An argument similar to ones we used in Sections 10 and 11 shows that $Z$ contains a normal subgroup $L$ centrally isomorphic to a direct product of simple finite groups of Lie $p$-type such that $Z/L$ satisfies the conclusions of (11.1).

Now a repetition of arguments of Sections 12 and 13 yields (14.1) in complete generality.   

\newpage
\section{\textbf{Maximality of some finite linear subgroups of $GL_n(\mathbb{C}).$}}  

Let $\varphi_m: \textnormal{Alt}_m\rightarrow GL_{m-1}(\mathbb{C})$ be the non-trivial component of the transitive permutation representation of $\textnormal{Alt}_m$ on $m$ letters and $\psi_m:\textnormal{Alt}_m\rightarrow GL_{\frac{m\,(m-3)}{2}}(\mathbb{C})$ the representation of $\textnormal{Alt}_m$ such that the permutation representation (of degree $m\,(m-3)/2$) of $\textnormal{Alt}_m$ on unordered pairs of $m$ distinct letters is equivalent to $\id\oplus \varphi_m\oplus \psi_m.$ 

We shall call a finite subgroup $G$ of $GL_n(\mathbb{C})$ \textit{nearly maximal} if for any finite subgroup $H$ of $GL_n(\mathbb{C})$ the inclusion $H\supseteq G$ implies that $H\subseteq N_{GL_n(\mathbb{C})}(G)$. In the cases we consider below near-maximality of $G$ means that $N_{GL_n(\mathbb{C})}(G)$ is modulo its center, a maximal finite subgroup of $GL_n(\mathbb{C}).$ But this interpretation does not hold in other examples. \\

\prop{} For any $r\in \mathbb{N}$ there exists $n_1=n_1(r)\in \mathbb{N}$ such that if $G:=\otimes^{r}_{i=1}{\varphi_{m_i}(\textnormal{Alt}_{m_i})}\subseteq GL_{\prod^{r}_{i=1}(m_i-1)}(\mathbb{C})$ and $m_i\geq n_1(r), i=1, \ldots, r,$ then $G$ is nearly maximal in $GL_{\prod^{r}_{i=1}(m_i-1)}(\mathbb{C})$.\\ 
\eprop
                                                                                                                                                       \prop{} There exists $n_2\in \mathbb{N}$ such that if $G:=\psi_m(\textnormal{Alt}_m)\subseteq GL_{m\,(m-3)/2}(\mathbb{C})$ and $m\geq n_2$ then $G$ is nearly maximal in $GL_{m\,(m-3)/2}(\mathbb{C})$.\\ 
\eprop                                                                                                                                                 

\prop{} If $G:=\otimes^{r}_{i=1}{\varphi_{m_i}(\textnormal{Alt}_{m_i})}\subseteq GL_{\sum^{r}_{i=1}{(m_i-1)}}(\mathbb{C})$ and $m_i\geq 10, i=1, \ldots, r,$ then $G$ is nearly maximal in $GL_{\sum^{r}_{i=1}{(m_i-1)}}(\mathbb{C}).$\\ 
\eprop

\prop{} There exists $n_4\in \mathbb{N}$ such that if $G:=\varphi_m(\textnormal{Alt}_m)\oplus \psi_n(\textnormal{Alt}_n)\subseteq GL_{(m-1)+n\,(n-3)/2}(\mathbb{C})$ and $m\geq n_4, n\geq n_4$ then $G$ is nearly maximal in $GL_{(m-1)+n\,(n-3)/2}(\mathbb{C}).$\\
\eprop

\subsection{}Remark. (15.3) and (15.4) give examples of nearly maximal reducible finite linear groups.\\

\subsection{}\textit{Proof of (15.1).} First, if $r=1$, then $\tilde{G}:=\varphi_{m_1}(\Sym_{m_1})$ is a group generated by reflections. This $G$ is known not to be maximal for $m_1=9$ (for $\Sym_9$ is contained in the Weyl group of type $E_8$). When $m_1> 9$ consider $H\supseteq G, H\subseteq GL_m(\mathbb{C}), H$ finite and then replace $H$ by $\tilde{H}:=<H, \tilde{G}>.$ $H$ contains then a normal subgroup generated by reflections and from the classification of finite groups generated by reflections we see that $G$ is nearly maximal if $m_1\geq 10$ (so that $n_1(1)=10$). 

Let us now choose $n_1=n_1(r)$ so that $n_1(r)\geq 49$ and $(n_1!/2)^r > (2\,(n_1-1)^r+1)^{2\,\log_3{(2\,(n_1-1)^r+1)}}$. Clearly, such $n_1(r)$ exists and it is easy to see that $n_1(r)> n_1(a)$ if $a<r$.

Let now $m:=\prod^r_{i=1}{(m_i-1)}$ and let $H\subseteq GL_m(\mathbb{C})$ be a finite group such that $H\supseteq G.$ $H$ is primitive since so is $G$. Let $\bar{S}$ be the socle of $H/{\textnormal{center}}$ and $S$ its preimage in $H$. Then, as in Section 11, $S$ is a central product of centrally simple groups $G_1, \ldots, G_t,$ extraspecial groups $E_1, \ldots, E_s$ and of the center $C$ of $H$. We have, of course, $G\subseteq S$ (need to be proven) whence at once $s=0$. We show now that each $\varphi_{m_i}(\textnormal{Alt}_{m_i})\otimes \id$ is contained in some $G_j.$ Suppose that is not so. Assume, for definiteness, that $G_1$ is such that the projections of $\varphi_{m_1}(\textnormal{Alt}_{m_1})\otimes \id$ both on $G_1$ and on $G_2\ldots G_t$ are non-trivial. The representation of the central product of $G_1$ and $G_2\ldots G_t$ on $k^m$ is equivalent to the tensor product of representations $\pi\otimes\omega$ of the two factors. Restricted to $\varphi_{m_1}(\textnormal{Alt}_{m_1})\otimes \id$ it implies that $\varphi_{m_1}$ is a tensor product of two representations of $\textnormal{Alt}_{m_1}$. But since $\varphi_{m_1}$ is a non-trivial representation of smallest dimension, it can not be a tensor product. Thus each $\varphi_{m_i}(\textnormal{Alt}_{m_i})\otimes \id$ is contained in some $G_i$. Thus $\{1,\ldots,r\}=\bigcup^{t}_{i=1}{I_i}$ with $I_i\bigcap I_j=\varnothing$ if $i\neq j$ and $G_j\supseteq\otimes_{i\in I_j}{\varphi_{m_i}(\textnormal{Alt}_{m_i})\otimes \id}$.

Let us now argue by induction on $r$. The case $r=1$ was dealt with before. Make the inductive assumption. Since $n_1(r)\geq n_i(|I_j|)$ for $j=1,\ldots, t$ it follows that either $|I_j|=1$ for $j=1,\ldots, t$ (and then $t=r,$ whence $H\subseteq N_{GL_{m}(\mathbb{C})}$) or $t=1.$ Our choice of $n_1(r),$ together with (6.1), (7.1) implies then that $G_1\simeq \textnormal{Alt}_d$ for some $d\leq m+1$. Let $\varphi_i$ denote the imbedding of $\textnormal{Alt}_{m_i}$ in $\textnormal{Alt}_d$. Let $\Omega_1, \ldots, \Omega_b$ be different orbits of $\varphi_1(\textnormal{Alt}_{m_1})$ on $\{1,\ldots,d\}$. Let $\{1, \ldots, d\}=\bigcup^c_{i=1}{J_i}$ so that the orbits $\Omega_{\alpha}$ and $\Omega_{\beta}$ are equivalent if and only if $\alpha, \beta\in J_i$ for some $i$. Then $Z_{\textnormal{Alt}_d}(\varphi_1(\textnormal{Alt}_{m_1}))\simeq \prod^c_{i=1}{\Sym{J_i}}$ where each $\Sym{J_i}$ permutes the orbits $\Omega_{\alpha}, \alpha\in J_i.$ Since $\varphi_i(\textnormal{Alt}_{m_i})\subseteq Z_{\textnormal{Alt}_d}(\varphi_1(\textnormal{Alt}_{m_1}))$ for $i=2, \ldots, r$ and by the inductive assumptions we must have $\varphi_i(\textnormal{Alt}_{m_i})=\Sym{J_{\alpha(i)}}$ for $i=2,\ldots, r$ and an appropriate $\alpha(i)=1,\ldots, c$.

Note that $\textnormal{Alt}\,J_i$ acts trivially on $\bigcup_{\alpha\notin J_i}{\Omega_{\alpha}}.$ This implies that each $\textnormal{Alt}\,J_i$ has at most two types of orbits on $\{1,\ldots, d\}$; and if exactly two then one type is trivial. By symmetry this, therefore, holds for all $\varphi_i(\textnormal{Alt}_{m_i})$, and, in particular, for $i=1$. Thus $c\leq 2, r\leq 3$. Suppose there is a trivial orbit, say $\Omega_1,$ of $\varphi_1(\textnormal{Alt}_{m_1})$. Let $\Omega_1\in J_1.$ Then $Z_{\textnormal{Alt}_d}(\textnormal{Alt}\,J_1)=\textnormal{Alt}(\{1,\ldots, d\}-J_1)$ whence again by inductive assumption we must have $Z_{\textnormal{Alt}_d}(\textnormal{Alt}\,J_1)$ is one of the $\varphi_i(\textnormal{Alt}_{m_i}).$ Thus we can assume in this case (when there is a trivial orbit) that $c=r=2$ and $\varphi_i(\textnormal{Alt}_{m_i}), i=1,2,$ acts through the natural representation on its non-trivial orbit. We have thus $d=m_1+m_2$. Then $m=(m_1-1)\,(m_2-1)\leq ((d-2)/2)^2$. Thus we are dealing with a representation of $\textnormal{Alt}_d$ of dimension $\leq ((d-2)/2)^2.$ By R. Rasala [p. 132, Result 2] this implies (since $d\geq 2\,n_1(2)\geq 20$) that $m\leq d-1$, i.e. $(m_1-1)\,(m_2-1)\leq m_1+m_2-1$ or $(m_1-2)\,(m_2-2)\leq 2.$ This latter inequality is false for $m_1, m_2\geq 10$. Thus our current assumption that there are trivial orbits is false as well. Thus there are no trivial orbits and therefore $J_1=\{1,\ldots, d\}$. We assume $r > 1$, so that $r=2$. Since $\varphi_2(\textnormal{Alt}_{m_2})$ acts as $\textnormal{Alt}\,J_1$ on the non-trivial orbits of $\varphi_1(\textnormal{Alt}_{m_1})$ we see that non-trivial orbits of $\varphi_i(\textnormal{Alt}_{m_i})$ are of length $m_i$ for $i=1,2$. Thus $d=m_1\cdot m_2.$ Since $m=(m_1-1)\,(m_2-1)< d-1$ it follows that $\textnormal{Alt}_d$ can have no representation of dimension $m$ whence a contradiction in this case. This concludes the proof of (15.1). $\Box$\\

\subsection{}\textit{Proof of (15.2).} Again, take $n_2\geq 49$ and such that $(n_2!)/2>$$(n_2\,(n_2-3)+1)^{2\,\log_3{(n_2\,(n_2-3)+1)+1}}$. 
Clearly, such $n_2$ exists. Set $n:=m\,(m-3)/2.$ Let $H\subseteq GL_n(\mathbb{C}),$ $H$ finite, $H\supseteq G$. Since $G$ is primitive so is $H$. Let $S$ be the preimage in $H$ of the socle of $H/\textnormal{center}.$ Then $S$ is a central product of centrally simple groups $G_1, \ldots, G_t,$ extraspecial groups $E_1, \ldots, E_s$ and of the center $C$ of $H$. We have again that $s=0.$ Clearly, $G\subseteq G_1\ldots G_t$. Let $G_{(i)}$ be the projection of $G$ on $G_i, i=1,\ldots t.$ The representation of $G_1\ldots G_t$ on $k^n$ is a tensor product $\otimes^t_{i=1}{\pi_i}$ of representation of the $G_i$. Therefore our representation $\psi_n$ of $G\simeq \textnormal{Alt}_m$ is a tensor product of the $\pi_i|G_{(i)}$. But since by R. Rasala \cite{?} the smallest non-trivial representation of $G$ has dimension $m-1$ and since $(m-1)^2 > m\,(m-1)/2$ it follows that $\psi_m$ is not a tensor product whence $t=1$. 

Then by (6.1), (7.1) and because of or choice of $n_2$ we see that $G_1\simeq \textnormal{Alt}_d$ for some $d$. We have since $\textnormal{Alt}_d\supseteq \textnormal{Alt}_m$ that $d\geq m$. If $d > m$ then $d\,(d-3)/2 > n$ and hence by R. Rasala [\cite{?}, Result 2] we have $n=d-1,$ i.e. $d=m\,(m-3)/2+1.$ 

Consider the action of $\textnormal{Alt}_m$ on $\Omega:=\{1,\ldots d\}$. If $\textnormal{Alt}_m$ has an orbit $\Omega_1\neq \Omega, \Omega_1\neq\varnothing,$ on $\Omega$ then $\textnormal{Alt}_m\subseteq \textnormal{Alt}\,\Omega_1\times \textnormal{Alt}\,(\Omega-\Omega_1)$ and the restriction of our representation of $\textnormal{Alt}_d$ (of dimension $d-1$) to $\textnormal{Alt}\,\Omega_1\times \textnormal{Alt}\,(\Omega-\Omega_1)$ would be irreducible. But this is not so if $\Omega_1\neq \Omega, \Omega_1\neq\varnothing.$ Thus $\textnormal{Alt}_m$ is transitive on $\Omega.$

Since any primitive permutation representation of $\textnormal{Alt}_m$ of degree $> m$ has degree $\geq m\,(m-1)/2$ (see ?) it follows from $m\,(m-3)/2+1 < m\,(m-1)/2$ that the orbit of $\textnormal{Alt}_m$ on $\Omega$ has length divisible by $m$. Thus $m\,(m-3)/2+1=r\,m$ for some $r\in \mathbb{N}$. Clearly $r > (m-3)/2.$ Thus $2=m\,(2\,r-(m-3))$ with $2\,r-(m-3) > 1.$ This is clearly impossible for $m\geq 49.$ This is a contradiction with the assumption $d>m$. Thus $d=m$ and (15.2) is proved. 

\subsection{}\textit{Proof of (15.3).} For $r=1$ (15.3) reduce to (15.1). So assume $r\geq 2.$

Let $m:=\sum^r_{i=1}{m_i-1}$ and let $H\subseteq GL_m(\mathbb{C})$ be a finite group such that $H\supseteq G$. Replace $H$ by $<H, \oplus^{r}_{i=1}{\varphi_{m_i}(\Sym_{m_i})}>$ and let $M$ be the normal subgroup of $H$ generated by $\oplus^{r}_{i=1}{\varphi_{m_i}(\Sym_{m_i})}.$ Set $V:=k^m$ and let $V=\oplus^t_{i=1}{V_i}$ be the decomposition of $V$ into $M$-simple modules. Then $M_i=M|V_i$ is generated by reflections and is irreducible. Then by the classification of finite irreducible groups generated by reflections (see, e.g., ?) and since $\dim{V_i}\geq \min{(m_j-1)}\geq 9$ we have that $M_i$ is a Weyl group of one of the types $A_s, B_s,$ or $D_s$ with $s=\dim{V_i}$. We can assume that $M_i=M, V_i=V.$ Thus we obtain in case $A_s$ that $\textnormal{Alt}_{m+1}$ contains a direct product of $\textnormal{Alt}_{m_i}$ with $\sum{m_i}=m+r,$ and in the case of $B_s$ and $D_s$, by taking quotient of $M$ by its radical, that $\textnormal{Alt}_m$ contains a direct product of $\textnormal{Alt}_{m_i}$ with $\sum{m_i}=m+r.$ This situation is easily seen to be impossible unless $r=1, M$ is of type $A_{m+1}.$ Returning to our original $M$ this means that each $V_i$ is irreducible for exactly one $\varphi_{m_j}(\textnormal{Alt}_{m_j})$ whence our claim.\\

\subsection{}\textit{Proof of (15.4).} Let us take $n_4\geq \max{\{n_1(1), n_2\}}$ and such that $(n_4!/2)^2 > (2\,n_4+n_4\,(n_4-3)+1)^{2\,\log_3{(2\,n_4+n_4\,(n_4-3)+1)+1}}.$ 

Let $d=m-1+n\,(n-3)/2$ and let $H\subseteq GL_d(k), H$ finite, $H\supseteq G$. If $H$ is reducible then the irreducible components clearly will have dimensions $m-1$ and $n\,(n-3)/2$ and then $H\subseteq N_{GL_d(k)}(G)$ by (15.1) and (15.2). If $H$ is irreducible it is easily seen to be primitive. Let $S$ be the preimage in $H$ of the socle of $H/\textnormal{center}.$ As before one shows that $S$ must be centrally simple and then, by (6.1), (7.1), $S\simeq \textnormal{Alt}_r$ for some $r$. One has at once, looking at length of orbits of $\textnormal{Alt}_m$ and $\textnormal{Alt}_n$ on $\Omega:=\{1,\ldots, r\}$ that $r\geq m+n.$ Using again R. Rasala [\cite{?}, Result 2] and noting that $m-1+n\,(n-3)/2 < (m+n)\,(m+n-3)/2\leq r\,(r-3)/2$ we see that $r=d+1=m+n\,(n-3)/2.$ Let $\Omega_1,\ldots, \Omega_t$ be different orbits of $\textnormal{Alt}_m$ on $\Omega$. Since the representation of $\textnormal{Alt}_r$ on $k^d$ is the smallest one and since its restriction to $\textnormal{Alt}_m$ is again the smallest one taken once plus a number of trivial ones, we have that one orbit, say $\Omega_1,$ is of length 1 and the rest are of length 1. Then $Z_{\textnormal{Alt}_r}(\textnormal{Alt}_m)=\textnormal{Alt}_{r-m}$ whence $\textnormal{Alt}_n\leq \textnormal{Alt}_{r-m}$ which contradicts near maximality of $\psi_n(\textnormal{Alt}_n)$ (i.e. contradicts (15.2)) unless $r-m=n$. But this latter variant is impossible as $n< n\,(n-3)/2.$   

\newpage
\section{\textbf{Pairs of centrally simple irreducible embedded subgroups.}}

Let $k$ be an algebraically closed field of characteristic exponent $p:=p(k).$ 

\subsection{}A pair $H\subseteq G$ of centrally simple non-commutative subgroups of $GL_n(k)$ will be called here \textit{tight} if both $H$ and $G$ are irreducible. The following pairs of subgroups of $GL_n(\mathbb{C})$ are examples\\

(a) $G=\textnormal{Alt}_{n+1}$ and $H$ a doubly transitive group of permutations of $n+1$ letters (see W. Feit [\cite{?}, $\S 9.1$] for concrete $n$ and $H$) in the $n$-dimensional representation of $\textnormal{Alt}_{n+1}$,\\

(b) $G=\textnormal{Alt}_{m\,d}, m,d > 1, H=\textnormal{Alt}_{m\,d-1}$, and the representation of $G$ corresponds to the rectangular Young diagram with height $m$ and width $d$ (this example was explained to me by A. Regev), $n$ is the dimension of the corresponding representation,\\

(c) $n=(q^m-1)/2, G=Sp_{2\,m}(\mathbb{F}_q), H=SL_2(\mathbb{F}_{q^m}), q$ a power of an odd prime.\\

To explain (a) note that since $H$ is doubly transitive its permutation representation is a direct sum of the trivial one and an irreducible one of dimension $n$ (see, e.g., ?). Since the same holds for $\textnormal{Alt}_{n+1}$ we see that the restriction of the $n$-dimensional component of the permutation representation of $G$ on $n+1$ letters is irreducible. \\

To explain (b) note that by the branching rule (see G. James \cite{17}) the restriction of an irreducible  representation of $\Sym_r$ with Young diagram $T$ to $\Sym_{r-1}$ consists of components whose Young diagrams are obtained from $T$ by removing exactly one square. For a rectangular diagram there is just one way to remove a square which leads to a Young diagram. Thus the restriction of our representation from $\Sym_r$ to $\Sym_{r-1}$ remains irreducible. That the same holds for $\textnormal{Alt}_r$ and $\textnormal{Alt}_{r-1}$ follows from, e.g., G. James \cite{17}.\\

For (c) we note that $SL_2(\mathbb{F}_{q^m})$ can be considered as acting on $\mathbb{F}_{q^m}^2=\mathbb{F}^{2\,m}_q.$ It preserves in this action a symplectic form whence embedding $SL_2(\mathbb{F}_{q^m})\subseteq Sp_{2\,m}(\mathbb{F}_q)$. $Sp_{2\,m}(\mathbb{F}_q)$ has an irreducible  representation, say 4 of dimension $n:=(q^m-1)/2$ (see, e.g. ?). Since $(q^m-1)/2$ is also the smallest degree of a non-trivial representation of 
$SL_2(\mathbb{F}_{q^m})$ we see $\varphi|SL_2(\mathbb{F}_{q^m})$ must be irreducible.

The following result says that example (c) is, in some sense, typical.\\

\th{} Let $G\supseteq H$ be a tight pair in $GL_n(k).$ Suppose that $G$ is of Lie $l$-type and $H$ is of Lie $r$-type with $l\neq p.$ Then $l=r$ unless $|H/{\textnormal{center}}|\leq 1.5\cdot 10^{33}.$\\
\eth

We also have\\

\th{} Let $G\supseteq H$ be a tight pair in $GL_n(k).$ Suppose that $G$ is of Lie $l$-type with $l\neq p$ and $H$ is a covering group of $\textnormal{Alt}_m$. Then, unless $m\leq 32$, $G$ is classical and the image of $\textnormal{Alt}_m$ in the natural representation of $G$ over $\bar{\mathbb{F}}_l$ is the smallest non-trivial irreducible representation of $\textnormal{Alt}_m$.\\
\eth

\subsection{} \textit{Proof.} Suppose $r\neq l, G/{\textnormal{center}}\simeq$$^{c}\bar{X}_a(m^c), m^c$ a power of $l$. Take an irreducible non-trivial representation $G\rightarrow GL_d(\bar{\mathbb{F}}_l)$ with\\

\begin{tabular}{cccccccccc}
$d=d(X_a)$ & $\leq$ & $a+1$ & $2\,a$ & $2\,a+1$ & $7$ & $26$ & $27$ & $56$ & $248$ \\
if $X_a $ & $=$ & $A_a$ & $C_a, D_a$ & $B_a$ & $G_2$ & $F_4$ & $E_6$ & $E_7$ & $E_8$ \\ 
\end{tabular}
\ \\

Since $|H/{\textnormal{center}}|\leq f(d)$ where, as before, $f(d):=(2\,d+1)^{2\,\log_3{(2\,d+1)+1}}$ we see that only a finite number (independent of $p$) of $H$ can join tight pairs with $G$ having $d\leq 213$ or $d=248$. (This excludes all exceptional $G$). Assume therefore that $d\geq 214, d\neq 248.$ From Table T4.4  we see that $b\geq a$. Since $^{2}B_2(2)$, $^{2}G_2(3)$, $^{2}F_4(2)$ are excluded we also have $m\geq 2$. Thus by (4.4.2) $n\geq (2^a-1)/2.$ Combining this with $a\geq (d-1)/2$ we obtain $n\geq (2^{(d-1)/2}-1)/2=2^{(d-3)/2}-0.5$. One easily verifies that $2^{(d-3)/2}-0.5>f(d)$ for $d\geq 214$. Thus $n > |H/{\textnormal{center}}|.$

By C. Curtis and I. Reiner \cite[(53.16)]{CR} degrees of ordinary representations of $H$ divide $|H/{\textnormal{center}}|.$ Since every representation in characteristic $p > 0$ is a composition factor of a reduction $\mod{p}$ of an ordinary representation we see that the degrees of all representations of $H$ are $\leq |H/{\textnormal{center}}|$. Thus if $d\geq 214, d\neq 248,$ our representation of $H$ can not be irreducible, a contradiction. Note, finally, that $f(248)\leq 1.5\cdot 10^{33}.$

\subsection{} Remark. Actually for almost all $H$ we have $|H|\leq f(d).$ Then instead of estimate
$$
\text{(degree of a representation)}\leq |H/\textnormal{center}| 
$$
we can use: 
$$
\text{(degree of a representation)}\leq \sqrt{|H|}
$$

Then it is sufficient to take $d\geq 78, d\neq 248.$ Additional restrictions are obtained if one takes into account (4.4.3)(b) but this only reduces 78 to 68. 

\cor{} Let $p$ and $l$ be odd primes, $p\neq l$. Let $\varphi: SL_2(\mathbb{F}_l)\rightarrow GL_{(l-1)/2}(\bar{\mathbb{F}}_p)$ be an irreducible representation. Then $\varphi(SL_2(\mathbb{F}_l))$ is maximal in $Sp_{(l-1)/2}(\mathbb{F}_p)$ if $l\equiv 1(\mod{?})$ in $SO_{(l-1)/2}(\mathbb{F}_p).$
\ecor

\newpage
\section{\textbf{Splitting fields.}}

Let $k$ be an algebraically closed field of characteristic $p>1.$ For a finite group $H$ we denote by $e(H)$ its exponent, the smallest integer $e$ such that $h^e=1$ for all $h\in H.$ By C. Curtis and I. Reiner \cite[(70.24)]{CR} every representation of $H$ over $\mathbb{F}_p(\sqrt[e]{1})(\leq \mathbb{F}_{p^{\varphi(e)}},$ where $\varphi(e)$ is the Euler function) is absolutely irreducible. In other words if $H\subseteq GL_n(k)$ then a conjugate of $H$ is contained in $GL_n(\mathbb{F}_p(\sqrt[e]{1})).$ Let $e_0$ be the least common multiple of the exponent of the universal central extensions of the sporadic simple groups and ??? groups of Lie type having ??? central extensions (see (4.3.3) for a lit). 

For a field $K$ and  a finite group $G$ we denote by $I(KG)$ a set of representatives of the equivalence classes of the irreducible $K$-representations of $G$. For $f\in I(KG)$ we denote by $K_0(f)$ the extension of the prime field $K_0$ of $K$ given by $K_0(f):=K_0(\Tr{f(G)})$, meaning the field generated by the character values of the $g\in G$ under $f$.\\

\th{} Let $G$ be a centrally simple finite group not isomorphic to a group of Lie $p$-type and $f: G\rightarrow GL_n(k)$ an irreducible representation. Then $[\mathbb{F}_p(\Tr{f}):\mathbb{F}_p]\leq \max{\{2, \varphi(e_0), n^2\}}.$

If $n$ is sufficiently large ($n\geq \sqrt{|F_1|}$, e.g. $n\geq 9\cdot 10^{26}$ would suffice) then $[\mathbb{F}_p(\Tr{f}):\mathbb{F}_p]\leq \max{\{2, n^2\}}.$ \\
\eth

\cor{} Let $H_0=$ $^{c}\bar{X}_a(p^{r\,c})$ be a classical finite simple group of Lie $p$-type and let $H\supseteq H_0, H\subseteq \Aut{H_0}.$ Suppose that $r > \max{\{2, \varphi(e_0), (2\,a+1)^2\}}$. Let $M$ be a maximal subgroup of $H$. Then either\\

\ \ \ \ (a) $M$ is one of the groups on the M. Aschbacher \cite{1} list $\underline{C}_H$\\

or (b) the socle of $M$ is simple of Lie $p$-type\\

\ecor

\subsection{} Remark. Actually (17.2) can be sharpened by adding divisibility properties. For example, let $M_1, \ldots, M_{N(n)}$ be the list of universal central extensions of simple groups which can have a representation of degree $n$. 

For each such $M_i$ and $f\in I(\mathbb{C}M_i)$ let $a_f$ be the degree of a maximal cyclic subextension of $\mathbb{Q}(\Tr{f}).$ Let $a_1, \ldots, a_{R(n)}$ be the collection of all numbers $a_f$ for the $M_i, i=1, \ldots, N(n),$ and $f\in I(\mathbb{C}M_i).$ Then the conclusion of (17.2) holds if $r\geq \max_{1\leq i\leq R(n)}{\{LCD(r, a_i), 2\}}.$\\

\textit{Proof of (17.3).} M. Aschbacher [\cite{?}, ?] states that unless $M$ belongs to $\underline{C}_H$ the socle $G$ of $M$ is simple and can not be written in field smaller than $\mathbb{F}_{p^r}.$ Since $H_0$ can be embedded into $GL_n(\bar{\mathbb{F}}_p), n\leq 2\,a+1$ (see, e.g., beginning of (16.3)) our claim follows from (17.1). $\Box$\\

\subsection{} \textit{Proof of (17.1).} First let us consider an ordinary representation $g: G\rightarrow GL_n(\mathbb{C})$. Suppose $K:=\mathbb{Q}(\Tr{g})$.

\newpage
\begin{center}
A. \textbf{Computational lemmas.}\\
\end{center}

Let $f(x):=(2\,x+1)^{2\,\log_3{(2\,x+1)}}, x\geq 0$, and $\tilde{f}(n):=f(n)$ if $n > 2, \tilde{f}(2):=60.$\\

\textbf{Lemma A1.} \textit{If $x, y, t\in \mathbb{N}$ then}\\

\textit{(a) $\tilde{f}(x)\,\tilde{f}(y)\leq \tilde{f}(x\,y)$ if $x\geq 2, y\geq 2,$}\\

\textit{(b) $t!\,(\tilde{f}(x))^t\leq \tilde{f}(x^t)$ if $x\geq 2, t\geq 1,$}\\

\textit{(c) $1.025\,\tilde{f}(4)\,\tilde{f}(y)\leq \tilde{f}(4\,y)$ if $y\geq 2$,}\\

\textit{(d) $(1.025\,\tilde{f}(4))^t\,t!\leq \tilde{f}(4^t)$ if $t\geq 2.$}\\ 

\prf{} Consider $F(x, y):=\ln{f(x)}+\ln{f(y)}-\ln{f(x\,y)}.$ We have 
$$
\frac{\partial F}{\partial x}=\left(\frac{8}{\ln{3}}\right)\,\frac{\ln(2\,x+1)}{(2\,x+1)}+\frac{2}{(2\,x+1)}-\left(\frac{8}{\ln{3}}\right)\,\frac{\ln(2\,x\,y+1)}{(2\,x\,y+1)}-\frac{2\,y}{(2\,x\,y+1)}
$$

Note that $1/{(2\,x+1)} < y/{(2\,x\,y+1)}=1/{(2\,x+1/y)}$ for $x\geq 0, y> 1,$ i.e., $\partial F(x, y)/{\partial x} < 0$ for $x\geq 0, y > 1$. Hence $F(x, y)\leq F(2, y)$ for $y > 1$. Then 
$$
\frac{d\,F(2, y)}{d\,y}=\left(\frac{8}{\ln{3}}\right)\,\frac{\ln(2\,y+1)}{(2\,y+1)}+\frac{2}{(2\,y+1)}-\left(\frac{16}{\ln{3}}\right)\,\frac{\ln(4\,y+1)}{(4\,y+1)}-\frac{4}{(4\,y+1)}
$$ 

and the inequality $2/{(2\,y+1)} < 4/{(4\,y+1)}$ for $y > 1$ whence $F(2, y) < F(2, 4).$ Now one checks directly that $F(2, 4) < -0.87$ thus proving (a) if $x$ or $y \geq 4$. One checks directly that $\tilde{f}(2)\,\tilde{f}(2)=3600 < f(4)=59,049, \tilde{f}(2)\,f(3)=414000 < f(6)= 2,067,423, f(3)\,f(3)=47,610,000 < f(9)=1.35\cdot 10^8.$ This proves (a).\\

To prove (b) it is sufficient to assume that $t\geq 2.$ We consider as above $F(t, x)=\ln{t!}+ t\,\ln{f(x)}- \ln{f(x^t)}.$ We have
$$
\frac{\partial F(t, x)}{\partial x}=\frac{t\,f'(x)}{f(x)}-\frac{t\,x^{t-1}\,f'(x^t)}{f(x^t)} < 0
$$ 

as above. Thus $F(t, x) < F(t, 3)$ for $x> 3, t\geq 1$. We consider $t$ now as a real variable. 
$$
F(t, 3)=\ln{\Gamma(1+t)}+t\,\left(\left(\frac{2}{\ln{3}}\right)\,\ln^2{7}+\ln{7}\right)-\left(\frac{2}{\ln{3}}\right)\,\ln^2{\left(2\cdot 3^t+1\right)}-\ln{\left(2\cdot 3^t+1\right)}.
$$

Write $2\cdot 3^t=3^{t+\log_3{2}} \leq 3^{t+0.63}$ and then $\ln{(2\cdot 3^t+1)}\leq \ln{(2\cdot 3^t)}\leq \ln{3}\,(t+0.63).$ Then using the Stirling inequality for the $\Gamma$-function (see ?) we get 
$$
F(t, 3)\leq (t+1/2)\,\ln{t} -t+0.92+1/{12\,t}+8.84\,t-2\,\ln{3}\,(t+0.63)^2-\ln{3}\,(t+0.63) \leq
$$
$$
\leq (t+1/2)\,\ln{t}-2.19\,t^2+4\,t-0.6=:h(t)
$$

(we replaced $1/{12\,t}$ by $1/{24}$ since $t\geq 2$).

Now
$$
h'(t)=\ln{t}+(t+1/2)/t-4.38\,t+4 \leq \ln{t}-4.38\,t+5.25,\ \ \ 
h''(t)=1/t-4.38.
$$

Since $h''(t) < 0$ for $t>2, h'(t)<h'(2)=-2.8<0$ for $t>2$, i.e., $h(t)<h(3)=-4.46$ for $t\geq 3$, i.e., (b) holds if $x\geq 3, t\geq 3.$ One checks $2\,\tilde{f}(3)\,f(3)=9.25\cdot 10^7 < f(9)=1.36\cdot 10^8$ whence (b) holds for $x=3$.

As before we have $F(t, x)\leq F(t, 4)$ for $x\geq 4, t\geq 1.$ So to prove (b) for $x\geq 4$ it suffices to prove it for $x=4$. Instead we will prove the stronger claim (d). We have 
$$
\ln{((1.025\,f(4))^t\cdot \Gamma(1+t)/{f(4^t)})}\leq
$$
$$
\leq(t+1/2)\,\ln{t}-t+0.92+1/{12\,t}+11.011\,t-(2/\ln{3})\,\ln^2{(2^{2\,t+1}+1)}-\ln{(2^{2\,t+1}+1)}\leq
$$
$$
\leq (t+1/2)\,\ln{t}-t+0.96+11.011\,t-(2/\ln{3})\,\ln^2{2\,(2t+1)^2}-(\ln{2})\,(2\,t+1)\leq
$$
$$
\leq(t+1/2)\,\ln{t}-3.48\,t^2+5.15\,t-0.6=:h_1(t).
$$

Now the same argument as above gives $h_1(t) < h_1(2)=-1.72 < 0$ whence (d) holds for $t\geq 2$ whence (b) holds for $x\geq 4$. It remains to check that (b) holds for $x=2$. We have 
$$
\ln{(\Gamma(1+t)\,(60)^t/{f(2^t)})}\leq (t+1/2)\,\ln{t}-0.87\,t^2+0.67\,t-0.6=:h_2(t)
$$

with $h_2(2)=-0.24 < 0$ and $h'_2(t)<0$ for $t\geq 2$ whence as above (b) for $x=2$.

We skip the proof of (c) as it is the same as that of (a).
\epr

Let $g(x):=\Gamma(x+3).$ Recall that $\Gamma(n+3)=(n+2)!$ if $n\in \mathbb{N}$.\\

\textbf{Lemma A2.} \textit{Let $x, y, t\in \mathbb{N}$}\\

\textit{(a) $g(x)\,g(y)\leq g(x\,y)$ if $x, y\geq 2,$}\\

\textit{(b) $t!\,(g(x))^t\leq g(x^t)$ if $x\geq 3, t=1.$}\\

\prf{} We have $g(x)\,g(y)\leq (x+2)!\,\prod_{x+3\leq i\leq x+y+4}{i}=(x+y+4)!$  One easily checks that $x\,y+2\geq x+y+2$ if $x\geq 3$ and $y\geq 3$. If $x=2$ then $g(x)\,g(y)=2\,(y+2)!\geq (2\,y+2)!$ Thus (a) holds if $x\geq 2, y\geq 2$ as claimed. 

To prove (b) take logarithms and consider $F(x, t):=\ln{t!}+t\,\ln{g(x)}-\ln{g(x^t)}.$ 

Thus $F(x, t)=\ln{\Gamma(1+t)}+t\,\ln{\Gamma(x+3)}-\ln{\Gamma(x^t+3)}$. 

Then $\partial F/ {\partial x}=t\,(\ln{\Gamma})'\,(x+3)-t\,x^{t-1}\,(\ln{\Gamma})'\,(x^t+3)$.

Since $(\ln{\Gamma})'(x)$ is strictly increasing for $x\geq 2$ (by \cite{?}) we see that $\partial F/{\partial x} < 0$ for any $t$. 

Thus $F(x, t) < F(3, t)=\ln{\Gamma(1+t)}+t\,(\ln{120})-\ln{\Gamma(3^t+3)}$ 

and $d\,F(3, t)/{d\,t}=(\ln{\Gamma})'\,(1+t)-(\ln{3})\,3^t\,(\ln{\Gamma})'\,(3^t+3)$ which is $<0$ for $t\geq 1$ as before (by ?). Hence $F(x, t) < F(3, t) < F(3, 2)=-7.234.$\\
\epr

\textbf{Lemma A3.} \\

\textit{(a) $f(x)\leq g(x)$ if $x\geq 10,$}\\

\textit{(b) $x^2\,f(x) \leq g(x)$ if $x\geq 14,$}\\

(\textit{c) $x\,f(x) \leq g(x)$ if $x\geq 13.$}\\

\prf{} Taking logarithms of both sides we see that the left ones produce a concave (for $x\geq 2$) function and the right one a convex one. One checks that reverse inequalities hold for $x=2$. Therefore there is just one intersection point of curves representing two sides. Since (a) (resp. (b)) can be checked to hold for $x=10$ (resp. $x=14$) the claim follows from the above remarks.\\  
\epr

\textbf{Lemma A5. } \textit{$2\cdot 2^{x^2+x}\,f(y)\leq f(2^x\,y)$ if $x\geq 2, y\geq 2^{\beta\,x}$ where $\beta=(\log{3}-1)/2=0.29248125.$}\\

\prf{} As usual set (with $y=2^{\alpha\,x}, \alpha$ variable)
$$
F(x):=\ln{2}+\left(2\,x^2+x\right)\,\ln{2}+\frac{2}{\ln{3}}\,\left(\ln{(2^{\alpha\,x+1}+1)}\right)^2+\ln{(2^{(\alpha+1)\,x+1}+1)}-
$$
$$
-\frac{2}{\ln{3}}\,\left(\ln{(2^{(\alpha+1)\,x+1}+1)}\right)^2-\ln{(2^{(\alpha+1)\,x+1}+1)}\leq
$$
$$
\leq \ln{2}+\left(2\,x^2+x\right)\,\ln{2}+\frac{2}{\ln{3}}\,\left(\ln{2^{\alpha\,x+1}}+\ln{(1+2^{-\alpha\,x-1})}\right)^2+\ln{2^{\alpha\,x+1}}+\ln{(1+2^{-\alpha\,x-1})}-
$$
$$
-\frac{2}{\ln{3}}\,\left(\ln{(2^{(\alpha+1)\,x+1}+1)}\right)^2-\ln{2^{(\alpha+1)\,x+1}}\leq
$$
$$
\leq\ln{2}+\left(2\,x^2+x\right)\,\ln{2}+\frac{2}{\ln{3}}\,\left((\ln{2}\,(\alpha\,x+1)+2^{-\alpha\,x-1})\right)^2+(\ln{2})\,(\alpha\,x+1)+2^{-\alpha\,x-1}-
$$
$$
-\frac{2}{\ln{3}}\,\left(\ln{2}\right)^2\,\left((\alpha+1)\,x+1\right)^2-(\ln{2})\,((\alpha+1)\,x+1)
$$

(we used first term approximation to $\ln{(1+\epsilon)}$, i.e. $\ln{(1+\epsilon)}<\epsilon$ for $\epsilon>0$ small).

??? Collecting the terms with like powers of $x$ we have 

$$
F(x)\leq (2\,\ln{2})\,(1-2\,(\ln{2}/\ln{3})\,\alpha-\ln{2}/\ln{3})\,x^2+\ln{2}\,((4/\ln{3})\,\alpha\,2^{-\alpha\,x-1}-\alpha-4\,\ln{2}/\ln{3})\,x+
$$
$$
+(4\,(\ln{2})\,2^{-\alpha\,x-1}+2^{-2\,\alpha\,x-1})/\ln{3}
$$

The coefficient of $x^2$ is non-positive if $\alpha\geq (\log{3}-1)/2$. Assuming $\alpha\geq (\log{3}-1)/2$ and recalling that $x\geq 2$ we get
$$
F(x)\leq \ln{2}\,((4/\ln{3})\,\alpha\,2^{-1.58}-\alpha-4\,\ln{2}/\ln{3})\,x+(4\,(\ln{2})\,2^{-1.58}+2^{-2.16})/\ln{3}\leq
$$
$$
\leq -1.697\,x+1.05
$$

Thus $F(x) < 0$ for $x\geq 2$ and $y=2^{\alpha\,x}\geq 2^{\beta\,x}$ as claimed.\\
\epr

\textbf{A5.1.} Remark. Our argument also shows that $y\geq 2^{\beta\,x}$ is the best one can get.\\

For a simple group $G$ isomorphic to Suz, $\cdot 1, \cdot 2,$ $^{2}\bar{A}_3(9),$ or $D_4(2)$ take the numbers $a_1, a_2$ from Table T7.2 (or from Table TA6 below). Set, as in (7.1), \\

$$F(G, n):=
\begin{cases}
1 & \text{if $n < a_1$}\\
|\Aut{G}| & \text{if $a_1 \leq n < a_2$}\\
f(n) & \text{if $n \geq a_2$}\\
\end{cases}
$$
\ \\

Set $F(H, n):=
\begin{cases}
1 & \text{if $n\leq 3$} \\
f(n) & \text{if $n\geq 4$} \\
\end{cases}$\  for other sporadic simple groups.\\

\textbf{Lemma A6. } \textit{Let $G_1$ and $G_2$ be two sporadic simple groups or isomorphic to $^{2}\bar{A}_3(9).$ Let $F_i(n):=F(G, n), i=1, 2.$ Then}\\

\textit{(a) $F_1(n)\,F_2(m)\leq f(n\,m)$ for $n, m\geq 2, F_1(n)\neq 1, F_2(m)\neq 1$,}\\

\textit{(b) $t!\,(F_1(n))^t\leq f(n^t)$ for $n\geq 2, t\geq 2,$}\\

\textit{(c) $f(m)\,F_1(n)\leq f(n\,m)$ for $m\geq 4, n\geq 2$,}\\

\textit{(d) $\tilde{f}(m)\,F_1(n)\leq f(n\,m)$ for $m\geq 2, n\geq 2$ unless $G_1\simeq \cdot 1$ and $n \leq 37,$}\\

\textit{(e) $1.025\,f(4)\,F_1(n)\leq f(4\,n)$}.\\

\prf{} The claims follow from (A1)(a), (b) if $G_1$ and $G_2$ are not isomorphic to Suz, $\cdot 1, \cdot 2,$ $^{2}\bar{A}_3(9),$ or $D_4(2)$. If $G_2$ is not isomorphic to the 5 above groups then (a) follows from (c). Thus we may and shall assume that both $G_1$ and $G_2$ are isomorphic to Suz, $\cdot 1, \cdot 2,$ $^{2}\bar{A}_3(9),$ or $D_4(2)$. 

Next, if in (a) $n\geq a_2$ where $a_2$ is for $G_1$ from Table T7.2 or TA6 then $F_1(n)=f(n)$ and (a) follows from (c). If in (a) we have further $m\geq (a_2\ \text{for}\ G_2)$ then (a) follows from (A1)(a). Thus to prove (a) (modulo (c)) it is sufficient to check that $F_1(n)\,F_2(m)\leq f(n\,m)$ if $(a_1\ \text{for}\ G_1) \leq n < (a_2\ \text{for}\ G_1), (a_1\ \text{for}\ G_2) \leq m < (a_2\ \text{for}\ G_2),$ i.e., that $|\Aut{G_1}|\cdot|\Aut{G_2}|\leq f((a_1\ \text{for}\ G_1)\cdot(a_1\ \text{for}\ G_2))$. This is verified directly (in the 25 cases). Thus it remains to prove (b), (c), and (d).

There again if $n\geq a_2$ (now there is only one group) then (b), (c), (d) follow from (A1)(a), (b). If $n < a_1$ they become trivial. Thus it remains to verify (b), (c), (d) for $a_1\leq n < a_2$. In these cases $F_1(n)=|\Aut{G_1}|$ and the claim will follow if they hold for $n=a_1.$ 

For (b) we have (as in the proof of (A1)(b)) 
$$
\ln{(\Gamma(1+t)\,|G_1|^t/{f(a_1^t)})}\leq
$$ 
$$
\leq(t+1/2)\,\ln{t}-t+0.96+t\,\ln{|G_1|}-(2/\ln{3})\,(t\,\ln{a_1}+\ln{2})^2-(t\,\ln{a_1}+\ln{2})=
$$
$$
=(t+1/2)\,\ln{t}-(2/\ln{3})\,(\ln^2{a_1})\,t^2+(\ln{|G_1|}-1-(4/\ln{3})\,(\ln{a_1})\,(\ln{2})-\ln{a_1})\,t-0.6
$$

The coefficients of the above expression are given explicitly in Table TA6. One easily establishes, as in the proof of (A1)(b) that the above function is negative for $t\geq 2$ whence (A6)(b).

To prove (c) and (d) consider $r(x):=\ln{(f(x)\,|\Aut{G_1}|/{f(a_1\,x)})}.$ One has 
$$
r'(x)=(2/\ln{3})\,\ln{(2\,x+1)}/{(2\,x+1)} + 2/{(2\,x+1)}-
$$
$$
-(2/\ln{3})\,\ln{(2\,a_1\,x+1)}/{(2\,a_1\,x+1)}-2\,a_1/{(2\,a_1\,x+1)} < 0
$$ 

for $x\geq 2$. Thus $r(x)\leq r(3)$ and from Table TA6 one sees that $r(3)<0$ of $G_1\simeq$ Suz, $\cdot 2,$ $^{2}\bar{A}_3(9),$ or $D_4(2)$. We also have $r(x)\leq r(4)$ and $r(4) < 0$ for $G_1\simeq \cdot 1.$ This proves (c) and (d) except when $m=2$. In this latter case one verifies $60\cdot|\Aut{G_1}|< f(2\,a_1)$ from Table TA6. Finally, it is sufficient to check (e) for $n=a_1$ and this is done directly.
\epr

\ \\

\begin{center}
Table TA.6 (the last 2 rows to be used later)\\
\end{center}
\ \\

\begin{center}
\begin{tabular}{|c|c|c|c|c|c|}

\hline
& \ & \ & \ & \ & \ \\
$G\simeq$ & $^{2}\bar{A}_3(9)$ & $D_4(2)$ & Suz & $\cdot 1$ & $\cdot 2$ \\

\hline
& \ & \ & \ & \ & \ \\
$a_1$ & $6$ & $8$ & $12$ & $24$ & $20$ \\

\hline
& \ & \ & \ & \ & \ \\
$a_2$ & $7$ & $9$ & $18$ & $49$ & $24$ \\

\hline
& \ & \ & \ & \ & \ \\
$|G|$ & $3.26\cdot 10^6$ & $1.74\cdot 10^8$ & $4.48\cdot 10^{11}$ & $4.16\cdot 10^{18}$ & $4.23\cdot 10^{13}$ \\

\hline
& \ & \ & \ & \ & \ \\
$|\Aut{G}|$ & $8\cdot|G|$ & $6\cdot|G|$ & $2\cdot|G|$ & $|G|$ & $|G|$ \\

\hline
& \ & \ & \ & \ & \ \\
$F(G, 2\,a_1)$ & $3.89\cdot 10^9$ & $1.53\cdot 10^9$ & $4.63\cdot 10^{13}$ & $4.16\cdot 10^{18}$ & $1.5\cdot 10^{17}$ \\

\hline
& \ & \ & \ & \ & \ \\
$F(G, 3\,a_1)$ & $7.53\cdot 10^{11}$ & $4.63\cdot 10^{13}$ & $2.62\cdot 10^{16}$ & $5.54\cdot 10^{21}$ & $1.85\cdot 10^{20}$ \\

\hline
& \ & \ & \ & \ & \ \\
$\ln{(f(3)\cdot|\Aut{G}|/{f(3\,a_1)})}$ & $-1.43$ & $-1.86$ & $-1.44$ & $1.65$ & $-6.45$ \\

\hline
& \ & \ & \ & \ & \ \\
$(2/\ln{3})\,\ln^2{a_1}$ & $5.84$ & $7.87$ & $11.24$ & $18.38$ & $16.34$ \\

\hline
& \ & \ & \ & \ & \ \\
$\ln{|\Aut{G}|}-1-\ln{a_1}-$ & $9.77$ & $12.44$ & $17.77$ & $30.68$ & $19.82$ \\
$-(4/\ln{3})\,(\ln{a_1})\,\ln{2}$ & \ & \ & \ & \ & \ \\

\hline
& \ & \ & \ & \ & \ \\
$24\,|\Aut{G}|$ & $6.26\cdot 10^8$ & $2.51\cdot 10^{10}$ & $2.15\cdot 10^{13}$ & $1\cdot 10^{20}$ & $10^{15}$ \\

\hline
& \ & \ & \ & \ & \ \\
$64\,a_1\,f(a_1)$ & $5.3\cdot 10^8$ & $1.93\cdot 10^{10}$ & $3\cdot 10^{12}$ & $7.10\cdot 10^{18}$ & $4.2^{15}$ \\

\hline
\end{tabular}

\end{center}

\newpage
\pagestyle{plain}

\textbf{Lemma A7.} \ \\

\textit{(a) $f(n) \geq 4.796\,n^{2\,\log_3{n}+3.5}$ for $n\geq 2$,\\}

\textit{(b) $n\,f(n)\leq n^{2\,\log{n}+5}$ for $n\geq 4,$\\}

\textit{(c) $n\,f(n)\leq 2\,n^{2\,\log{n}+1}$ for $n\geq 37$, \\ }

\textit{(d) $2\,n^{2\,\log{n}+1}\cdot 2\,m^{2\,\log{m}+1}\leq 2\,(n\,m)^{2\,\log(n\,m)+1}$ for $n, m\geq 2$\\}

\prf{} We have 
$$
(2\,n+1)^{2\,\log_3{(2\,n+1)}}\geq (2\,n)^{2\,\log_3{2\,n}+1}=2^{(2/\log{3})\,(1+\log{n})+1}\cdot n^{(2/\log{3})\,(1+\log{n})+1}=
$$
$$
=2^{(2/\log{3})+1}\cdot n^{2\,\log_3{n}+4/\log{3}+1}\geq 4.796\,n^{2\,\log_3{n}+3.5}
$$

whence (a).

To prove (b) write 
$$
\ln{(f(x))}=(2/\ln{3})\,(\ln{(2\,x+1)})^2+\ln{(2\,x+1)}=
$$
$$
=(2/\ln{3})\,(\ln{2\,x}+\ln{(1+1/2x)})^2+\ln{2\,x}+\ln{(1+1/2x)}\leq
$$
$$
\leq (2/\ln{3})\,(\ln{2}+\ln{x}+1/2x)^2+\ln{2}+\ln{x}+1/2x=:h(x)
$$

Then 
$$
h'(x)=4/\ln{3})\,(1/x-1/2x^2)\,(\ln{2}+\ln{x}+1/2x)+1/x-1/2x^2\leq
$$
$$
\leq (4/\ln{3}+1)/x+4\,\ln{x}/x\ln{3}+2/x^2\ln{3}.
$$

On the other hand, if $r(x):=x^{2\,\log{x}+4}$ then $(\ln{r})'(x)=4\,\ln{x}/\ln{2}+4/x.$ 

We have 
$$
4\,\ln{x}/x\ln{3}-4\,\ln{x}/x\ln{2}=4\,\ln{x}/x(1/\ln{3}-1/\ln{2})\geq 2\ln{x}/x\geq 4/x.
$$

for $x\geq 4$. Therefore for $x\geq 4$
$$
(\ln{r})'(x)-h'(x)\geq 8/x-(4/\ln{3}+1)/x-2/4\,x\,\ln{3}=2.9/x.
$$

Thus $\ln{r}$ increases faster than $h(x)$ for $x\geq 4$. Since $h(4)\leq 11.06$ and $\ln{r(4)}\geq 11.09$ we get that $\ln{r(x)} > h(x)$ for $x\geq 4.$ Since $h(x)\leq \ln{f(x)}$ this implies (b).

To prove (c) we set $s(x):=2\,x^{2\,\log{x}}$. Then $(\ln{s})'(x)=(4/x\ln{2})\,\ln{x}.$ Therefore 
$$
h'(x)-(\ln{s})'(x)=4(1/\ln{3}-1/\ln{2})\,(\ln{x})/x+(4/\ln{3}+1)/x+2/x^2\ln{3}\leq
$$
$$
\leq -2.1298\,\ln{x}/x+4.641/x+1.8205/x^2
$$

and for $x\geq 37$,
$$
\leq (-7.69+4.641+0.05)/x < 0.
$$

Thus $\ln{s(x)}$ grows faster than $h(x)$ for $x\geq 37$. We have $\ln{s(37)}=38.31$ and $h(37)=38.25$ whence (c).

For (d) dividing right-hand side by the left-hand side we have
$$
2^{-1}\,n^{2\,\log{n}+2\,\log{m}+1}\cdot m^{2\,\log{n}+2\,\log{m}+1}/n^{2\,\log{n}+1}\,m^{2\,\log{m}+1}=2^{-1}\,n^{2\,\log{m}}\,m^{2\,\log{n}}
$$

which is, clearly, $> 1$ if $n\geq 2, m\geq 2,$ as claimed.\\
\epr

Set $s(x):=2\,x^{2\,\log{x}+1}$.\\

\textbf{Lemma A8.} \textit{Let $H$ be a sporadic group or centrally isomorphic to $^{2}\bar{A}_3(9)$ and $F(n)$ the function associated to $H$ in ? and (A6)(e). Then $s(m)\,F(n)\leq s(n\,m)$ for $n\geq 2, m\geq 128.$}\\

The proof is, essentially, the same as that of (A6)(c). We have to use only that $s(m)\,s(n)\leq s(m\,n)$ if $m\geq 128, n\geq 6,$ $s(m)\,f(n)\leq s(m\,n)$ if $m\geq 128, n\geq 6,$ and then to check the claim for $y_{a_1, a_2, |H|}$ for each of the groups in question.\\

\textbf{Lemma A9.} \ \\

\textit{(a) $t!\,((n+2)!)^t\leq (n\,t+2)!$ for $t\geq 1, n > 12$, }\\

\textit{(b) $(n+2)!\,(m+2)!\leq (n+m+2)!$ for $n, m > 12.$}\\

\prf{} Consider $F(t, x)=\ln{\Gamma(t+1)}+t\,\ln{\Gamma(x+3)}-\ln{\Gamma(t\,x+3)}.$ We have 
$$
\partial F/\partial x=t\,(\ln{\Gamma})'(x+3)-t\,(\ln{\Gamma})'(t\,x+3)\leq 0
$$

if $x\geq 2, t\geq 1$ since $(\ln{\Gamma})'(s+1)$ increases as a function of $s$ for $s\geq ?$ by AO. We can assume $x\geq 13, t\geq 2.$ Thus for these $t$ and $x$
$$
F(t, x)\geq F(t, 13)=\ln{\Gamma}(t+1)+t\,\ln{\Gamma}(16)-\ln{\Gamma}(13\,t+3).
$$

We have $(\ln{\Gamma})''(t+1)=\sum_{i\geq 1}{(t+i)}^{-2}$ by AO. Therefore 
$$
d^2\,F(t, 13)/d\,t^2=\sum_{i\geq 1}{(t+i)}^{-2}-169\,\sum_{i\geq 1}{(13\,t+2+i)^{-2}} <
$$
$$
< \sum_{i\geq 1}{(t+i)^{-2}}-169\,\sum_{i\geq 1}{(13\,t+2+13\,i)^{-2}}-169\,(13\,t+2+i)^{-2}=
$$
$$
=\sum_{i\geq 1}{(t+i)^{-2}}-\sum_{i\geq 1}{(t+i+2/13)^{-2}}-(t+3/13)^{-2}<
$$
$$
< (t+1)^{-2}-(t+3/13)^{-2} < 0.
$$

Thus $d\,F(t, 13)/d\,t$ strictly decreases for $t\geq 2$. We have 
$$
d\,F/d\,t(2, 13)=(\ln{\Gamma})'(3)+\ln{\Gamma}(16)-13\,(\ln{\Gamma})'(29)=
$$
$$
=(-\gamma+ \sum^{2}_{i=1}{1/i})+27.9+13\,(\gamma-\sum^{28}_{i=1}{1/i})=
$$
$$
=12\,\gamma+1.5+27.9-13\cdot 3.927=-14,7 < 0
$$

(we used AO ?). Thus since $d\,F(t, 13)/d\,t$ decreases for $t\geq 2$ we have that $d\,F(t, 13)/d\,t < 0$ for $t\geq 2$. Thus $F(t, 13)$ decreases for $t\geq 2$. Since $F(2, 13)=\ln{(3!)}+2\,\ln{(15!)}-\ln{(28!)}=-11<0$ our claim (a) follows.

To prove (b) note that 
$$
(n+m+2)!/(n+2)!\,(m+2)!=(n+3)\,(n+4)\ldots(n+m+2)/(m+2)!=
$$
$$
=1\cdot 2^{-1}\cdot((n+3)/3)\,((n+4)/4)\ldots\geq 2^{-1}\,(13+3)/3=16/6 > 1
$$

whence (b).\\
\epr

\textbf{Lemma A10.} \textit{If $a^m\,m!/{(b\,m+2)!} < 1$ for some $m$ then it is $< 1$ for all larger $m$}.

\end{document}